\documentclass[a4paper,12pt]{amsart}
\usepackage[frenchb]{babel}
\usepackage[latin1]{inputenc}
\usepackage{pdfsync}
\usepackage[T1]{fontenc}
\usepackage{amsmath}
\usepackage{amssymb}
\usepackage{amsxtra}
\usepackage{amscd}
\usepackage{amsthm}
\usepackage{amsfonts}
\usepackage{eucal}
\usepackage[all]{xy}
\usepackage{graphicx}
\usepackage{comment}
\usepackage{epsfig}
\usepackage{psfrag}
\usepackage{mathrsfs}
\usepackage{amscd}
\usepackage{rotating}
\usepackage{lscape}
\usepackage{amsbsy}
\usepackage{verbatim}
\usepackage{moreverb}
\usepackage{url}

\newcommand{\nc}{\newcommand}
\nc{\renc}{\renewcommand}

\nc\restr[2]{{ 
  \left.\kern-\nulldelimiterspace    #1  
  \vphantom{\big|}  
  \right|_{#2}  
  }}

\newtheorem{thm}{Théorème}[section]
\newtheorem{prop}[thm]{Proposition}
\newtheorem{lem}[thm]{Lemme}
\newtheorem{sous-lem}[thm]{Sous-lemme}
\newtheorem{cor}[thm]{Corollaire}

\newtheorem{def-prop}[thm]{Définition-Proposition}
 \theoremstyle{definition}
\newtheorem{defi}[thm]{Définition}
\newtheorem{rem}[thm]{Remarque}

\newtheorem{observation}[thm]{Observation}
\numberwithin{equation}{section}

\renc{\sec}{\section}
\nc{\ssec}{\subsection}
\nc{\sssec}{\subsubsection}

\nc{\thmref}[1]{théorème~\ref{#1}}
\nc{\secref}[1]{paragraphe~\ref{#1}}
\nc{\lemref}[1]{lemme~\ref{#1}}
\nc{\defiref}[1]{définition~\ref{#1}}
\nc{\propref}[1]{proposition~\ref{#1}}
\nc{\corref}[1]{corollaire~\ref{#1}}
\nc{\constructionref}[1]{construction~\ref{#1}}
\nc{\conjref}[1]{conjecture~\ref{#1}}
\nc{\remref}[1]{remarque~\ref{#1}}
\nc{\questref}[1]{question~\ref{#1}}
\nc\Omegasour{\hbox{$\buildrel\smile\over{\vrule height 6pt depth 0pt width 0pt \smash \Omega}$}}

\nc{\on}{\operatorname}

\nc\wt{\widetilde}
\nc\wh{\widehat}
\nc\ol{\ov}
\nc{\oc}[1]{{\overset{\circ}{#1}}}
\nc{\ov}[1]{{\overline{#1}}}
\nc{\isor}[1]{{\xrightarrow[\raisebox{0.25 em}{\smash{\ensuremath{\sim}}}]{#1}}}
\nc{\modmod}{/ \! \! /}

\nc{\mc}{\mathcal}
\nc{\mf}{\mathfrak}
\nc{\mr}{\mathrm}
\nc{\mb}{\mathbb}
\nc{\mbf}{\mathbf}

\nc{\R}{{\mathbb R}}
\nc{\Z}{{\mathbb Z}}
\nc{\N}{{\mathbb N}}
\nc{\C}{{\mathbb C}}
\nc{\Q}{{\mathbb Q}}

\nc{\Fq}{{\mathbb F}_q}
\nc{\Fl}{{\mathbb F}_\ell}
\nc{\Fqbar}{\ol{{\mathbb F}_q}}
\nc{\Flbar}{\ol{{\mathbb F}_\ell}}
\nc{\Zl}{{\mathbb Z}_\ell}
\nc{\Zlbar}{\ol{{\mathbb Z}_\ell}}
\nc{\Ql}{{\mathbb Q}_\ell}
\nc{\Qlbar}{\ol{{\mathbb Q}_\ell}}
\nc{\hl}{\overset{\leftarrow}h{}}
\nc{\hr}{\overset{\rightarrow}h{}}
\nc{\Gr}{{\on{Gr}}}
\nc{\Hecke}{\on{Hecke}}
 \nc{\Hom}{\on{Hom}}
 \nc{\Coker}{\on{Coker}}
 \nc{\Ker}{\on{Ker}}
 \nc{\Lie}{\on{Lie}}
\nc{\Loc}{\on{Loc}}
\nc{\Pic}{\on{Pic}}
\nc{\Bun}{\on{Bun}}
\nc{\IC}{\on{IC}}
\nc{\Aut}{\on{Aut}}
\nc{\Perv}{\on{Perv}}
\nc{\pos}{{\on{pos}}}
\nc{\Sym}{\on{Sym}}

\nc{\ta} {{}^\tau}
\nc {\tu}[1]{{}^{\tau^{#1}}\!}

\nc{\Id}{\on{Id}}
\nc{\Fil}{\on{Fil}}
\nc{\pr}{\on{pr}}
\nc{\Res}{\on{Res}}
\nc{\cusp}{\on{cusp}}
\nc{\Frob}{\on{Frob}}
\nc{\diag}{\Delta}
\nc{\gr}{\on{gr}}
\nc{\Inj}{\on{Inj}}
\nc{\Bl}{\on{Bl}}
\nc{\dem}{\noindent {\bf Démonstration. }}
\nc{\cqfd}{{\ }\hfill $\square$ \vskip 1mm}
\nc{\s}[1]{\langle #1 \rangle}
\nc{\Cht}{\on{Cht}}
\nc{\isom}{\overset {\thicksim}{\to}}
\nc{\sm}{\smallsetminus}

\begin{document}

\title[Chtoucas et paramétrisation de Langlands]{
Introduction aux chtoucas pour les groupes réductifs et à la  paramétrisation de  Langlands globale}

\author{Vincent Lafforgue}
\thanks{L'auteur fait partie de l'ANR-13-BS01-0001-01}
\address{Vincent Lafforgue: CNRS et Institut Fourier, UMR 5582, Université Grenoble Alpes, 
 100 rue des Maths, 38610 Gières, France.}

\date{\today}
\maketitle

      \section*{Introduction}
     
     Cet article introductif a pour but d'expliquer le résultat de 
      \cite{coh} et 
     de donner {\it toutes les idées}
     de la démonstration. De \cite{coh} il  reprend essentiellement l'introduction (en la rallongeant un peu) et les paragraphes 12.1 et 12.2 (en les résumant).      
      
         On montre le  sens ``automorphe vers Galois'' de la correspondance de Langlands  globale \cite{langlands67} pour tout groupe réductif $G$ sur un corps de fonctions.              
      De plus  on construit une décomposition {\it canonique} de l'espace des formes automorphes cuspidales, indexée par des  paramètres  de Langlands globaux.
  On n'obtient  pas de résultat nouveau dans   le cas où $G=GL_{r}$ puisque tout était déjà connu par Drinfeld \cite{drinfeld78,Dr1,drinfeld-proof-peterson,drinfeld-compact} pour $r=2$  et Laurent Lafforgue \cite{laurent-inventiones} pour $r$ arbitraire. 

          La méthode  est totalement indépendante de la formule des traces d'Arthur-Selberg. Elle  utilise les deux   ingrédients suivants : 
        \begin{itemize}
     \item les champs classifiants de chtoucas, introduits par   Drinfeld pour $GL_{r}$  
     \cite{drinfeld78,Dr1} et généralisés à tous les groupes  réductifs par Varshavsky
     \cite{var} 
   \item l'équivalence de  Satake géométrique de   Lusztig, Drinfeld, Ginzburg, et  Mirkovic--Vilonen 
   \cite{lusztig-satake,ginzburg,hitchin,mv}. 
     \end{itemize}

\noindent {\bf Remerciements. }
Ce travail est issu d'un travail en cours en commun avec Jean-Benoît Bost, 
et il n'existerait pas non plus sans de nombreuses discussions avec 
Alain Genestier. Je les remercie aussi beaucoup pour toute l'aide qu'ils m'ont apportée. 

Cet article repose doublement sur des travaux de Vladimir Drinfeld,  sur les chtoucas et  sur le programme de 
 Langlands géométrique, et je lui exprime ma grande reconnaissance.  
 
Je remercie Vladimir Drinfeld, Dennis Gaitsgory,  Jochen Heinloth, Laurent Lafforgue,  Sergey Lysenko et Yakov Varshavsky pour   leurs explications. Je remercie 
  Gebhard B\"ockle,  Michael Harris, Chandrashekhar Khare et Jack Thorne pour avoir trouvé l'énoncé du \lemref{prop-harris} (qui rend la rédaction plus limpide et est de plus utilisé dans  leur article \cite{boeckle-harris...}). 
Je remercie aussi 
    Ga\"etan Chenevier, Pierre Deligne,  Ed Frenkel,  Nicholas Katz,  Erez Lapid, Gérard Laumon,   Colette Moeglin, Sophie Morel, Ngô Bao Châu, Ngo Dac Tuan, Jean-Pierre Serre,   Jean-Loup Waldspurger  et Cong Xue.

 Je tiens aussi à remercier vivement le CNRS. Le programme de Langlands est  un sujet très différent de ma spécialité d'origine et je n'aurais  pas pu  m'y consacrer sans la grande liberté laissée aux chercheurs pour mener à bien leurs travaux. 
 
\noindent {\bf Plan. } Dans le paragraphe \ref{para-enonces} on donne l'énoncé principal (pour $G$ déployé), l'idée de la démonstration et deux énoncés intermédiaires qui structurent la démonstration. Le paragraphe \ref{para-commentaires}, qui 
résume les paragraphes 12.1 et 12.2   de \cite{coh}, discute le cas non déployé, quelques compléments et un certain nombre de problèmes ouverts. Le paragraphe \ref{defi-chtou-intro} introduit les champs de chtoucas. Le paragraphe \ref{section-esquisse-abc} montre que la partie ``Hecke-finie'' de leur cohomologie vérifie certaines propriétés qui avaient été énoncées dans la \propref{prop-a-b-c}  du paragraphe \ref{para-enonces}.  A l'aide de ces propriétés on montre dans le paragraphe \ref{intro-idee-heurist} le théorème principal, qui avait été énoncé dans le paragraphe \ref{para-enonces}. Dans les paragraphes \ref{subsection-crea-annihil-intro} à \ref{subsection-intro-decomp} on montre des propriétés de la cohomologie des champs de chtoucas, qui avaient été admises dans le paragraphe \ref{section-esquisse-abc} et à la fin du paragraphe \ref{intro-idee-heurist}.   
Le paragraphe \ref{subsection-link-langl-geom}
explique le lien avec le programme de Langlands géométrique. 
Enfin le  paragraphe \ref{intro-previous-works}
discute les rapports avec les travaux antérieurs.   

Suivant le temps dont il dispose, le lecteur peut  se limiter au paragraphe \ref{para-enonces}, aux paragraphes \ref{para-enonces} et 
\ref{para-commentaires}, ou bien aux paragraphes \ref{para-enonces} à \ref{intro-idee-heurist}. Les paragraphes \ref{para-enonces}, \ref{para-commentaires} et  \ref{intro-idee-heurist} ne nécessitent pas de connaissances en géométrie algébrique.

\section{Enoncé et idées principales}\label{para-enonces}

   \subsection{Préliminaires}        Soit  $\Fq$ un corps fini.   Soit  $X$ une courbe projective lisse géométriquement irréductible  sur $\Fq$ et  $F$ son corps de fonctions. 
      Soit  $G$ un groupe réductif connexe sur $F$. Soit   $\ell$ un nombre premier ne divisant pas $q$.

    Pour énoncer le théorème principal on suppose que $G$ est déployé (le cas non déployé sera expliqué dans le paragraphe \ref{cas-non-deploye}).  
        On note  $\wh G$ le groupe dual de  Langlands de  $G$,  considéré comme un groupe déployé sur  $\Ql$. 
        Ses racines et ses poids sont les coracines et les copoids de $G$, et vice-versa (voir \cite{borel-corvallis}  pour plus de  détails).

Pour des raisons topologiques on doit travailler avec des extensions finies de $\Ql$ plutôt que $\Qlbar$.  Soit $E$ une extension finie de $\Ql$ contenant une racine carrée de $q$ et $\mc O_{E}$ son anneau d'entiers.

        Soit  $v$ une place de $X$. On note  $\mc O_{v}$ l'anneau local complété en  $v$ et  $F_{v}$ son corps de fractions. 
   On a l'isomorphisme de   Satake    $[V]\mapsto h_{V,v}$ de l'anneau des représentations de  $\wh G$ (à coefficients dans $E$) vers l'algèbre de  Hecke  
      $C_{c}(G(\mc O_{v})\backslash G(F_{v})/G(\mc O_{v}),E)$ (voir  
      \cite{satake,cartier-satake,gross}). En fait les $h_{V,v}$  pour $V$ irréductible forment une base sur $\mc O_{E}$ de $C_{c}(G(\mc O_{v})\backslash G(F_{v})/G(\mc O_{v}),\mc O_{E})$. 
      On note   $\mb A=\prod_{v\in |X|} ' F_{v}$ l'anneau des adèles  de $F$ et $\mb O=
      \prod_{v\in |X|} \mc O_{v}$. 
             Soit   $N$ un sous-schéma fini de  $X$.   
   On note  
   $\mc O_{N}$ l'anneau des fonctions sur $N$, et 
   \begin{gather}\label{def-K-N}K_{N}=\on{Ker}(G(\mb O)\to G(\mc O_{N}))\end{gather} le sous-groupe compact ouvert de $G(\mb A)$ associé au niveau  $N$. On fixe un réseau   $\Xi\subset Z(F)\backslash Z(\mb A)$ (où $Z$ est le centre de $G$).  
   Une fonction $f\in C_{c} (G(F)\backslash G(\mb A)/K_N \Xi,E)$ 
   est dite cuspidale si pour tout 
    parabolique $P\subsetneq G$, de  Levi $M$ et de radical unipotent $U$,    le terme constant $f_{P}: g\mapsto \int_{U(F)\backslash U(\mb A)}f(ug)$ est  nul comme fonction sur $U(\mb A)M(F)\backslash G(\mb A)/K_{N}\Xi$.  
 On rappelle que $C_{c}^{\rm{cusp}}(G(F)\backslash G(\mb A)/K_N \Xi,E)$ est un $E$-espace vectoriel  de dimension finie. Il est muni d'une structure de module sur l'algèbre de Hecke $C_{c}(K_{N}\backslash G(\mb A)/K_{N},E)$ : on convient  que la fonction caractéristique de $K_{N}$ est une unité et agit par l'identité et pour $f\in   C_{c}(K_{N}\backslash G(\mb A)/K_{N},E)$ on note $T(f)\in \on{End}(C_{c}^{\rm{cusp}}(G(F)\backslash G(\mb A)/K_N \Xi,E))$ l'opérateur de Hecke correspondant. 
    
 \subsection{Enoncé du théorème principal} 
   On construira  les opérateurs suivants, dits ``d'excursion''. 
    Soient $I$ un ensemble fini, $f$ une fonction sur 
$\wh G \backslash (\wh G)^{I}/\wh G    $ 
(quotient grossier de $(\wh G)^{I}$ par les translations à gauche et à droite par $\wh G $ diagonal), et 
$(\gamma_{i})_{i\in I}\in \on{Gal}(\ov F/F)^{I}$. 
On construira 
l'opérateur  d'excursion $$S_{I,f,(\gamma_{i})_{i\in I}}\in \mr{End}_{C_{c}(K_{N}\backslash G(\mb A)/K_{N},E)}( C_{c}^{\mr{cusp}}(G(F)\backslash G(\mb A)/K_{N}\Xi,E)). 
    $$ On montrera que ces opérateurs  engendrent une sous-algèbre commutative $\mc B$. 
    
    On ne sait pas si $\mc B$ est réduite mais 
par décomposition spectrale  on obtient néanmoins  une décomposition canonique 
 \begin{gather}\label{intro1-nu-dec-canonique}
 C_{c}^{\mr{cusp}}(G(F)\backslash G(\mb A)/K_{N}\Xi,\Qlbar)=\bigoplus_{\nu}
 \mf H_{\nu} \end{gather}
    où la somme directe dans le membre de droite est indexée par des  caractères $\nu$ de $\mc B$, et où 
   $\mf H_{\nu}$ est l'espace propre généralisé (ou ``espace caractéristique'') associé à $\nu$. On montrera ensuite qu'à tout caractère $\nu$ de $\mc B$ correspond un {\it unique} paramètre  de Langlands $\sigma$ (au sens du théorème suivant), caractérisé par \eqref{relation-fonda} ci-dessous. 
    En posant $ \mf H_{\sigma} =\mf H_{\nu}$, on en déduira le théorème suivant.

 \begin{thm}  \label{intro-thm-ppal}  On possède  une décomposition canonique de  
   $C_{c}(K_{N}\backslash G(\mb A)/K_{N},\Qlbar)$-modules
 \begin{gather}\label{intro1-dec-canonique}
 C_{c}^{\mr{cusp}}(G(F)\backslash G(\mb A)/K_{N}\Xi,\Qlbar)=\bigoplus_{\sigma}
 \mf H_{\sigma},\end{gather}
 où la somme directe dans le membre de droite est indexée par des paramètres de Langlands globaux, c'est-à-dire des classes de  $\wh G(\Qlbar)$-conjugaison de  morphismes 
       $\sigma:\on{Gal}(\ov F/F)\to \wh G(\Qlbar)$ 
       définis sur une extension finie de  $\Ql$, continus,  semi-simples et non ramifiés en dehors de  $N$.
       
     Cette décomposition  est caractérisée par la propriété suivante : 
       $ \mf H_{\sigma}$ est égal à l'espace propre généralisé $\mf H_{\nu}$ 
        associé au caractère $\nu$ de $\mc B$ défini par 
       \begin{gather}\label{relation-fonda}\nu(S_{I,f,(\gamma_{i})_{i\in I}})=f((\sigma(\gamma_{i}))_{i\in I}. \end{gather}
       
      Elle est compatible avec l'isomorphisme de  Satake en toute place  $v$  de $X\sm N$, c'est-à-dire que pour toute représentation irréductible $V$ de $\wh G$, 
       $T(h_{V,v})$ agit sur   $\mf H_{\sigma}$
       par multiplication par le scalaire  $\chi_{V}(\sigma(\Frob_{v}))$, où $\chi_{V}$ est le caractère de $V$ et $\Frob_{v}$ est un relèvement arbitraire d'un élément de  Frobenius   en $v$. 
     Elle   est aussi compatible avec la limite sur  $N$.   
    \end{thm}

La compatibilité avec l'isomorphisme de  Satake en les places     de $X\sm N$ montre que ce théorème réalise la ``correspondance'' de Langlands globale dans le sens ``automorphe vers Galois''. 
En fait, sauf dans le cas de $GL_{r}$ \cite{laurent-inventiones},   les conjectures  de Langlands consistent  plutôt en 
\begin{itemize}
\item une paramétrisation, obtenue dans le théorème ci-dessus, 
\item des formules de multiplicités d'Arthur pour les  $\mf H_{\sigma}$, que nous ne savons pas calculer avec les méthodes de cet article. \end{itemize}

\subsection{Idées principales de la démonstration} 
Pour construire les opérateurs d'excursion et montrer ce théorème la stratégie sera la suivante. Les champs de chtoucas, qui jouent un rôle analogue aux variétés de Shimura sur les corps de nombres, existent dans une  généralité beaucoup plus grande. En effet, alors que les variétés de Shimura  sont définies sur un ouvert du spectre d'un anneau d'entiers d'un  corps de nombres et sont associées à {\it un copoids  minuscule} du groupe dual, on possède
pour tout ensemble fini $I$, pour tout  niveau $N$ et pour toute représentation irréductible $W$  de $(\wh G)^{I}$  un champ de chtoucas $\on{Cht}_{N,I,W}$ qui est défini  sur $(X\sm N)^{I}$. 

On construira alors un $E$-espace vectoriel $H_{I,W}$
(où la lettre $N$ est omise pour raccourcir les formules) 
 comme un certain   {\it sous-espace}   de la cohomologie d'intersection 
 de la fibre de $\on{Cht}_{N,I,W}$   au-dessus d'un point géométrique générique de 
$(X\sm N)^{I}$. Par cohomologie d'intersection on entend ici  la cohomologie d'intersection 
 à support compact,  à coefficients dans $E$ et  en degré $0$ (pour la normalisation perverse). Grâce aux ``morphismes de Frobenius partiels'' introduits par Drinfeld,   on munira $H_{I,W}$  d'une action de $\on{Gal}(\ov F/F)^{I}$.

 \begin{rem} \label{rem-HIW-dim-finie} Dans cet article on définira ce {\it sous-espace}  $H_{I,W}$ par une condition technique (de finitude sous l'action des opérateurs de Hecke). 
 En fait Cong Xue a montré dans \cite{these-cong} 
 qu'il admet une définition équivalente comme   partie ``cuspidale'' de la cohomologie d'intersection et qu'il est de dimension finie. \end{rem}
 
L'équivalence de Satake géométrique permet de raffiner 
cette construction de $H_{I,W}$ (défini ci-dessus pour toute classe d'isomorphisme de représentation irréductible $W$ de $(\wh G)^{I}$) en 
celle de {\it foncteurs}  $W\mapsto  
    H_{I,W}$, munis de la {\it donnée} des isomorphismes 
    \eqref{isom-chi-zeta-b} ci-dessous. Autrement dit   $H_{I,W}$ est {\it canonique} et se comporte naturellement quand on fait varier $W$ ou $I$. 
    
   \begin{prop}\label{prop-a-b-c}    Les $H_{I,W}$ vérifient   les 
  propriétés   suivantes : 
        \begin{itemize}
    \item[] {\bf a)} pour tout ensemble fini   $I$,      $$W\mapsto  
    H_{I,W},  \ \ u\mapsto \mc H(u)$$  est un foncteur  $E$-linéaire  de la  catégorie des  représentations $E$-linéaires de dimension finie de  $(\wh G)^{I}$ vers la  catégorie des  limites inductives de représentations $E$-linéaires continues de dimension finie de     $\on{Gal}(\ov F/F)^{I}$,    
              \item[] {\bf b)} pour toute application   $\zeta: I\to J$, 
 on possède  un isomorphisme 
      \begin{gather}\label{isom-chi-zeta-b}
           \chi_{\zeta}: H_{I,W}\isom 
 H_{J,W^{\zeta}},\end{gather} 
 qui est 
 \begin{itemize}
 \item  fonctoriel en   $W$, où  $W$ est une  représentation de $(\wh G)^{I}$ et  $W^{\zeta}$ désigne la   représentation de $(\wh G)^{J}$ sur  $W$ obtenue en composant avec le  morphisme  diagonal $$  (\wh G)^{J}\to (\wh G)^{I}, (g_{j})_{j\in J}\mapsto (g_{\zeta(i)})_{i\in I} $$ 
 \item $\on{Gal}(\ov F/F)^{J}$-équivariant, où $\on{Gal}(\ov F/F)^{J}$ agit sur le membre de gauche par le   morphisme diagonal  
 \begin{gather}\label{morp-diag-Gal-intro0}\on{Gal}(\ov F/F)^{J}\to \on{Gal}(\ov F/F)^{I},  \ (\gamma_{j})_{j\in J}\mapsto (\gamma_{\zeta(i)})_{i\in I}, 
\end{gather}
   \item   et compatible avec la  composition, c'est-à-dire   que pour 
 $I\xrightarrow{\zeta} J\xrightarrow{\eta} K$ on a 
 $\chi_{\eta\circ \zeta}=\chi_{\eta}\circ\chi_{\zeta}$,
  \end{itemize}
     \item[] {\bf c)} pour $I=\emptyset$ et  $W=\mbf  1$, on a un isomorphisme     \begin{gather}\label{c-de-la-prop}
       H_{\emptyset,\mbf  1}=C_{c}^{\mr{cusp}}(G(F)\backslash G(\mb A)/K_{N}\Xi,E). \end{gather}
    \end{itemize}
    
  Par ailleurs les $H_{I,W}$ sont des modules sur 
  $C_{c}(K_{N}\backslash G(\mb A)/K_{N},E)$, de fa\c con compatible  avec 
les propriétés a), b), c) ci-dessus. 
  \end{prop}

\begin{rem}
Pour tout ensemble fini $J$, en appliquant b) à l'application évidente $\zeta:\emptyset \to J$, on obtient un isomorphisme 
$\chi_{\zeta}: H_{\emptyset,\mbf  1}\isom H_{J,\mbf  1}$ 
(où $\mbf 1$ désigne la représentation triviale de  $(\wh G)^{J}$) 
et donc  l'action de $\on{Gal}(\ov F/F)^{J}$ sur $H_{J,\mbf  1}$ est triviale. 
\end{rem}

Grâce à \eqref{c-de-la-prop} la décomposition \eqref{intro1-dec-canonique} que l'on cherche est équivalente à une décomposition 
\begin{gather}\label{dec-emptyset-1-sigma}H_{\emptyset,\mbf  1}=\bigoplus_{\sigma}
 \mf H_{\sigma}\end{gather}
 (où l'on suppose, quitte à augmenter $E$, que les $\sigma$ et $ \mf H_{\sigma}$ sont définis sur $E$).

La définition suivante, où l'on construit  les opérateurs d'excursion,  sera répétée dans le paragraphe \ref{intro-idee-heurist}. 
Le lecteur peut consulter dès à présent, s'il le souhaite, le paragraphe \ref{descr-conj-HIW} ci-dessous pour une description heuristique des $H_{I,W}$, qui éclaire {\it a posteriori}
la définition des opérateurs d'excursion. 

 On a besoin de considérer un ensemble à un élément et on le note $\{0\}$. 
Pour tout ensemble fini $I$ on note $\zeta_{I}:I\to \{0\}$ l'application évidente,  de sorte que $W^{\zeta_{I}}$ est simplement $W$ muni de l'action diagonale de $\wh G$. 

Grâce à la remarque ci-dessus (appliquée à $J=\{0\}$, de sorte que $\zeta:\emptyset \to \{0\}$ est noté maintenant $\zeta_{\emptyset}$), et à   \eqref{c-de-la-prop}, on a 
 \begin{gather}\label{egalite-vide-0-cusp}H_{\{0\},\mbf 1}\isor{\chi_{\zeta_{\emptyset}}^{-1}} H_{ \emptyset ,\mbf 1}= C_{c}^{\mr{cusp}}(G(F)\backslash G(\mb A)/K_{N}\Xi,E).\end{gather}    
On va définir maintenant les opérateurs d'excursion comme des endomorphismes de \eqref{egalite-vide-0-cusp}  en utilisant  $H_{\{0\},\mbf 1}$. 

\begin{defi} \label{defi-constr-excursion-intro} Pour toute fonction $f\in\mc O(\wh G\backslash (\wh G)^{I}/\wh G)$ on peut trouver 
    une représentation $W$ de $(\wh G)^{I}$ et $x\in W$ et $\xi\in W^{*}$ invariants par l'action diagonale de $\wh G$ tels que 
    \begin{gather}\label{f-W-x-xi-intro}    f ((g_{i})_{i\in I})=\s{\xi, (g_{i})_{i\in I}\cdot x}.\end{gather}
      Alors l'endomorphisme 
    $S_{I,f,(\gamma_{i})_{i\in I}}$ de  \eqref{egalite-vide-0-cusp} 
   est défini comme la composée 
  \begin{gather}\label{excursion-def-intro1}
  H_{\{0\},\mbf  1}\xrightarrow{\mc H(x)}
 H_{\{0\},W^{\zeta_{I}}}\isor{\chi_{\zeta_{I}}^{-1}} 
  H_{I,W}
  \xrightarrow{(\gamma_{i})_{i\in I}}
  H_{I,W} \isor{\chi_{\zeta_{I}}} H_{\{0\},W^{\zeta_{I}}}  
  \xrightarrow{\mc H(\xi)} 
  H_{\{0\},\mbf  1} 
    \end{gather}
    où $x:\mbf 1\to W^{\zeta_{I}}$ et $\xi: W^{\zeta_{I}}\to \mbf 1$ sont considérés ici comme des morphismes de représentations de $\wh G$. 
\end{defi}

On montrera   dans le paragraphe \ref{intro-idee-heurist}, à l'aide des propriétés a) et b) de la proposition  \ref{prop-a-b-c}, que 
$S_{I,f,(\gamma_{i})_{i\in I}}$   ne dépend pas du choix de $W,x,\xi$ vérifiant 
\eqref{f-W-x-xi-intro} et qu'il est donc bien défini. 
On montrera dans le \lemref{prop-harris} 
(trouvé par B\"ockle,    Harris,  Khare et  Thorne) 
        que    $S_{I,f,(\gamma_{i})_{i\in I}}$ dépend seulement de  l'image de 
 $(\gamma_{i})_{i\in I}$ dans   $\pi_{1}(X\sm N, \ov\eta)^{I}$ (où $   \ov\eta=\on{Spec} \ov F$).      
On montrera que ces opérateurs d'excursion 
commutent entre eux, et en notant $\mc B$ la sous-algèbre commutative de 
$\on{End}(C_{c}^{\mr{cusp}}(G(F)\backslash G(\mb A)/K_{N}\Xi,E))$ qu'ils engendrent, qu'ils 
 vérifient les propriétés suivantes (qui sont bien celles que l'on attend car elles sont tautologiquement  vérifiées par le membre de droite de \eqref{relation-fonda}). 

 \begin{prop} \label{prop-SIf-i-ii-iii} Les  opérateurs d'excursion  
   $S_{I,f,(\gamma_{i})_{i\in I}}$ vérifient les propriétés suivantes:   
  \begin{itemize}
  \item [] (i) pour tout   $I$ et 
 $(\gamma_{i})_{i\in I}\in  \pi_{1}(X\sm N, \ov\eta)^{I}$, 
  $$f\mapsto 
  S_{I,f,(\gamma_{i})_{i\in I}}$$ est un  morphisme 
  d'algèbres  commutatives  $\mc O(\wh G\backslash (\wh G)^{I}/\wh G)\to \mc B$, 
  \item [] (ii) pour toute application 
  $\zeta:I\to J$, toute fonction  $f\in \mc O(\wh G\backslash (\wh G)^{I}/\wh G)$ et  tout 
  $(\gamma_{j})_{j\in J}\in  \pi_{1}(X\sm N, \ov\eta)^{J}$, on a 
  $$S_{J,f^{\zeta},(\gamma_{j})_{j\in J}}=S_{I,f,(\gamma_{\zeta(i)})_{i\in I}}$$
   où $f^{\zeta}\in \mc O(\wh G\backslash (\wh G)^{J}/\wh G)$ est définie par    $$f^{\zeta}((g_{j})_{j\in J})=f((g_{\zeta(i)})_{i\in I}),$$
   \item [] (iii) 
  pour tout   $f\in \mc O(\wh G\backslash (\wh G)^{I}/\wh G)$
  et  $(\gamma_{i})_{i\in I},(\gamma'_{i})_{i\in I},(\gamma''_{i})_{i\in I}\in  \pi_{1}(X\sm N, \ov\eta)^{I}$ on a     $$S_{I\cup I\cup I,\wt f,(\gamma_{i})_{i\in I}\times (\gamma'_{i})_{i\in I}\times (\gamma''_{i})_{i\in I}}=
  S_{I,f,(\gamma_{i}(\gamma'_{i})^{-1}\gamma''_{i})_{i\in I}}$$
   où  $I\cup I\cup I$ est une réunion disjointe et 
   $\wt f\in \mc O(\wh G\backslash (\wh G)^{I\cup I\cup I}/\wh G)$ est définie  par  
   $$\wt f((g_{i})_{i\in I}\times (g'_{i})_{i\in I}\times (g''_{i})_{i\in I})=f((g_{i}(g'_{i})^{-1}g''_{i})_{i\in I}).$$
      \item [] (iv) pour tout $I$ et tout $f$, le morphisme 
      \begin{gather}\label{mor-excursion-f}\pi_{1}(X\sm N, \ov\eta)^{I}\to \mc B, \ \ (\gamma_{i})_{i\in I}\mapsto S_{I,f,(\gamma_{i})_{i\in I}}\end{gather} est continu, avec $\mc B$ munie de la topologie $E$-adique, 
      \item [] (v) pour toute place $v$ de $X\sm N$, pour toute représentation irréductible $V$ de $\wh G$, l'opérateur de Hecke $$T(h_{V,v})\in \on{End}(C_{c}^{\mr{cusp}}(G(F)\backslash G(\mb A)/K_{N}\Xi,E))$$ est égal à l'opérateur d'excursion $ S_{\{1,2\}, f,(\Frob_{v},1)}$ où $f\in \mc O(\wh G\backslash (\wh G)^{2}/\wh G)$ est donnée par $f(g_{1},g_{2})=\chi_{V}(g_{1}g_{2}^{-1})$, et $\Frob_{v}$ est un élément de  Frobenius   en $v$. 
            \end{itemize}
         \end{prop}

En fait les propriétés (i), (ii), (iii) et (iv) résulteront formellement de la proposition  \ref{prop-a-b-c} et (v) s'obtiendra par un argument géométrique
(le calcul de la composée de deux correspondances cohomologiques entre champs de chtoucas). 
    
 \begin{rem} \label{rem-gamma-en-plus}
Pour tout $\gamma\in  \pi_{1}(X\sm N, \ov\eta) $, on a 
  $S_{I,f,(\gamma_{i})_{i\in I}}=S_{I,f,(\gamma_{i}\gamma)_{i\in I}}$.   On  le vérifie très facilement à l'aide de la définition des opérateurs d'excursion, ou bien  on le déduit de la proposition précédente. 
  De même on a  $S_{I,f,(\gamma_{i})_{i\in I}}=S_{I,f,(\gamma\gamma_{i})_{i\in I}}$. 
\end{rem}

    Dans le paragraphe \ref{intro-idee-heurist} on  déduira assez facilement de la proposition précédente qu'à tout  caractère $\nu$ de $\mc B$   correspond un paramètre de Langlands
    $\sigma:\pi_{1}(X\sm N, \ov\eta) \to \wh G(\Qlbar)$ vérifiant \eqref{relation-fonda}, unique à conjugaison près par 
    $\wh G(\Qlbar)$. En effet  la connaissance de $\nu(S_{I,f,(\gamma_{i})_{i\in I}})$ (qui doit être égal à $f((\sigma(\gamma_{i}))_{i\in I}$)   pour toute fonction $f$ détermine l'image du $I$-uplet 
    $(\sigma(\gamma_{i}))_{i\in I}\in (\wh G(\Qlbar))^{I}$ comme point défini sur $\Qlbar$ du quotient grossier $\wh G\backslash (\wh G)^I/\wh G$. Or en prenant $I=\{0,..,n\}$ on voit que 
       $$  (\wh G)^{n}\modmod \wh G\isom \wh G\backslash (\wh G)^{\{0,...,n\}}/\wh G ,  (g_{1},...,g_{n})\mapsto 
   (1,g_{1},...,g_{n})$$ est  un  isomorphisme, en notant $ (\wh G)^{n}\modmod \wh G$ le quotient grossier de $ (\wh G)^{n}$ par conjugaison diagonale. Donc pour tout entier $n$ et pour tout $n$-uplet 
   $(\gamma_{1},...,\gamma_{n})\in  \pi_{1}(X\sm N, \ov\eta)^{n}$, la connaissance de $\nu$ détermine l'image de $(\sigma(\gamma_{1}), ..., \sigma(\gamma_{n}))$ comme point défini sur $\Qlbar$ du quotient grossier    $ (\wh G)^{n}\modmod \wh G$.  
    Grâce à des résultats de \cite{richardson} fondés sur la théorie géométrique des invariants, cela signifie que $\nu$ détermine 
    $(\sigma(\gamma_{1}), ..., \sigma(\gamma_{n}))\in (\wh G(\Qlbar))^{n}$ à semi-simplification et conjugaison diagonale près. Comme on demande que $\sigma$ soit semi-simple, il est clair 
    (en choisissant $n$ et $(\gamma_{1},...,\gamma_{n})$ pour que 
    le sous-groupe engendré par $\sigma(\gamma_{1}), ..., \sigma(\gamma_{n})$ soit Zariski dense dans l'image de $\sigma$) 
    que ces données déterminent $\sigma$ à conjugaison près. Inversement les relations (i), (ii), (iii) et (iv) satisfaites par les opérateurs d'excursion permettront, dans la \propref{intro-Xi-n},  de montrer l'existence de $\sigma$ vérifiant \eqref{relation-fonda} et la propriété
   (v) assurera  la compatibilité avec l'isomorphisme de Satake en les places de $X\sm N$. 
    
     \subsection{Une remarque heuristique}  \label{descr-conj-HIW}
    Ce paragraphe propose une description conjecturale
  des $H_{I,W}$, qui justifie {\it a posteriori} la définition des opérateurs d'excursion donnée dans \eqref{excursion-def-intro1}. Bien sûr cette description conjecturale des $H_{I,W}$ n'intervient pas dans les raisonnements
  et   ce paragraphe  ne sera utilisé  nulle part dans le reste de l'article.

  On conjecture qu'il existe un ensemble fini $\Sigma$ (dépendant de $N$) 
 de paramètres de Langlands semi-simples (bien déterminés à conjugaison près), et que, quitte à augmenter $E$,    on possède pour tout $\sigma\in \Sigma$
     une  représentation  $E$-linéaire   $A_{\sigma}$ 
      du centralisateur   $S_{\sigma}$  de l'image de $\sigma$ dans $\wh G$  (triviale sur $Z(\wh G)$), de telle sorte que 
  pour tout    $I$ et  $W$  \begin{gather}\label{intro-dec-I-W}
 H_{I,W}\overset{?}{=}\bigoplus_{\sigma\in \Sigma}
\Big(A_{\sigma} \otimes _{E}W_{\sigma^{I}}\Big)^{S_{\sigma}}, \end{gather}
 où 
$W_{\sigma^{I}}$ désigne la   représentation de  
$\pi_{1}(X\sm N, \ov\eta)^{I}$     obtenue en composant la   représentation $W$ avec le  morphisme   
$\sigma^{I}: \pi_{1}(X\sm N, \ov\eta)^{I} \to (\wh G(E))^{I}$. 
De plus $A_{\sigma}$ doit être un module sur 
$C_{c}(K_{N}\backslash G(\mb A)/K_{N},E)$, 
et \eqref{intro-dec-I-W} doit être un isomorphisme de 
$C_{c}(K_{N}\backslash G(\mb A)/K_{N},E)$-modules. 
Dans le cas particulier où $I= \emptyset$  et $W=1$, 
\eqref{intro-dec-I-W} doit être la décomposition  
\eqref{dec-emptyset-1-sigma}  et   on doit avoir 
$ \mf H_{\sigma}=(A_{\sigma})^{S_{\sigma}}$. 
  
   Ces conjectures sont bien connues des experts, par extrapolation des 
conjectures  de   \cite{kottwitz3}    sur les multiplicités dans la cohomologie des variétés de  Shimura et grâce au résultat de Cong Xue mentionné dans la \remref{rem-HIW-dim-finie}.  
Dans le cas de $GL_{r}$ on s'attend à ce que 
$\Sigma$ soit l'ensemble des représentations irréductibles de dimension $r$ de $\pi_{1}(X\sm N, \ov\eta)$ et que pour tout 
$\sigma\in \Sigma$, 
$S_{\sigma}=\mathbb G_{m}=Z(\wh G)$  et $A_{\sigma}=(\pi_{\sigma})^{K_{N}}$ où $\pi_{\sigma}$ est la représentation automorphe cuspidale correspondant à $\sigma$ (voir \cite{laurent-inventiones} et  la conjecture 2.35 de \cite{var}). En général si $\sigma$ est associé à un paramètre d'Arthur elliptique  $\psi$ (comme dans le paragraphe \ref{para-param-Arthur} ci-dessous), $A_{\sigma}$ devrait être induit d'une  représentation de dimension finie du sous-groupe de  $S_{\sigma}$ engendré par le centralisateur de $\psi$ et par le sous-groupe diagonal  $\mathbb{G}_{m}\subset SL_{2}$ (parce que nous considérons seulement la cohomologie en degré $0$).

On conjecture de plus que   \eqref{intro-dec-I-W} est fonctoriel en   $W$ et que pour toute application $\zeta:I\to J$ il entrelace   $\chi_{\zeta}$ avec  
  \begin{gather}\nonumber \Id:  \bigoplus_{\sigma} \Big(A_{\sigma} \otimes _{E}W_{\sigma^{I}}\Big)^{S_{\sigma}}\to
 \bigoplus_{\sigma} 
 \Big(A_{\sigma} \otimes _{E}(W^{\zeta})_{\sigma^{J}}\Big)^{S_{\sigma}}
   \end{gather}
   (comme $W_{\sigma^{I}}$ et $(W^{\zeta})_{\sigma^{J}}$ sont tous les deux égaux à $W$ comme $E$-espaces vectoriels, $\Id$ a bien un sens et il est $\pi_{1}(X\sm N, \ov\eta)^{J}$-équivariant). 
 Sous ces hypothèses,   la composée \eqref{excursion-def-intro1} (qui   définit  $S_{I,f,(\gamma_{i})_{i\in I}}$) 
 agit sur   
    $ \mf H_{\sigma}= (A_{\sigma})^{S_{\sigma}} \subset  H_{\{0\},\mbf 1}$ par   la composée  
      \begin{gather} \nonumber (A_{\sigma})^{S_{\sigma}}
 \xrightarrow{\Id_{A_{\sigma}}\otimes x} 
\Big(A_{\sigma} \otimes _{E}W_{\sigma^{I}}\Big)^{S_{\sigma}}
\xrightarrow{\on{Id}_{A_{\sigma}}\otimes (\gamma_{i})_{i\in I}}
\Big(A_{\sigma} \otimes _{E}W_{\sigma^{I}}\Big)^{S_{\sigma}}
\xrightarrow{\Id_{A_{\sigma}}\otimes \xi} 
 (A_{\sigma})^{S_{\sigma}}
    \end{gather}
   c'est-à-dire  par le  produit par le scalaire  $\s{\xi, (\sigma(\gamma_{i}))_{i\in I} \cdot x }=f((\sigma(\gamma_{i}))_{i\in I} )$. Cela justifie 
  {\it   a posteriori} la définition des $S_{I,f,(\gamma_{i})_{i\in I}}$
   (et suggère que ces opérateurs sont diagonalisables, mais nous ne savons pas le démontrer). 
   
   La conjecture \eqref{intro-dec-I-W} n'est pas démontrée mais en suivant des idées de Drinfeld on peut montrer que les propriétés a) et b) de la \propref{prop-a-b-c} 
      impliquent   une décomposition 
 assez proche de \eqref{intro-dec-I-W} 
   (mais plus difficile à énoncer, car rempla\c cant la donnée de $\Sigma$ et des $A_{\sigma}$ par celle d'un ``$\mc O$-module sur le champ des paramètres de Langlands'').   Pour plus détails sur cette construction, voir \cite{texte-ICM}.

              Ce paragraphe  était heuristique et à partir de maintenant on oublie la conjecture \eqref{intro-dec-I-W}. 
     
\section{Cas non déployé, compléments et questions ouvertes }\label{para-commentaires}

\subsection{Cas où $G$ n'est pas pas nécessairement déployé}\label{cas-non-deploye}  On ne donne que les énoncés et on renvoie au chapitre 12 de \cite{coh} pour les démonstrations, qui ne font pas intervenir d'idées nouvelles par rapport au cas déployé. 

Soit $G$ un groupe réductif connexe sur $F$. 
Soit $\wt F$ une extension finie de $F$ déployant $G$ et ${}^{L}G=\wh G \rtimes \on{Gal}(\wt F/F)$ (où le produit semi-direct est pris pour l'action de $\on{Gal}(\wt F/F)$ par automorphismes de $\wh G$ préservant un épinglage). Soit $U$ un ouvert de $X$ sur lequel $G$ est réductif. En chaque point de $X\sm U$ on choisit  un modèle parahorique de Bruhat-Tits \cite{bruhat-tits} pour $G$,  si bien que $G$ est un schéma en groupes lisses sur $X$. Par commodité on suppose le niveau $N$ assez grand pour que $X\sm N \subset U$. 
On note $\Bun_{G,N}$ le champ sur $\Fq$  classifiant les $G$-torseurs sur $X$ avec structure de niveau $N$, autrement dit pour tout schéma $S$ sur $\Fq$, $\Bun_{G,N}(S)$ est le groupoïde classifiant les $G$-torseurs $\mc G$ sur $X\times S$ munis d'une trivialisation de $\restr{\mc G}{N\times S}$. Ce champ est lisse   (\cite{heinloth-unif}). On définit $K_{N}$ comme précédemment dans \eqref{def-K-N}. 
On a 
 \begin{gather}\label{dec-alpha-general-ker1}  \Bun_{G,N}(\Fq)=
  \bigcup_{\alpha\in \ker^{1}(F,G)}G_{\alpha}(F)\backslash G_{\alpha}(\mb A)/K_{N}\end{gather} où la réunion est  disjointe, $\ker^{1}(F,G)$ est fini et $G_{\alpha}$ est la  forme intérieure pure  de $G$  obtenue par torsion  par $\alpha$. Pour tout  $\alpha\in \ker^{1}(F,G)$ on a   $G_{\alpha}(\mathbb A)=G(\mathbb A)$ et donc le  quotient par $K_{N}$ dans le membre de droite a bien un sens.     
     
      On   fixe  un réseau $\Xi\subset Z(\mb A)/Z(F)$. 
   On définit    \begin{gather}\label{dec-alpha}C_{c}^{\mr{cusp}}(\Bun_{G,N}(\Fq)/\Xi,E)= \bigoplus_{\alpha\in \ker^{1}(F,G)}C_{c}^{\mr{cusp}}(G_{\alpha}(F)\backslash G_{\alpha}(\mb A)/K_{N}\Xi,E).\end{gather} 
   Alors les opérateurs d'excursion sont des endomorphismes 
  $$S_{I,f,(\gamma_{i})_{i\in I}}\in \mr{End}_{C_{c}(K_{N}\backslash G(\mb A)/K_{N},E)}( C_{c}^{\mr{cusp}}(\Bun_{G,N}(\Fq)/\Xi,E))$$
où  $I$ est un ensemble fini, $(\gamma_{i})_{i\in I}\in  \pi_{1}(X\sm N, \ov\eta)^{I}$   
 et    $f$ est une fonction sur le quotient grossier 
   $\wh G \backslash ({}^{L } G)^{I}/\wh G    $. 
  La méthode  pour les construire  est la même que  dans le cas déployé, grâce à   une variante tordue  sur $X\sm N$ de l'équivalence  de Satake géométrique  (cas non ramifié de \cite{richarz,zhu}). Cette variante 
    fait  intervenir ${}^{L} G$ parce que  l'épinglage de $\wh G$ apparaît naturellement 
   dans le foncteur fibre de Mirkovic-Vilonen
   (en effet ce foncteur fibre est donné par la cohomologie totale et  l'épinglagle est déterminé par la graduation 
   par le degré cohomologique et par 
   le cup-produit par le $c_{1}$ d'un fibré très ample sur la grassmannienne affine). 
Les opérateurs d'excursion  engendrent une sous-algèbre commutative $\mc B$ et par décomposition spectrale suivant les caractères de $\mc B$ on obtient une décomposition 
 \begin{gather}\label{intro2-dec-canonique}
 C_{c}^{\mr{cusp}}
 (\Bun_{G,N}(\Fq)/\Xi,\Qlbar)
 =\bigoplus_{\sigma}  \mf H_{\sigma}. \end{gather} 
  La somme directe dans le membre de droite est indexée par des paramètres de Langlands globaux, c'est-à-dire des classes de  $\wh G(\Qlbar)$-conjugaison de  morphismes 
       $\sigma:\on{Gal}(\ov F/F)\to {}^{L} G(\Qlbar)$ 
       définis sur une extension finie de  $\Ql$, continus,  semi-simples,  non ramifiés en dehors de  $N$ et rendant commutatif le diagramme  
        \begin{gather}\label{diag-sigma}
 \xymatrix{
\on{Gal}(\ov F/F) \ar[rr] ^{\sigma}
\ar[dr] 
&& {}^{L} G(\Qlbar) \ar[dl] 
 \\
& \on{Gal}(\wt F/F) }\end{gather}
Comme dans le 
\thmref{intro-thm-ppal}  
la décomposition \eqref{intro2-dec-canonique} est caractérisée par 
\eqref{relation-fonda}, et elle est compatible avec l'isomorphisme de Satake  (tordu) \cite{satake,cartier-satake,borel-corvallis,blasius-rogawski-pspm} en toutes les places de $X\sm N$.   

\begin{rem} Habituellement (par exemple dans les formules de multiplicités d'Arthur) on quotiente l'ensemble des morphismes $\sigma$ 
par une relation d'équivalence plus faible, qui, en plus  de la conjugaison par $\wh G(\Qlbar)$, autorise à tordre $\sigma$ par des éléments de $\ker^{1}(F,Z(\wh G)(\Qlbar))$. 
 D'après   Kottwitz~\cite{kottwitz1,kottwitz2} (et le théorème 2.6.1 de  Nguyen Quoc Thang~\cite{thang}   pour l'adaptation en caractéristique $p$), ce groupe fini est le dual de $\ker^{1}(F,G)$. Il a donc le même cardinal  et du point de vue des formules de multiplicités d'Arthur notre relation d'équivalence plus fine sur $\sigma$ compense exactement  le fait que l'espace que nous décomposons est une somme indexée par $ \alpha\in \ker^{1}(F,G)$. 
 Par exemple, si $G$ est un tore, les $\mf H_{\sigma}$ de \eqref{intro2-dec-canonique} sont de dimension $1$. 
 \end{rem}
 
 \begin{rem} Lorsque $G$ est déployé, $\ker^{1}(F,Z(\wh G)(\Qlbar))$ est nul par le théorème de Tchebotarev, donc 
  $\ker^{1}(F,G)$ est nul aussi et  
$\Bun_{G,N}(\Fq)=G (F)\backslash G (\mb A)/K_{N}$. C'est pourquoi le quotient $G(F)\backslash G(\mb A)/K_{N}\Xi$ intervenait dans le paragraphe \ref{para-enonces}, et interviendra de nouveau à partir du
 paragraphe \ref{defi-chtou-intro}  (où l'on reviendra au cas où $G$ est déployé pour simplifier la rédaction).  
\end{rem}

 On renvoie à la proposition 12.5 
 de \cite{coh} pour le fait que la décomposition \eqref{intro2-dec-canonique} est compatible avec les  isogénies de $G$ (et plus généralement avec  tous les morphismes $G\to G'$ dont l'image est distinguée). 

On espère  que la décomposition \eqref{intro2-dec-canonique}  est  compatible avec tous les cas connus de fonctorialité où l'on dispose d'un noyau explicite. En particulier on  devrait pouvoir montrer qu'elle est compatible avec 
  la    correspondance theta, grâce à la géométrisation du noyau theta  par Lysenko \cite{sergey-theta,sergey-theta-SO-Sp} et au 
 lien entre notre construction et  \cite{brav-var} (expliqué dans le paragraphe  \ref{subsection-link-langl-geom}).

\subsection{Paramètres d'Arthur} \label{para-param-Arthur} On aimerait montrer que 
   les  paramètres de Langlands $\sigma$ qui apparaissent dans la décomposition \eqref{intro1-dec-canonique} (ou \eqref{intro2-dec-canonique} dans le cas non déployé)   proviennent de  paramètres d'Arthur elliptiques.  On rappelle qu'un paramètre d'Arthur  est une classe  de  $\wh G(\Qlbar)$-conjugaison de     morphisme  $$
  \psi : \on{Gal}(\ov F/F) \times SL_{2}(\Qlbar)\to 
 {}^{L} G(\Qlbar) \text{ (algébrique sur  $SL_{2}(\Qlbar)$),}$$ dont la restriction à  $\on{Gal}(\ov F/F)$ prend ses valeurs  dans une extension finie de  $\Ql$, est continue, non ramifiée presque partout, semi-simple, pure de poids $0$,      et  fait commuter le triangle similaire à \eqref{diag-sigma}.  
 De plus  $ \psi$ est dit elliptique  si 
le centralisateur  de $\psi$ dans  $\wh  G(\Qlbar)$ est fini  modulo $(Z(\wh  G)(\Qlbar))^{\on{Gal}(\wt F/F)}$.  
  
 Le  paramètre de Langlands associé à $\psi$ est  $\sigma_{\psi}: \on{Gal}(\ov F/F) \to {}^{L} G(\Qlbar)$ défini par 
 $$\sigma_{\psi}(\gamma)=\psi\Big(\gamma, \begin{pmatrix} |\gamma|^{1/2} & 0 \\
 0 & |\gamma|^{-1/2}
 \end{pmatrix}\Big)$$  où $|\gamma|^{1/2}$ est bien défini grâce au choix d'une racine carrée  de $q$.   On conjecture que tout  paramètre 
 de Langlands  $\sigma$ apparaissant
  dans la  décomposition 
 \eqref{intro1-dec-canonique} (ou \eqref{intro2-dec-canonique} dans le cas non déployé)
    est de la forme   $\sigma_{\psi}$ 
     avec  $\psi$ un  paramètre d'Arthur 
     elliptique non ramifié sur $X\sm N$. 
  D'après \cite{kostant-betti} la classe de conjugaison   par  $\wh G(\Qlbar)$ de $\psi$ est déterminée de manière unique par celle de $\sigma$. 
    
    On aimerait en fait  obtenir une  décomposition canonique comme \eqref{intro1-dec-canonique} (ou \eqref{intro2-dec-canonique} dans le cas non déployé)  pour {\it toute la partie discrète} (et  pas seulement la partie cuspidale) et cette décomposition devrait être indexée par des paramètres  d'Arthur elliptiques. 
       
\subsection{Signification de la décomposition}   La décomposition 
\eqref{intro1-dec-canonique} 
est certainement plus fine  en général 
que celle obtenue  par diagonalisation  
des opérateurs de Hecke en les places  non ramifiées. Même en prenant en compte les classes d'isomorphisme  de représentations de 
 $C_{c}(K_{N}\backslash G(\mb A)/K_{N},\Qlbar)$ on ne récupère pas en général la décomposition \eqref{intro1-dec-canonique}, et bien que les  formules de multiplicités d'Arthur fassent intervenir une somme  sur les  paramètres d'Arthur, une telle  décomposition canonique semble inconnue  en général dans le cas des corps de nombres. 
    En effet d'après   Blasius \cite{blasius}, Lapid  \cite{lapid} et Larsen 
    \cite{larsen1,larsen2}, pour certains groupes $ G $ (y compris déployés)  la même représentation de $C_{c}(K_{N}\backslash G(\mb A)/K_{N},\Qlbar)$ peut apparaître dans des espaces $\mf H_{\sigma}$ différents, à cause du phénomène suivant. Il y a des exemples de groupes finis $ \Gamma $ et de morphismes $ \tau, \tau ': \Gamma \to \wh G (\Qlbar) $ tels que $ \tau $ et $ \tau' $ ne soient pas conjugués mais que pour tout  $ \gamma \in \Gamma $, $ \tau (\gamma) $ et $ \tau '(\gamma) $ soient conjugués.
 On s'attend alors à ce qu'il existe un    morphisme surjectif  $ \rho: \on {Gal} (\ov F / F) \to \Gamma $  partout non ramifié et une représentation $ (H_ {\pi}, \pi) $  de $ G (\mb A) $ tels que $ (H_ {\pi}) ^ {K_ {N}} $ apparaisse à la fois dans  $ \mf H_ {\tau \circ \rho} $ et $ \mf H_ {\tau' \circ \rho} $.
 Les exemples de Blasius et Lapid sont pour $ G = SL_{r} $, $ r \geq 3$
  (en fait dans ce cas on peut retrouver {\it a posteriori} la décomposition \eqref{intro1-dec-canonique} à l'aide du plongement  $SL_{r}\hookrightarrow GL_{r}$). 
    Mais pour certains  groupes (par exemple $E_{8}$)  
 nous ne savons pas comment retrouver la décomposition \eqref{intro1-dec-canonique} autrement que par les méthodes du présent  article, qui ne fonctionnent que sur les corps de fonctions.

 \subsection{Indépendance de $\ell$} On espère que  la décomposition   \eqref{intro1-dec-canonique} (ou \eqref{intro2-dec-canonique} dans le cas non-déployé) est définie sur $\ov \Q$, indépendante de $\ell$ (et du plongement $\ov \Q\subset \Qlbar$), et indexée par des paramètres de Langlands motiviques
 (la notion de paramètre de Langlands motivique est claire si l'on admet les conjectures standard mais on note que dans \cite{drinfeld-pro-completion} Drinfeld en a donné une définition inconditionnelle). 
 Cette conjecture paraît hors d'atteinte pour le moment. On renvoie à la conjecture   12.12  de \cite{coh} pour un énoncé plus précis. 
   
  \subsection{Cas des corps de nombres} 
 Il est évidemment  hors d'atteinte d'appliquer les méthodes de cet article aux  corps de nombres. Cependant on peut se demander s'il est raisonnable d'espérer une décomposition analogue à la  décomposition  {\it canonique} 
  \eqref{intro1-dec-canonique} (ou \eqref{intro2-dec-canonique} dans le cas non déployé) et  une formule de multiplicités d'Arthur pour chacun des espaces  $\mf H_{\sigma}$. 
Quand  $F$ est un corps de fonctions comme dans cet article, 
la limite   $ \varprojlim_{N } \Bun_{G,N}(\Fq)$ est égale à  $\big(G(\ov F)\backslash G(\mb A\otimes_{F}\ov F)   \big)^{\on{Gal}(\ov F/F)}$. Or cette dernière expression garde un sens 
pour les corps de nombres et elle  ne pose pas de  problèmes topologiques car on peut remarquer qu'elle est aussi égale 
\begin{itemize}\item à  
    $\big(G(\check F)\backslash G(\mb A\otimes_{F}\check F)\big)^{\on{Gal}(\check F/F)}$ où  $\check F$ est une  extension finie galoisienne de $F$ sur laquelle  $G$ est déployé
    \item à 
    $ \bigcup_{\alpha\in \ker^{1}(F,G)}G_{\alpha}(F)\backslash G_{\alpha}(\mb A)
$ 
où $\ker^{1}(F,G)$ est fini et $G_{\alpha}$ est une forme intérieure de $G$.     \end{itemize} 
  On peut espérer que si   $F$ est un corps de nombres  et si 
$\Xi$ est un réseau  dans  $Z(F)\backslash Z(\mb  A)$, 
la partie discrète    \begin{gather}\label{L2disc}L^{2}_{\mr{disc}}\Big(\big(G(\ov F)\backslash G(
\mb A\otimes_{F}\ov F)\big)^{\on{Gal}(\ov F/F)}/\Xi,\C\Big)\end{gather} admette une décomposition   {\it canonique}   indexée par les classes de conjugaison par  $\wh G(\C)$ de 
paramètres d'Arthur elliptiques. 
 Le cas particulier qui ressemble le plus au cas des corps de fonctions est celui des formes automorphes cohomologiques. 
En effet la partie cohomologique  de \eqref{L2disc} est définie sur  $\ov\Q$ et on peut se demander si elle admet une  décomposition {\it canonique}   sur  $\ov \Q$ indexée par les  classes d'équivalence de paramètres  d'Arthur  elliptiques motiviques (avec les subtilités de \cite{buzzard-gee} concernant la différence entre $L$-algébricité et $C$-algébricité). 
   
\subsection{Importance de la somme sur $\ker^{1}$} 
 Dans le cas non déployé on ne sait pas si l'inclusion 
 $$
  C_{c}^{\mr{cusp}}(G (F)\backslash G (\mb A)/K_{N}\Xi,\Qlbar)\subset C_{c}^{\mr{cusp}}(\Bun_{G,N}(\Fq)/\Xi,\Qlbar)$$  (correspondant au terme $\alpha=0$ dans le membre de droite de \eqref{dec-alpha-general-ker1}) 
   est compatible avec la décomposition \eqref{intro1-dec-canonique}, même après avoir regroupé les $\sigma$ différant par un élément de   $\ker^{1}(F,Z(\wh G)(\Qlbar))$. En dehors du cas où  
   la somme \eqref{dec-alpha-general-ker1} est réduite à un seul terme
 (c'est-à-dire lorsque $\ker^{1}(F,Z(\wh G)(\Qlbar))=0$ et donc en particulier lorsque   $G$  est une forme intérieure d'un groupe déployé), 
   on n'obtient donc pas de  décomposition canonique 
   pour l'espace $
  C_{c}^{\mr{cusp}}(G (F)\backslash G (\mb A)/K_{N}\Xi,\Qlbar)$. 
   
\subsection{Coefficients dans les corps finis} Grâce au fait que l'equivalence de Satake géométrique est définie sur $\mc O_{E}$ \cite{mv,ga-de-jong}, on peut montrer que 
lorsque la fonction $f$ est définie sur $\mc O_{E}$, $S_{I,f,(\gamma_{i})_{i\in I}}$ préserve 
  $ C_{c}^{\mr{cusp}}(\Bun_{G,N}(\Fq)/\Xi,\mc O_{E})$. Par décomposition spectrale des réductions de ces opérateurs modulo l'idéal maximal de $\mc O_{E}$ on obtient une décomposition  de 
   $ C_{c}^{\mr{cusp}}(\Bun_{G,N}(\Fq)/\Xi,\Flbar))$ indexée par les classes de $\wh G(\Flbar)$-conjugaison de morphismes $\sigma : \pi_{1}(X\sm N, \ov\eta)\to {}^{L} G(\Flbar)$ définis sur un corps fini, continus, faisant commuter un  diagramme analogue  à 
      \eqref{diag-sigma}, et complètement réductibles au sens de Jean-Pierre Serre \cite{bki-serre,bmr} (c'est-à-dire que si l'image est incluse dans un parabolique de ${}^{L} G$, elle est incluse dans un Levi associé). On renvoie au chapitre 13 de \cite{coh} pour plus de détails. 

\subsection{Paramétrisation locale}  La paramétrisation de Langlands locale (à semi-simplification près) et la compatibilité local-global sont traitées dans un article  avec Alain Genestier \cite{genestier-lafforgue}. 
Elles résultent de l'énoncé suivant : si $v$ est une place de $X$   et si tous les $\gamma_{i}$ appartiennent à $\on{Gal}(\ov F_{v}/F_{v})$, alors $S_{I,f,(\gamma_{i})_{i\in I}}$ est égal à l'action d'un élément du (complété $\ell$-adique du) centre de Bernstein de $G(F_{v})$ (évidemment le cas intéressant est celui où $v\in N$ car le cas non ramifié   est rendu  totalement explicite par la compatibilité avec l'isomorphisme de Satake dans le \thmref{intro-thm-ppal}).  
   
\subsection{Multiplicités} Ce travail ne détermine pas les multiplicités, et en particulier il ne dit pas pour quels  paramètres de Langlands $\sigma$ 
   l'espace $\mf H_{\sigma}   $ est non nul. 
          
     \subsection{Cas des groupes métaplectiques} 
     Le chapitre 14 de \cite{coh}  indique sommairement comment étendre 
  les résultats de cet article aux groupes métaplectiques, grâce à la variante  métaplectique de l'équivalence de Satake géométrique établie dans 
     \cite{finkelberg-lysenko, lysenko-red, dennis-sergey}. 
          
    \section{
    Chtoucas de Drinfeld pour les  groupes réductifs, d'après Varshavsky 
   } \label{defi-chtou-intro}
       
       Dans toute la suite de l'article $G$ est déployé
       (le cas non déployé, qui a été discuté dans le paragraphe \ref{cas-non-deploye} ci-dessus,  
       fait l'objet du chapitre 12 de \cite{coh} mais ne nécessite  pas d'idée supplémentaire). 
       
       Les ingrédients géométriques de notre construction sont expliqués dans ce paragraphe et le paragraphe \ref{subsection-crea-annihil-intro}. 
       En voici un aper\c cu rapide. La cohomologie d'intersection à support compact des champs de chtoucas  fournit pour tout ensemble fini $I$, pour tout niveau $N$ et pour toute représentation $W$ de $(\wh G)^{I}$ un système inductif  $\varinjlim_{\mu}\mc H_{N,I,W}^{\leq\mu}$ de  $E$-faisceaux constructibles sur  $(X\sm N)^{I}$. Le but de ce paragraphe est de construire ce système inductif, de fa\c con  fonctorielle en $W$,   de le munir d'actions des opérateurs de Hecke et des morphismes de Frobenius partiels $F_{\{i\}}$, et d'établir les isomorphismes de coalescence \eqref{intro-isom-coalescence} qui concernent sa restriction par un morphisme  diagonal  $(X\sm N)^{J}\to (X\sm N)^{I}$ 
             (associé à une application arbitraire $I\to J$). 
       Dans le paragraphe \ref{subsection-crea-annihil-intro} on montrera 
       que les opérateurs de Hecke aux places non ramifiées s'expriment à l'aide des isomorphismes de coalescence et des morphismes de Frobenius partiels. C'est cette propriété qui assurera la compatibilité de notre construction avec l'isomorphisme de Satake en les places non ramifiées.  Elle jouera aussi   un rôle technique fondamental en permettant d'étendre les opérateurs de Hecke en des morphismes de faisceaux sur $(X\sm N)^{I}$ tout entier, et en  fournissant les relations d'Eichler-Shimura. Ces relations   (qui seront énoncées dans la \propref{Eichler-Shimura-intro} ci-dessous) 
affirment que pour toute place  de $X\sm N$ et pour tout $i\in I$ la restriction en $x_{i}=v$ du  morphisme de Frobenius partiel $F_{\{i\}}$ est annulé par un polynôme dont les coefficients sont des opérateurs de Hecke en $v$ (à coefficients dans $\mc O_{E}$).  Elles  serviront    dans le paragraphe \ref{section-sous-faisceaux-constr} pour montrer que la propriété de finitude sous  l'action des opérateurs de Hecke (qui sert à définir les $H_{I,W}$) implique une finitude sous l'action des morphismes de Frobenius partiels, d'où, grâce à  un lemme fondamental de Drinfeld, l'action de $\on{Gal}(\ov F/F)^{I}$ sur $H_{I,W}$. Pour que le lecteur ne soit pas surpris nous signalons que cet usage des relations d'Eichler-Shimura 
est tout à fait 
         inhabituel. 
       
    Les chtoucas ont été introduits par  Drinfeld  \cite{drinfeld78,Dr1} pour  $GL_{r}$    et  généralisés  à tous les  groupes réductifs (et à tous les  copoids) par  Varshavsky dans  
    \cite{var} (entre-temps le cas des algèbres à  division a été étudié par  Laumon-Rapoport-Stuhler, Laurent Lafforgue, Ngô Bao Châu et  Eike Lau, voir les  réferences au début du paragraphe \ref{intro-previous-works}).
        Soit $I$  un ensemble fini 
 et  $W$ une    représentation $E$-linéaire irréductible de $(\wh G)^{I}$. On écrit   $W=\boxtimes_{i\in I}W_{i}$ où $W_{i}$ est une représentation irréductible de  $\wh G$.  
  Le champ $\Cht_{N,I,W}^{(I)}$ classifiant les  $G$-chtoucas avec structure de  niveau $N$, a été étudié dans  \cite{var}.  
  Contrairement à \cite{var} nous imposons dans la définition suivante  qu'il soit réduit 
  (cela ne change rien pour  la cohomologie étale). 
  
  \noindent{\bf Notation. }  Pour tout schéma  $S$ sur $\Fq$  et pour tout 
  $G$-torseur $\mc G$ sur  $ X\times S$ on note   $\ta \mc G=(\Id_{X}\times \Frob_{S})^{*}(\mc G)$. 
\begin{defi} On définit $\Cht_{N,I,W}^{(I)}$ comme le champ de Deligne-Mumford {\it réduit} sur $(X\sm N)^{I}$ dont 
 les points   sur un schéma  $S$ sur $\Fq$ classifient   
  \begin{itemize}
  \item des points  $(x_{i})_{i\in I}: S\to (X\sm N)^{I}$,  
\item un    $G$-torseur $\mc G$ sur  $ X\times S$, 
\item  un isomorphisme 
$$\phi :\restr{\mc G }{(X\times S)\sm(\bigcup_{i\in I }\Gamma_{x_i})}\isom \restr{\ta \mc G}{(X\times S)\sm(\bigcup_{i\in I }\Gamma_{x_i})}$$ 
 où   $\Gamma_{x_i}$ désigne le graphe de $x_{i}$, tel que la position relative en   $x_{i}$ soit bornée par le  copoids dominant de  $G$ correspondant au poids dominant $\omega_{i}$ de   $W_{i}$, 
 \item une trivialisation de   $(\mc G,\phi)$  sur  $N\times S$. 
 \end{itemize}
 \end{defi}
 
 Cette définition sera généralisée dans la \defiref{defi-Cht-I1-Ik} ci-dessous, et la condition sur la position relative sera précisée dans la \remref{position-relative}. 
   
 On note  $\Cht_{ I,W}^{(I)}$ lorsque $N$ est vide et on remarque que  $\Cht_{N,I,W}^{(I)}$ est un $G(\mc O_{N})$-torseur 
 sur $\restr{\Cht_{ I,W}^{(I)}}{(X\sm N)^{I}}$.

 \begin{rem}  Les personnes connaissant le programme de Langlands géométrique noteront que $\Cht_{N,I,W}^{(I)}$ est le produit fibré sur  $\Bun_{G,N}\times \Bun_{G,N}$   d'un champ de Hecke (considéré comme une correspondance entre $\Bun_{G,N}$ et lui-même) avec le graphe du morphisme de Frobenius de $\Bun_{G,N}$. 
 \end{rem}

Les  $x_{i}$ seront appelés les pattes du chtouca. On notera    $$\mf p_{N,I,W} ^{(I)}: \Cht_{N,I,W} ^{(I)}\to (X\sm N)^{I}$$ le morphisme correspondant. 
 
 Pour tout copoids dominant  $\mu$ de $G^{\mr{ad}}$ on note  
$\Cht_{N,I,W}^{(I),\leq\mu}$ l'ouvert de  $\Cht_{N,I,W}^{(I)}$ défini par la condition que le polygone de   Harder-Narasimhan   de $\mc G$ (ou plutôt, pour être précis, du $G^{\mr{ad}}$ -torseur associé) est $\leq \mu$ (au sens où la différence est  une combinaison à coefficients rationnels positifs  de  coracines simples).   

On fixe un réseau   $\Xi\subset Z(F)\backslash Z(\mb A)$. 
Alors $\Xi$ s'envoie dans $\Bun_{Z,N}(\Fq)$ qui agit sur $\Cht_{N,I,W}^{(I)}$ par torsion, et préserve les ouverts $\Cht_{N,I,W}^{(I),\leq\mu}$. 

On peut montrer que   $\Cht_{N,I,W}^{(I),\leq\mu}/\Xi$ est un champ de  Deligne-Mumford  de type fini. Il est  muni du morphisme 
$$\mf p_{N,I,W} ^{(I),\leq\mu}: \Cht_{N,I,W}^{(I),\leq\mu}/\Xi\to (X\sm N)^{I}$$
 qui se déduit de  $\mf p_{N,I,W} ^{(I)}$. 

 On note $\on{IC}_{\Cht_{N,I,W}^{(I),\leq \mu}/\Xi}$ le faisceau d'intersection  de  $\Cht_{N,I,W}^{(I),\leq \mu}/\Xi$ à coefficients dans  $E$, normalisé relativement à    $(X\sm N)^{I}$. 
La définition suivante sera rendue plus canonique (et notamment fonctorielle en $W$) dans la \defiref{defi-HNIW-can}. 

\begin{defi} \label{defi-HNIW-naive}
On pose 
$$\mc H_{N,I,W}^{\leq\mu}=R^{0}\big(\mf p_{N,I,W} ^{(I),\leq\mu}\big)_{!}
 \Big(\on{IC}_{\Cht_{N,I,W}^{(I),\leq \mu}/\Xi}\Big).$$
  \end{defi}

 Par rapport à   \cite{coh} on a allégé les notations, en enlevant les indices qui rappelaient que $\mc H_{N,I,W}^{\leq\mu}$ est une cohomologie   en degré $0$ (pour la normalisation perverse) et  à coefficients dans $E$. 
 La cohomologie est prise au sens de \cite{laumon-moret-bailly,laszlo-olsson} mais en fait la cohomologie étale des schémas suffirait. 
 En effet, dès que le degré de $N$ est suffisamment grand en fonction de $\mu$, $\Cht_{N,I,W}^{(I),\leq \mu}/\Xi$ est un schéma de type fini. 
  Donc, pour tout ouvert  $U\subset X\sm N$ tel que $U\subsetneq X$ (pour pouvoir augmenter $N$ sans changer    $U$),
  $\restr{\Cht_{N,I,W}^{(I),\leq \mu}/\Xi}{U^{I}}$   est le quotient d'un schéma de type fini par un groupe fini.

 Quand  $I$ est vide   et $W=\mbf 1$,  on a  
\begin{gather}\label{I-vide-intro} \restr{\varinjlim_{\mu}\mc H_{N,\emptyset,\mbf 1}^{\leq\mu}}{\Fqbar}=C_{c}(G(F)\backslash G(\mb A)/K_N \Xi,E)\end{gather}
  car  $\Cht_{N,\emptyset,\mbf 1}$ est le champ discret 
  $\Bun_{G,N}(\Fq)$, considéré comme un champ constant  sur $\Fq$, et  par ailleurs 
 $\Bun_{G,N}(\Fq)=G(F)\backslash G(\mb A)/K_N$
 (on utilise ici  l'hypothèse que   $G$ est déployé, en général par \eqref{dec-alpha-general-ker1}, $\Bun_{G,N}(\Fq)$ serait une réunion finie de quotients adéliques pour des formes intérieures de $G$).  
 
 \begin{rem} \label{rem-cht-W1} Plus généralement pour tout $I$ et  $W=\mbf  1$, le champ  $\Cht_{N,I,{\mbf  1} }^{(I)}/\Xi$ est simplement le champ constant 
  $G(F)\backslash G(\mb A)/K_{N}\Xi$ sur $(X\sm N)^{I}$.  \end{rem}
    
On considère $\varinjlim_{\mu}\mc H_{N,I,W}^{\leq\mu}$ comme un système inductif de    $E$-faisceaux constructibles sur  $(X\sm N)^{I}$.
On va  introduire maintenant les actions des morphismes de  Frobenius partiels et des  opérateurs de Hecke sur ce système inductif (on notera que ces actions augmentent  $\mu$).  
Pour toute partie  $J\subset I$ on note   $$
 \Frob_{J}:(X\sm N)^{I}\to (X\sm N)^{I}$$  le  morphisme qui à  
 $(x_{i})_{i\in I}$ associe   $(x'_{i})_{i\in I}$ avec   $$x'_{i}=\Frob(x_{i})\text{ \  si \  }
  i\in J\text{   \ et  \  }x'_{i}=x_{i} \text{ \ sinon.}$$  
 Alors on possède 
 \begin{itemize}
 \item pour $\kappa$ assez grand  et pour tout   $i\in I$, un  morphisme  
 \begin{gather}\label{intro-action-Frob-partiel}F_{\{i\}}:\Frob_{\{i\}}^{*}(\mc H_{N,I,W}^{\leq\mu})\to 
 \mc H_{N,I,W}^{\leq\mu+\kappa}\end{gather} 
 de faisceaux constructibles sur   $  (X\sm N)^{I}$, de sorte que les 
 $F_{\{i\}}$ commutent entre eux et que leur produit pour $i\in I$ est l'action naturelle du morphisme de Frobenius total de $(X\sm N)^{I}$ sur le faisceau $\mc H_{N,I,W}^{\leq\mu}$, 
  \item pour tout   $f\in C_{c}(K_{N}\backslash G(\mb A)/K_{N},E)$  et pour $\kappa$ assez grand, un  morphisme 
 \begin{gather}\label{defi-Tf}T(f):\restr{\mc H_{N,I,W}^{\leq\mu}}{(X\sm \mf P)^{I}}\to 
 \restr{\mc H_{N,I,W}^{\leq\mu+\kappa}}{(X\sm \mf P)^{I}}\end{gather} de faisceaux constructibles sur   $  (X\sm \mf P)^{I}$ où $\mf P$ est un ensemble fini de places contenant  $|N|$ et en dehors duquel   $f$ est triviale.  
 \end{itemize}
 
 Les morphismes $T(f)$ sont appelés des ``opérateurs de Hecke'' bien que ce soient des  morphismes de faisceaux. 
Ils sont obtenus grâce à la construction, assez évidente, de correspondances de Hecke entre les champs de chtoucas. 
  On verra après la \propref{prop-coal-frob-cas-part-intro}   que   $T(f)$ peut être étendu naturellement en un  morphisme de faisceaux sur   $(X\sm N)^{I}$,  mais cela n'est pas trivial. 
  Bien sûr lorsque $I=\emptyset$ et $W=\mbf 1$, les morphismes $T(f)$ sont les opérateurs de Hecke habituels sur \eqref{I-vide-intro}. 
 
  Pour construire les actions \eqref{intro-action-Frob-partiel} des morphismes de Frobenius partiels, on a besoin d'une petite généralisation des champs 
 $\Cht_{N,I,W} ^{(I)}$ où l'on demande une factorisation de  $\phi$ en une suite de plusieurs modifications. Soit $(I_{1},...,I_{k})$ une partition (ordonnée)  de $I$. 
 
 \begin{defi}\label{defi-Cht-I1-Ik}
 On définit     $\Cht_{N,I,W} ^{(I_{1},...,I_{k})}$ comme le champ de Deligne-Mumford {\it réduit }      dont les points sur 
 un schéma   $S$ sur  $\Fq$  classifient  les données    \begin{gather}\label{intro-donnee-chtouca}\big( (x_i)_{i\in I}, (\mc G_{0}, \psi_{0}) \xrightarrow{\phi_{1}}  (\mc G_{1}, \psi_{1}) \xrightarrow{\phi_{2}}
\cdots\xrightarrow{\phi_{k-1}}  (\mc G_{k-1}, \psi_{k-1}) \xrightarrow{ \phi_{k}}    (\ta{\mc G_{0}}, \ta \psi_{0})
\big)
\end{gather}
avec 
 \begin{itemize}
\item $x_i\in (X\sm N)(S)$ pour $i\in I$, 
\item pour $i\in \{0,...,k-1\}$, $(\mc G_{i}, \psi_{i})\in \Bun_{G,N}(S)$ (c'est-à-dire que  $\mc G_{i}$ est un   $G$-torseur sur  $X\times S$ et 
$\psi_{i} : \restr{\mc G_{i}}{N\times S} 
   \isom 
   \restr{G}{N\times S}$ est une trivialisation au-dessus de  $N\times S$) et on note  
   $(\mc G_{k}, \psi_{k})=(\ta{\mc G_{0}}, \ta \psi_{0})$
 \item  
pour   $j\in\{1,...,k\}$
 $$\phi_{j}:\restr{\mc G_{j-1}}{(X\times S)\sm(\bigcup_{i\in I_{j}}\Gamma_{x_i})}\isom \restr{\mc G_{j}}{(X\times S)\sm(\bigcup_{i\in I_{j}}\Gamma_{x_i})}$$ est un   isomorphisme  tel que la position relative de  $\mc G_{j-1}$ par rapport à  $\mc G_{j}$ en  $x_{i}$ (pour  $i\in I_{j}$) soit bornée par le  copoids dominant de  $G$ correspondant au  poids dominant de  $W_{i}$, 
 \item les $\phi_{j}$, qui induisent des isomorphismes sur $N\times S$, respectent les structures de niveau, c'est-à-dire que 
$\psi_{j}\circ \restr{\phi_{j}}{N\times S}=\psi_{j-1}$ pour tout  $j\in\{1,...,k\}$. 
 \end{itemize}
 \end{defi}
 La condition sur la position relative sera précisée dans la \remref{position-relative}. On note $\Cht_{N,I}^{(I_{1},...,I_{k})}$ l'ind-champ obtenu en oubliant cette condition. 
 
Pour tout copoids dominant $\mu$ de $G^{\mr{ad}}$  on note  
 $\Cht_{N,I,W}^{(I_{1},...,I_{k}),\leq\mu}$   l'ouvert de  $\Cht_{N,I,W}^{(I_{1},...,I_{k})}$ défini par la condition que  le polygone de Harder-Narasimhan  de  $\mc G_{0}$ est  $\leq\mu$. On note   $$\mf p_{N,I,W} ^{(I_{1},...,I_{k})}: \Cht_{N,I,W} ^{(I_{1},...,I_{k})}\to (X\sm N)^{I}$$ le morphisme   qui à un chtouca associe la famille de ses pattes. 
  
    \noindent{\bf Exemple.} Lorsque  $G=GL_r$,  $I=\{1,2\}$ et $W=\mr{St}\boxtimes \mr{St}^{*}$, les champs 
 $\Cht_{N,I,W}^{(\{1\}, \{2\})}$, {\it resp.} 
  $\Cht_{N,I,W}^{(\{2\}, \{1\})}$ 
 sont les champs de chtoucas à gauche, {\it resp.} à droite  introduits par  Drinfeld (et utilisés aussi dans  \cite{laurent-inventiones}), 
 et $x_{1}$ et  $x_{2}$ sont le pôle et le zéro.

On  construit maintenant un     morphisme  lisse  
\eqref{lisse-chtouca-grass-intro} de $\Cht_{N,I,W}^{(I_{1},...,I_{k})}$ vers le  quotient d'une  strate  fermée  d'une  grassmannienne  affine  de Beilinson-Drinfeld  par un schéma   en groupes lisse (cela permettra en plus dans la \remref{position-relative} de préciser la condition sur les positions relatives dans la \defiref{defi-Cht-I1-Ik}).   
Les lecteurs familiers avec les variétés de Shimura peuvent considérer ce morphisme comme un ``modèle local'' à condition de noter
\begin{itemize}
\item que l'on est dans une situation de bonne réduction puisque les $x_{i}$ appartiennent à $X\sm N$, 
\item et que pourtant ce modèle local n'est pas lisse (sauf si tous les $I_{j}$ sont des singletons et tous  les copoids sont minuscules). 
\end{itemize}

\begin{defi}
La grassmannienne  affine  de Beilinson-Drinfeld est l'ind-schéma  $\mr{Gr}_{I }^{(I_{1},...,I_{k})}$ sur $X^{I}$ dont les $S$-points classifient la donnée de  
\begin{gather}\label{formule-rem-grassm-intro}\big((x_{i})_{i\in I}, \mc G_{0} \xrightarrow{\phi_{1}}  
\mc G_{1}\xrightarrow{\phi_{2}}
\cdots\xrightarrow{\phi_{k-1}}  \mc G_{k-1} \xrightarrow{ \phi_{k}}   \mc G_{k}\isor{\theta} G_{X\times S} \big)  \end{gather}
où les $\mc G_{i}$ sont des $G$-torseurs sur $X\times S$, $\phi_{i}$ est un isomorphisme sur $(X\times S)\sm(\bigcup_{i\in I_{j}}\Gamma_{x_i})$ et $\theta$ est une trivialisation de $\mc G_{k}$. 
 La strate  fermée   $\mr{Gr}_{I,W}^{(I_{1},...,I_{k})}$ est le  sous-schéma fermé réduit de 
  $\mr{Gr}_{I }^{(I_{1},...,I_{k})}$ défini par la condition que 
  la position relative de  $\mc G_{j-1}$ par rapport à  $\mc G_{j}$ en  $x_{i}$ (pour  $i\in I_{j}$) est  bornée par le  copoids dominant de  $G$ correspondant au  poids dominant $\omega_{i}$ de  $W_{i}$. Plus précisément 
  au-dessus de l'ouvert $\mc U$ de $X^{I}$ où les $x_{i}$ sont deux à deux distincts, $\mr{Gr}_{I }^{(I_{1},...,I_{k})}$ est un produit de grassmanniennes affines usuelles et 
  \begin{itemize}
  \item on définit 
  la restriction de $\mr{Gr}_{I,W}^{(I_{1},...,I_{k})}$ au-dessus de $\mc U$  comme le produit des strates fermées habituelles (notées $\ov {\mr{Gr}_{\omega_{i}}}$ dans  \cite{mv,brav-gaitsgory}), 
  \item puis 
    on   définit  $\mr{Gr}_{I,W}^{(I_{1},...,I_{k})}$ comme l'adhérence de Zariski (dans $\mr{Gr}_{I }^{(I_{1},...,I_{k})}$) de 
  sa restriction au-dessus de $\mc U$. 
  \end{itemize}
  \end{defi}
  
   D'après Beauville-Laszlo \cite{BL} (voir aussi le premier paragraphe   de \cite{coh} pour des références complémentaires  dans \cite{hitchin}), 
   $\mr{Gr}_{I }^{(I_{1},...,I_{k})}$ peut aussi être défini comme l'ind-schéma dont les $S$-points classifient 
       \begin{gather}\label{formule-rem-grassm-intro-loc}\big((x_{i})_{i\in I}, \mc G_{0} \xrightarrow{\phi_{1}}  
\mc G_{1}\xrightarrow{\phi_{2}}
\cdots\xrightarrow{\phi_{k-1}}  \mc G_{k-1} \xrightarrow{ \phi_{k}}   \mc G_{k}\isor{\theta} G_{\Gamma_{\sum \infty x_i}} \big)  \end{gather}
 où les $\mc G_{i}$ sont des $G$-torseurs  sur le voisinage formel 
  $\Gamma_{\sum \infty x_i}$ de la réunion des graphes des $x_{i}$ dans
 $  X\times S$, $\phi_{i}$ est un isomorphisme sur $\Gamma_{\sum \infty x_i}\sm(\bigcup_{i\in I_{j}}\Gamma_{x_i})$ et $\theta$ est une trivialisation de $\mc G_{k}$. 
La restriction à la Weil $G_{\sum \infty x_i}$ de $G$ de $\Gamma_{\sum \infty x_i}$ à $S$ agit donc sur  $\mr{Gr}_{I }^{(I_{1},...,I_{k})}$  et $\mr{Gr}_{I,W}^{(I_{1},...,I_{k})}$ par changement de la trivialisation $\theta$. 

On a un morphisme naturel 
    \begin{gather}\label{lisse-chtouca-grass-intro0}
    \Cht_{N,I,W}^{(I_{1},...,I_{k})}\to \mr{Gr}_{I,W}^{(I_{1},...,I_{k})}/
  G_{\sum \infty  x_i}\end{gather} 
  qui associe à un chtouca \eqref{intro-donnee-chtouca} le 
  $G_{\sum \infty x_{i}}$-torseur $\restr{\mc G_{k}}{\Gamma_{\sum \infty x_{i}}}$
  et, pour toute trivialisation $\theta$ de celui-ci, le  point de $\mr{Gr}_{I,W}^{(I_{1},...,I_{k})}$ égal à  \eqref{formule-rem-grassm-intro-loc}.

 \begin{rem} \label{position-relative} La meilleure fa\c con d'énoncer la condition sur les positions relatives dans la \defiref{defi-Cht-I1-Ik} est de {\it définir}   $\Cht_{N,I,W}^{(I_{1},...,I_{k})}$ comme l'image inverse de 
  $\mr{Gr}_{I,W}^{(I_{1},...,I_{k})}/
  G_{\sum \infty x_i}$ par le morphisme   
  $\Cht_{N,I}^{(I_{1},...,I_{k})}\to \mr{Gr}_{I}^{(I_{1},...,I_{k})}/
  G_{\sum \infty x_i}$ construit  comme \eqref{lisse-chtouca-grass-intro0}. 
    \end{rem}

Pour $(n_{i})_{i\in I}\in \N^{I}$ on note $\Gamma_{\sum n_{i} x_i}$ le sous-schéma fermé de $X\times S$ associé au diviseur de Cartier $\sum n_{i} x_i$ qui est  effectif et relatif  sur $S$. On note $G_{\sum n_{i} x_i}$ le schéma en groupes lisse sur $S$ obtenu par restriction à la Weil de $G$ de $\Gamma_{\sum n_{i} x_i}$ à $S$. 
 Alors si les entiers $n_{i}$ sont assez grands en fonction de $W$, l'action de $G_{\sum \infty x_i}$ sur $\mr{Gr}_{I,W}^{(I_{1},...,I_{k})}$ se factorise par 
  $G_{\sum n_{i} x_i}$. 
 Le  morphisme \eqref{lisse-chtouca-grass-intro0} fournit donc un morphisme 
    \begin{gather}\label{lisse-chtouca-grass-intro}
    \Cht_{N,I,W}^{(I_{1},...,I_{k})}\to \mr{Gr}_{I,W}^{(I_{1},...,I_{k})}/
  G_{\sum n_{i} x_i}\end{gather} 
  (qui associe à un chtouca \eqref{intro-donnee-chtouca} le 
  $G_{\sum n_{i}x_{i}}$-torseur $\restr{\mc G_{k}}{\Gamma_{\sum n_{i}x_{i}}}$
  et, pour toute trivialisation $\lambda$ de celui-ci, le  point de $\mr{Gr}_{I,W}^{(I_{1},...,I_{k})}$ égal à  \eqref{formule-rem-grassm-intro-loc} pour toute trivialisation $\theta$ de $\restr{\mc G_{k}}{\Gamma_{\sum \infty x_i}}$ prolongeant  $\lambda$ de $\Gamma_{\sum n_{i} x_i}$ à 
  $\Gamma_{\sum \infty   x_i}$).

 On montre  dans la proposition 2.8 de \cite{coh}  que le morphisme 
  \eqref{lisse-chtouca-grass-intro} 
    est  lisse de  dimension $\dim G_{\sum n_{i}x_{i}}=(\sum_{i\in I} n_{i})\dim G$.  
    
 On en déduit  que le  morphisme d'oubli des modifications intermédiaires 
 \begin{gather}
 \label{oubli-Cht}
 \Cht_{N,I,W}^{(I_{1},...,I_{k})}\to \Cht_{N,I,W}^{(I)} \\ \nonumber 
 \text{   qui envoie \eqref{intro-donnee-chtouca} sur  } 
\big( (x_i)_{i\in I}, (\mc G_{0}, \psi_{0}) \xrightarrow{\phi_{k} \cdots \phi_{1}}      (\ta{\mc G_{0}}, \ta \psi_{0})
\big) \end{gather}
 est petit. En effet il est connu que le morphisme analogue 
\begin{gather}\label{mor-Gr-oubli}\mr{Gr}_{I,W}^{(I_{1},...,I_{k})}\to \mr{Gr}_{I,W}^{(I)}  \text{    qui envoie \eqref{formule-rem-grassm-intro-loc} sur }
\big( (x_i)_{i\in I},  \mc G_{0}  \xrightarrow{\phi_{k} \cdots \phi_{1}}       \mc G_{k}\isor{\theta} G_{\Gamma_{\sum \infty x_i}} 
\big)  \end{gather}   
est petit, et d'ailleurs cela joue un rôle essentiel dans   \cite{mv}. De plus l'image inverse de $\Cht_{N,I,W}^{(I),\leq\mu}$ 
par \eqref{oubli-Cht} est exactement $\Cht_{N,I,W}^{(I_{1},...,I_{k}),\leq\mu}$ 
puisque les troncatures ont été définies à l'aide du polygône de Harder-Narasimhan de $\mc G_{0}$. 
 On a donc  $$\mc H_{N,I,W}^{\leq\mu}=
   R^{0}(\mf p_{N,I,W} ^{(I_{1},...,I_{k}),\leq\mu})_{!}\Big(\on{IC}_{\Cht_{N,I,W}^{(I_{1},...,I_{k}),\leq\mu}/\Xi}\Big)$$ pour {\it toute} partition $(I_{1},...,I_{k})$ de $I$ (alors que la \defiref{defi-HNIW-naive} utilisait la partition grossière $(I)$).  

Le  morphisme de  Frobenius partiel 
$$\on {Fr}_{I_{1}} ^{(I_{1},...,I_{k})}: \Cht_{N,I,W} ^{(I_{1},...,I_{k})}\to \Cht_{N,I,W} ^{(I_{2},...,I_{k},I_{1})},$$ 
 défini par 
\begin{gather}\on {Fr}_{I_{1}} ^{(I_{1},...,I_{k})}\big( (x_i)_{i\in I}, (\mc G_{0}, \psi_{0}) \xrightarrow{\phi_{1}}  (\mc G_{1}, \psi_{1}) \xrightarrow{\phi_{2}}
\cdots\xrightarrow{\phi_{k-1}}  (\mc G_{k-1}, \psi_{k-1}) \xrightarrow{ \phi_{k}}    (\ta{\mc G_{0}}, \ta \psi_{0})
\big)\nonumber \\
= \nonumber
\big( \Frob_{I_{1}}\big((x_i)_{i\in I}\big), (\mc G_{1}, \psi_{1}) \xrightarrow{\phi_{2}}  (\mc G_{2}, \psi_{2}) \xrightarrow{\phi_{3}}
\cdots \xrightarrow{ \phi_{k}}    (\ta{\mc G_{0}}, \ta \psi_{0}) \xrightarrow{\ta  \phi_{1}   } (\ta{\mc G_{1}}, \ta \psi_{1})
\big)  
\end{gather}
est au-dessus du morphisme  $\Frob_{I_{1}}:(X\sm N)^{I}\to (X\sm N)^{I}$. 
La composée des morphismes $\on {Fr}_{I_{1}} ^{(I_{1},...,I_{k})}$, $\on {Fr}_{I_{2}} ^{(I_{1},...,I_{k})}$, ..., $\on {Fr}_{I_{k}} ^{(I_{1},...,I_{k})}$ est égale au morphisme de Frobenius total de $\Cht_{N,I,W} ^{(I_{1},...,I_{k})}$ sur $\Fq$. 
Comme   $\on {Fr}_{I_{1}} ^{(I_{1},...,I_{k})}$
est un  homéomorphisme local totalement radiciel, on a un isomorphisme    \begin{gather}\label{action-Frob-partiels}\Big(\on {Fr}_{I_{1}} ^{(I_{1},...,I_{k})}\Big)^{*}\Big(\on{IC}_{\Cht_{N,I,W} ^{(I_{2},...,I_{k},I_{1})}}\Big)=\on{IC}_{ \Cht_{N,I,W} ^{(I_{1},...,I_{k})}}   \end{gather}
que l'on normalise par un coefficient $q^{-d/2}$, où $d$ est la dimension relative de $\mr{Gr}_{I_{1},\boxtimes _{i\in I_{1}}W_{i}}^{(I_{1})}$ sur $X^{I_{1}}$ 
(cette normalisation sera justifiée dans la \remref{rem-Frob-partiel}). 
Grâce à l'isomorphisme de changement de base propre 
(et au  fait que les homéomorphismes locaux totalement radiciels ne changent pas la topologie étale), ou par une correspondance cohomologique facile donnée dans \cite{coh}, 
on en déduit alors de \eqref{action-Frob-partiels} un morphisme 
\begin{gather}\label{morph-F1}F_{I_{1}}:\Frob_{I_{1}}^{*}(\mc H_{N,I,W}^{\leq\mu})\to 
 \mc H_{N,I,W}^{\leq\mu+\kappa}\end{gather} pour $\kappa$ assez grand
 (en effet, dans les notations ci-dessus, si le polygône de Harder-Narasimhan de $\mc G_{1}$ est $\leq \mu$, comme la modification entre $\mc G_{0}$ et $\mc G_{1}$ est bornée en fonction de $W$, le polygône de Harder-Narasimhan de $\mc G_{0}$ est $\leq \mu+\kappa$ où $\kappa$ dépend de $W$). En prenant  n'importe quelle partition $(I_{1},...,I_{k})$  telle que  $I_{1}=\{i\}$ on obtient $F_{\{i\}}$ dans  \eqref{intro-action-Frob-partiel}.

 Pour l'instant on a défini $\mc H_{N,I,W}^{\leq\mu}$ pour les classes d'isomorphismes de représentations irréductibles $W$ de $(\wh G)^{I}$. On va raffiner  cette construction en celle plus {\it canonique} d'un {\it  foncteur} $E$-linéaire
 \begin{gather}\label{fonc-W-HIW-intro}W\mapsto \mc H_{N,I,W}^{\leq\mu} \end{gather} 
 de la catégorie des représentations $E$-linéaires de dimension finie de $(\wh G)^{I}$ vers la catégorie des $E$-faisceaux constructibles sur $(X\sm N)^{I}$. En particulier pour tout    morphisme $u:W\to W'$ de  représentations $E$-linéaires de $(\wh G)^{I}$  on notera 
 $$\mc H(u): \mc H_{N,I,W}^{\leq\mu} \to \mc H_{N,I,W'}^{\leq\mu}$$  
le morphisme de faisceaux constructibles  associé. 

 Le foncteur \eqref{fonc-W-HIW-intro} sera compatible à la coalescence des pattes, au sens suivant. Dans tout cet article nous appelons coalescence la situation où des pattes fusionnent entre elles. 
 Nous aurions pu employer  le mot fusion plutôt que coalescence
 mais nous avons préféré utiliser le mot coalescence pour les pattes
 (qui ne sont que des points sur la courbe) en gardant le mot fusion 
 pour le produit de fusion (qui intervient dans l'équivalence de Satake géométrique et concerne les faisceaux pervers sur les grassmanniennes affines de Beilinson-Drinfeld). 
    Soit  $\zeta: I \to J$ une application.  On note $W^{\zeta}$ la représentation de  $\wh G^{J}$ qui est la composée de la  représentation $W$ avec le morphisme diagonal
 $$ \wh G^{J}\to \wh G^{I}, (g_{j})_{j\in J}\mapsto (g_{\zeta(i)})_{i\in I}.$$
   On note 
   \begin{gather}\label{morph-giad-X-intro}\Delta_{\zeta}: X^{J}\to X^{I},(x_{j})_{j\in J}\mapsto (x_{\zeta(i)})_{i\in I}\end{gather} le   morphisme diagonal. 
On va construire, d'après  \cite{var} et \cite{brav-var}   un isomorphisme  de faisceaux constructibles sur $(X\sm N)^{J}$, dit de  {\it coalescence}:  \begin{gather}\label{intro-isom-coalescence}\chi_{\zeta}: \Delta_{\zeta}^{*}(\mc H_{N,I,W}^{\leq\mu})\isom 
 \mc H^{\leq\mu}_{N,J,W^{\zeta}}.\end{gather} 
 Cet isomorphisme sera {\it canonique} au sens où ce sera un isomorphisme de foncteurs en $W$, compatible avec la composition de $\zeta$.  
 On explique  maintenant  la construction du foncteur \eqref{fonc-W-HIW-intro} et de l'isomorphisme de coalescence \eqref{intro-isom-coalescence}. 
 
 Lorsque $W$ n'est pas irréductible on 
 note $\mr{Gr}_{I,W }^{(I_{1},...,I_{k})}$  la réunion des $\mr{Gr}_{I,V }^{(I_{1},...,I_{k})}\subset \mr{Gr}_{I }^{(I_{1},...,I_{k})}$ pour $V$ constituant irréductible de $W$. On fait de même avec 
  $\Cht_{N,I,W}^{(I_{1},...,I_{k})}$.

   On rappelle maintenant l'équivalence de Satake géométrique, due à Lusztig, Drinfeld, Ginzburg et Mirkovic-Vilonen. 
   Pour des références on cite  \cite{lusztig-satake,ginzburg,hitchin,mv,ga-iwahori,ga-de-jong,richarz,zhu}. On va utiliser  ici     la  forme  expliquée par Gaitsgory dans \cite{ga-de-jong}.  Habituellement l'équivalence de Satake géométrique s'exprime de la manière suivante. Pour tout corps $k$ algébriquement clos de caractéristique première à $\ell$, 
    la catégorie des faisceaux pervers $G(k[[z]])$-équivariants sur la grassmannienne affine $G(k((z)))/G(k[[z]])$ est munie 
    \begin{itemize}
    \item
    d'une structure tensorielle par le produit de fusion (ou de convolution), 
    avec une modification des signes dans la contrainte de commutativité rappelée dans la \remref{modif-comm} ci-dessous, 
\item      et d'un foncteur fibre, donné par la cohomologie totale, 
qui grâce à cette modification est tensoriel à valeurs dans la catégorie des espaces vectoriels, et non des super-espaces vectoriels  (pour être plus canonique  il faut introduire en plus  une torsion à la Tate, c'est-à-dire tensoriser par $E(\frac{i}{2})$ la partie de   degré cohomologique $i$ pour tout $i\in \Z$). 
\end{itemize}
    Elle est donc équivalente à la catégorie des représentations 
   du groupe des automorphismes du foncteur fibre, qui s'avère être isomorphe à  
   $\wh G$ (et muni d'un épinglage canonique). 
   De plus ses objets sont naturellement équivariants par le groupe des automorphismes de $k[[z]]$. Ceci permet de remplacer 
   $\on{Spf}(k[[z]])$ par un disque formel arbitraire, et en particulier un disque  formel variant sur une courbe. 
   
   Ici on   utilise  seulement un sens de l'équivalence, à savoir le foncteur de la catégorie des représentations de $\wh G$ vers la catégorie des faisceaux pervers 
   $G(k[[z]])$-équivariants  sur la grassmiannenne affine. En revanche  on l'énonce en utilisant la grassmannienne affine de Beilinson-Drinfeld. 
   Le fait que le foncteur fibre qui fournit l'équivalence de Satake géométrique est   donné par la cohomologie totale 
  implique, dans les notations du théorème ci-dessous, que $W$ est canoniquement égal à la cohomologie totale  de $\mc S_{I,W }^{(I_{1},...,I_{k})}$ dans les fibres de $\mr{Gr}_{I }^{(I_{1},...,I_{k})}$ au-dessus de $X^{I}$. Cependant nous n'utilisons pas cette propriété en tant que telle, seulement à travers le fait qu'elle fournit le foncteur fibre dans l'équivalence de Satake géométrique et   donc la canonicité de nos constructions.

 \begin{thm}  \label{satake-geom-thm} (un sens de l'équivalence de Satake géométrique \cite{hitchin, mv, ga-de-jong}). On possède pour tout ensemble fini $I$ et pour toute partition $(I_{1},...,I_{k})$ de $I$ un   foncteur   $E$-linéaire   
   $$W\mapsto \mc S_{I,W }^{(I_{1},...,I_{k})}$$
de la catégorie des  représentations $E$-linéaires de dimension finie de $(\wh G)^{I}$ vers la  catégorie des    $E$-faisceaux pervers  
$G_{\sum \infty x_{i}}$-équivariants sur $\mr{Gr}_{I }^{(I_{1},...,I_{k})}$ (pour la normalisation perverse  relative à $X^{I}$). De plus $\mc S_{I,W}^{(I_{1},...,I_{k})}$ est supporté par 
$\mr{Gr}_{I,W}^{(I_{1},...,I_{k})}$ et est universellement localement acyclique relativement à $X^{I}$. Ces foncteurs vérifient les propriétés suivantes.  

\begin{itemize}
\item [] a) Compatibilité aux morphismes d'oubli (des modifications intermédiaires)  : $S_{I,W}^{(I)}$ est canoniquement isomorphe à l'image directe de 
$\mc S_{I,W}^{(I_{1},...,I_{k})}$ par le morphisme d'oubli  
$\mr{Gr}_{I}^{(I_{1},...,I_{k})}\to \mr{Gr}_{I}^{(I)}$ (défini dans \eqref{mor-Gr-oubli}).   
\item [] b)  Compatibilité à la convolution : si $W=\boxtimes_{j\in \{1,...,k\}} W_{j}$ où $W_{j}$ est une  représentation de $(\wh G)^{I_{j}}$, 
  $\mc S_{I,W}^{(I_{1},...,I_{k})}$  est canoniquement isomorphe à l'image inverse de $\boxtimes _{j\in \{1,...,k\}} \mc S_{I_{j},W_{j}}^{(I_{j})}$ par le morphisme 
  \begin{align*}
  \mr{Gr}_{I}^{(I_{1},...,I_{k})}/G_{\sum_{i\in I} \infty x_{i}} &\to 
  \prod_{j=1}^{k}\Big( \mr{Gr}_{I_{j}}^{(I_{j})}/G_{\sum_{i\in I_{j}} \infty x_{i}}\Big) \\
  (\mc G_{0}\to \mc G_{1} \to  \cdots \to  \mc G_{k})& 
  \mapsto \Big(\big(\restr{\mc G_{j-1}}{\Gamma_{\sum_{i\in I_{j}} \infty x_{i}}}\to \restr{\mc G_{j}}{\Gamma_{\sum_{i\in I_{j}} \infty x_{i}}} \big) \Big)_{j=1,...,k}
  \end{align*}
  où les $\mc G_{i}$ sont des $G$-torseurs sur 
  $\Gamma_{\sum_{i\in I} \infty x_{i}}$. 
   \item [] c)  Compatibilité à la fusion: soient  $I,J$ des ensembles finis et  $\zeta: I\to J$ une application. Soit $(J_{1},...,J_{k})$ une partition de $J$. Son image inverse 
   $(\zeta^{-1}(J_{1}), ..., \zeta^{-1}(J_{k}))$ est une partition de $I$. 
On note  $$\Delta_{\zeta} : X^{J}\to X^{I},  \ \ (x_{j})_{j\in J}\mapsto (x_{\zeta(i)})_{i\in I}$$ le   morphisme diagonal associé à 
  $\zeta$. On note encore 
  $\Delta_{\zeta}$    l'inclusion  
 $$\mr{Gr}_{J}^{(J_{1},...,J_{k})} =\mr{Gr}_{I}^{(\zeta^{-1}(J_{1}), ..., \zeta^{-1}(J_{k}))}\times _{X^{I}}X^{J}\hookrightarrow \mr{Gr}_{I}^{(\zeta^{-1}(J_{1}), ..., \zeta^{-1}(J_{k}))}.$$
   Soit  $W$ une   représentation  $E$-linéaire de dimension finie   de  $\wh G^{I}$. On note  $W^{\zeta}$ la   représentation de $\wh G^{J}$ qui est la composée de la   représentation $W$ avec le    morphisme diagonal $$\wh G^{J}\to \wh G^{I},  \ \ (g_{j})_{j\in J}\mapsto (g_{\zeta(i)})_{i\in I}  .$$   On a alors un  isomorphisme canonique 
 \begin{gather}\label{coalescence-Gr-section1-thm} \Delta_{\zeta}^{*}\Big( \mc S_{I,W}^{(\zeta^{-1}(J_{1}), ..., \zeta^{-1}(J_{k}))} \Big)\simeq 
 \mc S_{J,W^{\zeta}}^{(J_{1},...,J_{k})}\end{gather}  
 qui est fonctoriel en  $W$ et compatible avec la composition pour $\zeta$. 
  \item [] d) 
 Lorsque  
$W$ est irréductible, le faisceau pervers  $\mc S_{I,W}^{(I_{1},...,I_{k})}$ sur 
 $\mr{Gr}_{I,W}^{(I_{1},...,I_{k})}$ est isomorphe au  faisceau d'intersection   (avec la  normalisation perverse relative à  $X^{I}$). 
 \end{itemize}  
\end{thm}
  
  Les propriétés a) et b) auraient pu être énoncées avec des partitions plus générales, mais au prix de notations beaucoup plus lourdes. 
  
 Dans le théorème précédent  $\mc S_{I,W}^{(I_{1},...,I_{k})}$ est supporté par 
$\mr{Gr}_{I,W}^{(I_{1},...,I_{k})}$ et on peut donc le considérer comme un faisceau pervers 
(à un décalage près)   sur  $\mr{Gr}_{I,W}^{(I_{1},...,I_{k})}/G_{\sum n_{i}x_{i}}$  (avec les entiers $n_{i}$   assez grands).

   \begin{rem}\label{modif-comm}
 La contrainte de commutativité est définie par la convention  {\it modifiée},  introduite dans la discussion qui précède la   proposition 6.3 de \cite{mv} (et expliquée aussi dans  \cite{hitchin}). Voici un bref rappel. 
  Le lemme 3.9 de   \cite{mv} montre que 
   \begin{itemize}
   \item pour toute composante  connexe  de la grassmannienne affine, les strates sont toutes de dimension paire ou toutes de dimension impaire (on parle alors de composante paire ou impaire), 
   \item si un faisceau $\mc S_{I_{j},W_{j},\mc A}^{(I_{j})}$ est supporté sur une composante  paire ({\it resp.} impaire)   sa cohomologie totale est concentrée en degrés cohomologiques pairs ({\it resp.} impairs). 
   \end{itemize}
  La contrainte de commutativité modifée consiste à   ajouter aux signes habituels donnés par les règles de Koszul  un signe moins lorsque l'on permute deux faisceaux 
$\mc S_{I_{j},W_{j},\mc A}^{(I_{j})}$ et $\mc S_{I_{j'},W_{j'},\mc A}^{(I_{j'})}$ supportés 
sur deux composantes connexes impaires de la grassmannienne affine.  
Autrement dit c'est la contrainte de commutativité naturelle que l'on aurait si on 
normalisait les faisceaux $\mc S_{I_{j},W_{j},\mc A}^{(I_{j})}$ pour que leur cohomologie totale   soit en degrés cohomologiques pairs.     \end{rem}

 \begin{rem}   La compatibilité entre l'isomorphisme de  Satake classique et l'équivalence de Satake géométrique s'exprime par le fait que, si $v$ est une place de $X$ de corps résiduel $k(v)$ et $V$ est une représentation irréductible de $\wh G$ de  plus haut poids $ \omega$  et 
 $\rho$ désigne la demi-somme des coracines positives de $\wh G$, 
 $(-1)^{\s{2\rho, \omega}}h_{V,v}$ est égale à la  trace de  $\Frob_{\mr{Gr}_{v}/k(v)} $  (où $\mr{Gr}_{v}$ est la grassmannienne affine en $v$)
 sur le faisceau pervers $\restr{\mc S_{\{0\},V,E}^{(\{0\})}}{\mr{Gr}_{v}}$ 
 (qui est le faisceau d'intersection  de la strate fermée $\mr{Gr}_{v,\omega}$).   Le signe $(-1)^{\s{2\rho, \omega}}$   dans la définition de $h_{V,v}$ ci-dessus 
 vient du fait que la strate $\mr{Gr}_{v,\omega}$ est de dimension $\s{2\rho, \omega}$. Autrement dit $h_{V,v}$ serait la trace de  $\Frob_{\mr{Gr}_{v}/k(v)} $ sur  $\restr{\mc S_{\{0\},V,E}^{(\{0\})}}{\mr{Gr}_{v}}$ 
si on l'avait normalisé pour que   sa cohomologie totale soit supportée en degrés cohomologiques pairs (mais sans changer la torsion à la Tate).
Ce choix  est cohérent avec  le fait que la contrainte de commutativité modifiée de \cite{mv} est celle que l'on obtiendrait naturellement avec de telles normalisations, ainsi qu'on l'a rappelé dans la remarque précédente. La cohérence  de ce choix est la raison pour laquelle   il n'y aura pas de signe   dans l'égalité de la \propref{prop-coal-frob-cas-part-intro}. 
 \end{rem}

   \begin{rem}\label{rem-Satake-Gad}
    Dans le théorème précédent $Z_{\sum_{i\in I} \infty x_{i}}\subset 
G_{\sum_{i\in I} \infty x_{i}}$ agit trivialement sur 
 $ \mr{Gr}_{I}^{(I_{1},...,I_{k})}$ et donc (par d)) sur tous les faisceaux $\mc S_{I,W }^{(I_{1},...,I_{k})}$. On note  $G^{\mr{ad}}=G/Z$. On peut donc considérer $\mc S_{I,W}^{(I_{1},...,I_{k})}$  comme un 
 faisceau pervers  (à un décalage près)
$G^{\mr{ad}}_{\sum \infty x_{i}}$-équivariant sur $\mr{Gr}_{I }^{(I_{1},...,I_{k})}$
 ou si on préfère comme  
 un faisceau pervers 
(à un décalage près)   sur  $\mr{Gr}_{I,W}^{(I_{1},...,I_{k})}/G^{\mr{ad}}_{\sum n_{i}x_{i}}$  (avec les entiers $n_{i}$   assez grands).   
    \end{rem}

Voici la construction du  foncteur \eqref{fonc-W-HIW-intro}. 
Le morphisme \eqref{lisse-chtouca-grass-intro} ne se factorise pas par le quotient par $\Xi$ (comme me l'a fait remarquer un rapporteur anonyme), mais c'est  le cas de sa composée avec le morphisme d'oubli $ \mr{Gr}_{I,W}^{(I_{1},...,I_{k})}/
  G_{\sum n_{i} x_i}\to  \mr{Gr}_{I,W}^{(I_{1},...,I_{k})}/
  G^{\mr{ad}}_{\sum n_{i} x_i}$ (puisque $\Xi$ agit en tordant par des $Z$-torseurs). Autrement dit on possède un morphisme 
  \begin{gather}\label{lisse-chtouca-grass-intro-ad}
    \Cht_{N,I,W}^{(I_{1},...,I_{k})}/\Xi\to \mr{Gr}_{I,W}^{(I_{1},...,I_{k})}/
  G^{\mr{ad}}_{\sum n_{i} x_i} \end{gather} 
et d'après la \remref{rem-Satake-Gad}, $\mc S_{I,W}^{(I_{1},...,I_{k})}$  est un faisceau pervers (à un décalage près) sur l'espace d'arrivée.

\begin{defi}\label{defi-HNIW-can}
On 
 définit le  faisceau pervers (avec la normalisation relative à $(X\sm N)^{I}$) 
 $\mc F_{N,I,W}^{(I_{1},...,I_{k})}$ sur 
 $\Cht_{N,I,W}^{(I_{1},...,I_{k}) }/\Xi$ 
 comme l'image inverse de 
 $\mc S_{I,W}^{(I_{1},...,I_{k})}$ par le morphisme 
   \eqref{lisse-chtouca-grass-intro-ad}.    
   On définit alors le foncteur \eqref{fonc-W-HIW-intro} en posant   
  \begin{gather}\label{defi-can-H}\mc H_{N,I,W}^{\leq\mu}=
   R^{0}(\mf p_{N,I,W} ^{(I_{1},...,I_{k}),\leq\mu})_{!}\Big(\restr{\mc F_{N,I,W}^{(I_{1},...,I_{k})}}{\Cht_{N,I,W}^{(I_{1},...,I_{k}), \leq\mu}/\Xi}\Big)\end{gather} pour toute partition $(I_{1},...,I_{k})$ de $I$. 
   \end{defi}
   Grâce au a) du théorème précédent la définition  \eqref{defi-can-H} ne dépend pas du choix de la partition  $(I_{1},...,I_{k})$. 
 
 Lorsque $W$ est irréductible, la lissité du morphisme \eqref{lisse-chtouca-grass-intro}, et donc celle du morphisme \eqref{lisse-chtouca-grass-intro-ad},  et le calcul de sa dimension impliquent  que 
   $\mc F_{N,I,W}^{(I_{1},...,I_{k})}$ est isomorphe au faisceau  d'intersection de 
   $\Cht_{N,I,W}^{(I_{1},...,I_{k})}/\Xi$ (avec la normalisation perverse relative à $(X\sm N)^{I}$).   Donc la définition précédente est cohérente avec 
   la définition \ref{defi-HNIW-naive} (et la raffine en la rendant plus canonique).

 A l'aide de la définition \eqref{defi-can-H} on peut reformuler 
 l'action des opérateurs de Hecke de fa\c con évidente, ainsi que celle  des morphismes  de Frobenius partiels  de la fa\c con indiquée dans la remarque suivante, que le lecteur peut sauter. 
 
 \begin{rem}\label{rem-Frob-partiel}
 Soit   $W=\boxtimes_{j\in \{1,...,k\}} W_{j}$ où $W_{j}$ est une  représentation de $(\wh G)^{I_{j}}$. 
 Alors l'action des morphismes de Frobenius partiels sur \eqref{defi-can-H} prend encore  la forme d'un morphisme comme \eqref{morph-F1},  et provient, de la même manière que \eqref{morph-F1}, d'un isomorphisme 
 \begin{gather}\label{action-Frob-partiels-general}\Big(\on {Fr}_{I_{1}} ^{(I_{1},...,I_{k})}\Big)^{*}\Big(\mc F_{N,I,W} ^{(I_{2},...,I_{k},I_{1})} \Big)
 \isom \mc F_{N,I,W} ^{(I_{1},...,I_{k})}  \end{gather}
que l'on avait déjà défini dans  \eqref{action-Frob-partiels} pour  $W$  irréductible. Grâce à la définition plus canonique de $\mc F_{N,I,W}^{(I_{1},...,I_{k})}$ donnée dans la \defiref{defi-HNIW-can}, l'isomorphisme \eqref{action-Frob-partiels-general} admet la définition alternative plus canonique suivante 
  
  Grâce à  \defiref{defi-HNIW-can} et à b) du \thmref{satake-geom-thm}, $\mc F_{N,I,W}^{(I_{1},...,I_{k})}$ est isomorphe à l'image inverse  de $\boxtimes _{j\in \{1,...,k\}} \mc S_{I_{j},W_{j}}^{(I_{j})}$ par le morphisme lisse naturel 
 $$  \Cht_{N,I,W}^{(I_{1},...,I_{k}), \leq\mu} /\Xi\to 
  \prod_{j=1}^{k}\Big( \mr{Gr}_{I_{j}}^{(I_{j})}/G^{\mr{ad}}_{\sum_{i\in I_{j}} n_{i} x_{i}}\Big) 
  $$ (où les entiers $n_{i}$ sont asez grands) et en notant 
  $\Frob_{1}$ 
  le morphisme de  Frobenius   de  $\mr{Gr}_{I_{1},W_{1}}^{(I_{1})}/G^{\mr{ad}}_{\sum _{i\in I_{1}}n_{i}x_{i}}$ et 
   \begin{gather}\label{Frob-Gr-1} \on F_{1}:  
 \Frob_{1}^{*}\big( \mc S_{I_{1},W_{1},E}^{(I_{1})}
\big) \isom \mc S_{I_{1},W_{1},E}^{(I_{1})} 
\end{gather} 
l'isomorphisme naturel,   le diagramme 
\begin{gather*}
\begin{CD}
 \Cht_{N,I,W} ^{(I_{1},...,I_{k})}/\Xi @>{\on  {Fr}_{I_{1}, N,I} ^{(I_{1},...,I_{k})}}>> \Cht_{N,I,W} ^{(I_{2},...,I_{1})} /\Xi\\
 @VVV @VVV \\
\prod_{j=1}^{k}
\mr{Gr}_{I_{j},W_{j}}^{(I_{j})}/G^{\mr{ad}}_{\sum _{i\in I_{j}}n_{i}x_{i}}
 @>{\Frob_{1}
\times \Id}>>  \prod_{j=1}^{k} \mr{Gr}_{I_{j},W_{j}}^{(I_{j})}/G^{\mr{ad}}_{\sum _{i\in I_{j}}n_{i}x'_{i}}
 \end{CD}
 \end{gather*}   est commutatif et rend compatibles 
 \eqref{action-Frob-partiels-general} et 
   $$F_{1}\times \Id : (\Frob_{1}\times \Id )^{*}\big(\boxtimes _{j\in \{1,...,k\}} \mc S_{I_{j},W_{j}}^{(I_{j})}
\big) \isom \boxtimes _{j\in \{1,...,k\}} \mc S_{I_{j},W_{j}}^{(I_{j})},  
   $$ ce qui fournit une définition alternative plus canonique de  \eqref{action-Frob-partiels-general} (en particulier le facteur $q^{-d/2}$ vient de $F_{1}$). 
 \end{rem}

 Cependant  la canonicité de la définition \eqref{defi-can-H} est surtout cruciale pour la construction des isomorphismes de coalescence 
 \eqref{intro-isom-coalescence}, que nous allons expliquer maintenant, 
 car les espaces de départ et d'arrivée ne sont pas les mêmes. 
     
 \begin{defi}  L'isomorphisme {\it canonique} \eqref{intro-isom-coalescence} est défini (grâce au théorème de changement de base propre) par l'isomorphisme {\it canonique} 
   entre 
   $\mc F_{N,J,W^{\zeta}}^{(J_{1},...,J_{k})}$   et l'image inverse de 
 $\mc F_{N,I,W}^{(\zeta^{-1}(J_{1}), ..., \zeta^{-1}(J_{k}))}$ par l'inclusion 
   $$\Cht_{N,J }^{(J_{1},...,J_{k})}
     =\Cht_{N,I }^{(\zeta^{-1}(J_{1}), ..., \zeta^{-1}(J_{k}))}\times _{(X\sm N)^{I}}(X\sm N)^{J}\hookrightarrow \Cht_{N,I }^{(\zeta^{-1}(J_{1}), ..., \zeta^{-1}(J_{k}))} 
     $$
     qui provient de l'isomorphisme \eqref{coalescence-Gr-section1-thm} dans le  c) du \thmref{satake-geom-thm}. 
     \end{defi}
     La définition précédente est  indépendante du choix de la partition $(J_{1},...,J_{k})$. Le fait d'avoir choisi cette partition arbitraire permet de  montrer 
    la proposition suivante. 
    
    \begin{prop} 
    \label{rem-coalescence-frob-compat} 
    Les  actions des morphismes de  Frobenius partiels et des opérateurs de Hecke sont  compatibles  avec les morphismes $\mc H(u)$ et avec les isomorphismes de  coalescence de la proposition précédente. 
\end{prop}
  \noindent  La      compatibilité entre les morphismes de Frobenius partiels  et l'isomorphisme de coalescence \eqref{intro-isom-coalescence} signifie 
    que pour tout $j\in J$, $\Delta_{\zeta}^{*}(F_{\zeta^{-1}(\{j\})})$ et 
    $F_{\{j\}}$ se correspondent par l'isomorphisme $\chi_{\zeta}$ 
     de \eqref{intro-isom-coalescence}. 
                
       \section{Démonstration de la \propref{prop-a-b-c} }
       \label{section-esquisse-abc}

 On appelle  point géométrique  $\ov x$ d'un  schéma  $Y$ la donnée d'un corps  algébriquement clos  $k(\ov x)$ et d'un   morphisme $\on{Spec}(k(\ov x))\to Y$.
 Dans cet article $k(\ov x)$ sera toujours une clôture algébrique 
 du corps résiduel $k(x)$ du point $x\in Y$ en-dessous de  $\ov x$. 
  On note  $Y_{(\ov x)}$ le localisé strict  (ou   hensélisé strict) de  $Y$ en   $\ov x$. 
  Autrement dit  $Y_{(\ov x)}$ est le spectre de l'anneau   $\varinjlim \Gamma(U,\mc O_{U})$, où la limite inductive est prise sur les  voisinage étales  $\ov x$-pointés de  $x$. C'est un anneau   local hensélien dont le corps résiduel est la  clôture séparable de  $k(x)$ dans  $k(\ov x)$. 
 Si  $\ov x$ et $\ov y$ sont deux points géométriques de $Y$,  on appelle  flèche de spécialisation $\on{\mf{sp}}:\ov x\to \ov y$   un  morphisme  
 $Y_{(\ov x)}\to Y_{(\ov y)}$, ou de fa\c con équivalente un    morphisme  
 $\ov x\to Y_{(\ov y)}$ (une telle flèche existe si et seulement si $y$ est dans l'adhérence de Zariski de $x$). D'après le paragraphe 7 de \cite{grothendieck-sga4-2-VIII} une   flèche de spécialisation $\on{\mf{sp}}:\ov x\to \ov y$ induit pour tout faisceau    $\mc F$ sur un ouvert de $Y$ contenant $y$ un homomorphisme de spécialisation $\on{\mf{sp}}^{*}: \mc F_{\ov y}\to \mc F_{\ov x}$.

  On fixe une clôture algébrique $\ov F$ de  $F$ et on  note 
      $\ov \eta=\on{Spec}(\ov F)$ le point géométrique correspondant au-dessus du   point  générique
 $\eta$    de $X$.  
 
  On note $\Delta:X\to X^{I}$ le   morphisme diagonal.  
 On fixe un  point géométrique $\ov{\eta^{I}}$ au-dessus du point  générique $\eta^{I}$ de $X^{I}$ et une flèche de spécialisation   
 $\on{\mf{sp}}: \ov{\eta^{I}}\to \Delta(\ov \eta)$. Le rôle de  $\on{\mf{sp}}$  est de rendre le foncteur fibre en  $\ov{\eta^{I}}$ plus canonique, et en particulier compatible avec la  coalescence des pattes (cette dernière affirmation est claire lorsque   $\on{\mf{sp}}^{*}$ est un isomorphisme et en pratique nous serons essentiellement dans cette situation). Donc $ \ov{\eta^{I}}$ et $\on{\mf{sp}}$ vont ensemble et ci-dessous les énoncés faisant intervenir $ \ov{\eta^{I}}$  dépendent du choix de $\on{\mf{sp}}$.

 Un résultat  fondamental de Drinfeld (théorème 2.1 de \cite{drinfeld78} et proposition 6.1 de \cite{drinfeld-compact}) est rappelé  dans le  
 lemme suivant (voir le chapitre 8 de \cite{coh} pour un rappel de la preuve et d'autres références).
  On notera toujours les  $\mc O_{E}$-modules et  les $\mc O_{E}$-faisceaux   par  des lettres gothiques. 
   
 \begin{lem} \label{lem-Dr-intro} (Drinfeld)  Si 
 $\mf E$ est un $\mc O_{E}$-faisceau   lisse   constructible 
 sur un ouvert dense     de  
 $(X\sm N)^{I}$, muni d'une action des  morphismes de  Frobenius partiels, c'est-à-dire  d'isomorphismes
 $$F_{\{i\}}:\restr{\Frob_{\{i\}}^{*}(\mf E)}{\eta^{I}}\to \restr{\mf E}{\eta^{I}}$$
 commutant entre eux et dont la composée est  l'isomorphisme naturel 
    $\Frob^{*}(\mf E)\isom \mf E$, 
alors il s'étend en un faisceau lisse $\wt {\mf E}$ sur  $U^{I}$, où  $U$ est  un ouvert dense assez petit  de  $X\sm N$, et la fibre  $\restr{\wt {\mf E}}{\Delta(\ov \eta)}$ est munie d'une action de   $\pi_{1}(U,\ov\eta)^{I}$. 

De plus, si on fixe  $U$,  le  foncteur
 $\mf E\mapsto \restr{ \mf E }{\Delta(\ov \eta)}$ 
 fournit  une équivalence 
 \begin{itemize}
 \item de la  catégorie des   $\mc O_{E}$-faisceaux
lisses  constructibles  sur   $U^{I}$ munis d'une action des morphismes de  Frobenius partiels 
\item
 vers la  catégorie des 
 représentations continues de   $\pi_{1}(U,\ov\eta)^{I}$ sur des  $\mc O_{E}$-modules de   type fini, 
 \end{itemize}
 de fa\c con  compatible avec   la  coalescence (c'est-à-dire l'image inverse par le morphisme $U^{J}\to U^{I}$ associé à toute application $I\to J$).   
 \end{lem}
       
       \begin{rem} \label{rem-lem-Drinfeld} Dans la situation du début du lemme précédent, $$\on{\mf{sp}}^{*}: 
       \restr{\wt {\mf E}}{\Delta(\ov \eta)}\to \restr{\wt {\mf E}}{\ov{\eta^{I}}}=\restr{\mf E}{\ov{\eta^{I}}}$$ est un isomorphisme, donc $\restr{\mf E}{\ov{\eta^{I}}}$ est muni lui aussi d'une action de $\pi_{1}(U,\ov\eta)^{I}$.        \end{rem}
  \dem   L'idée  pour montrer le \lemref{lem-Dr-intro} est qu'il suffit de traiter le  cas où 
     $\wt {\mf E}=\boxtimes _{i\in I} \mf E_{i}$ avec $ \mf E_{i}$ lisse sur $U$, et dans ce cas l'action de  $\pi_{1}(U,\ov\eta)^{I}$
     sur  $ \restr{\wt {\mf E}}{\Delta(\ov \eta)}=  \otimes _{i\in I} \restr{ \mf E_{i}}{ \ov \eta}$  est évidente.  \cqfd
       
 On ne peut pas appliquer directement ce lemme, car l'action des  morphismes de  Frobenius partiels augmente $\mu$, et d'autre part la limite  inductive  $\varinjlim_{\mu}\mc H_{N,I,W}^{\leq\mu}$ n'est pas constructible  (car ses fibres sont de  dimension infinie). 
 Mais on pourra l'appliquer à la partie ``Hecke-finie'', au sens suivant.

 \begin{defi}
Soit  $\ov x$ un  point géométrique de  $(X\sm N)^{I}$. 
Un élément de  $\varinjlim _{\mu}\restr{\mc H _{N, I, W}^{\leq\mu}}{\ov{x}}$ est dit   Hecke-fini s'il appartient à un  sous-$\mc O_{E}$-module  de type fini de  $\varinjlim _{\mu}\restr{\mc H _{N, I, W}^{\leq\mu}}{\ov{x}}$ qui est 
stable par  $T(f)$ pour tout  $f\in C_{c}(K_{N}\backslash G(\mb A)/K_{N},\mc O_{E})$. 
\end{defi}
On note  $\Big( \varinjlim _{\mu}\restr{\mc H _{N, I, W}^{\leq\mu}}{\ov{x}}\Big)^{\mr{Hf}}$ l'ensemble de tous les  éléments Hecke-finis. 
C'est un sous-$E$-espace vectoriel de   $ \varinjlim _{\mu}\restr{\mc H _{N, I, W}^{\leq\mu}}{\ov{x}}$ et il est stable 
 par  $\pi_{1}(x,\ov{x})$ et  $C_{c}(K_{N}\backslash G(\mb A)/K_{N},E)$. 

\begin{rem}
La définition ci-dessus servira avec  $\ov x$ égal à 
$\Delta(\ov{\eta})$ ou $\ov{\eta^{I}}$. Dans ce cas l'action des opérateurs de Hecke $T(f)$ sur $\varinjlim _{\mu}\restr{\mc H _{N, I, W}^{\leq\mu}}{\ov{x}}$ est évidente (et ne nécessite pas leur extension en des morphismes de faisceaux sur $(X\sm N)^{I}$ qui sera réalisée après la \propref{Eichler-Shimura-intro}). 
\end{rem}
 
 On possède l'homomorphisme de spécialisation
 \begin{gather}\label{sp*-sans-Hf-intro}\on{\mf{sp}}^{*}: 
 \varinjlim _{\mu}\restr{\mc H _{N, I, W}^{\leq\mu}}{\Delta(\ov{\eta})} \to
  \varinjlim _{\mu}\restr{\mc H _{N, I, W}^{\leq\mu}}{\ov{\eta^{I}}}\end{gather} 
 où les deux membres sont considérés comme des $E$-espaces vectoriels
 (limites inductives de $E$-espaces vectoriels de dimension finie). 

On admet maintenant deux résultats, qui seront justifiés dans les 
paragraphes \ref{subsection-crea-annihil-intro}, \ref{section-sous-faisceaux-constr} et \ref{para-homom-spe-Hf}.

 \noindent{\bf Premier résultat admis provisoirement} (\lemref{lem-Hf-union-stab} et \propref{cor-action-Hf-intro}). L'espace   
   $\Big( \varinjlim _{\mu}\restr{\mc H _{N, I, W}^{\leq\mu}}{\ov{\eta^{I}}}\Big)^{\mr{Hf}}$ est une réunion de sous-$\mc O_{E}$-modules de type fini stables par  les morphismes de Frobenius partiels et    le lemme de Drinfeld le munit   donc  d'une action  de $\on{Gal}(\ov F/F)^{I}$, plus précisément $\Big( \varinjlim _{\mu}\restr{\mc H _{N, I, W}^{\leq\mu}}{\ov{\eta^{I}}}\Big)^{\mr{Hf}}$ est une limite inductive de représentations continues de dimension finie de  $\on{Gal}(\ov F/F)^{I}$. 
   
  \noindent{\bf Second résultat admis provisoirement}  (\corref{bijectivite-Hecke-fini}). La restriction de l'homomorphisme $\on{\mf{sp}}^{*}$ de \eqref{sp*-sans-Hf-intro} aux parties Hecke-finies est  un isomorphisme 
   \begin{gather}\label{isom-sp*-avec-Hf-intro} 
  \Big( \varinjlim _{\mu}\restr{\mc H _{N, I, W}^{\leq\mu}}{\Delta(\ov{\eta})}\Big)^{\mr{Hf}} \isor{\on{\mf{sp}}^{*}}
 \Big(  \varinjlim _{\mu}\restr{\mc H _{N, I, W}^{\leq\mu}}{\ov{\eta^{I}}}\Big)^{\mr{Hf}}. \end{gather}

   Ces deux résultats permettent de définir maintenant les $E$-espaces vectoriels $H_{I,W}$  (on omet $N$ dans la notation $ H_{I,W}$ pour limiter la taille des diagrammes dans le paragraphe suivant). 
   
   \begin{defi}
    On définit $H_{I,W}$ comme le membre de gauche de \eqref{isom-sp*-avec-Hf-intro}. 
  \end{defi}
 
 L'action de $\on{Gal}(\ov F/F)^{I}=\pi_{1}(\eta,\ov\eta)^{I}$ sur 
 $H_{I,W}$ ne dépend pas du choix de $\ov{\eta^{I}}$ et $\on{\mf{sp}}$. 
 En effet on peut reformuler ce qui précède en disant que d'après le premier résultat admis    on peut trouver 
 \begin{itemize}
 \item une réunion croissante 
 (indexée par $\lambda\in \N$) 
 de sous-$\mc O_{E}$-faisceaux constructibles 
 $\mf F_{\lambda}\subset \varinjlim _{\mu}\restr{\mc H _{N, I, W}^{\leq\mu}}{\eta^{I}}$ stables par les morphismes de Frobenius partiels
 (auxquels s'applique donc le lemme de Drinfeld) 
 \item une suite décroissante d'ouverts  denses $U_{\lambda}\subset X\sm N$ tels que 
 $\mf F_{\lambda}$ se prolonge en un faisceau lisse sur $(U_{\lambda})^{I}$
 \end{itemize}
 de sorte que $\bigcup _{\lambda\in \N}  \restr{\mf F_{\lambda}}{\ov{\eta^{I}}} 
 =
  \Big(  \varinjlim _{\mu}\restr{\mc H _{N, I, W}^{\leq\mu}}{\ov{\eta^{I}}}\Big)^{\mr{Hf}}$. 
   Alors le second résultat admis    
   implique   que  le morphisme naturel 
  \begin{gather}\label{mor-HIW-Flambda}
  H_{I,W}= \Big( \varinjlim _{\mu}\restr{\mc H _{N, I, W}^{\leq\mu}}{\Delta(\ov{\eta})}\Big)^{\mr{Hf}} \to \bigcup _{\lambda\in \N}  \restr{\mf F_{\lambda}}{\Delta(\ov{\eta})} \end{gather}
(qui vient de la lissité de $\mf F_{\lambda}$ sur $(U_{\lambda})^{I}\ni 
\Delta(\ov{\eta})$)
  est un isomorphisme. Or l'action de $\on{Gal}(\ov F/F)^{I}$   sur le membre de droite de \eqref{mor-HIW-Flambda}, qui est donnée par le lemme de Drinfeld, ne dépend pas  du choix de $\ov{\eta^{I}}$ et $ \mf{sp} $, et donc l'action  de $\on{Gal}(\ov F/F)^{I}$ sur le membre de gauche n'en dépend pas non plus. 
  
 \begin{rem}
 Dans cet article nous montrons seulement que  $H_{I,W}$ est une limite inductive de $E$-espaces vectoriels de dimension finie munis de représentations continues de $\on{Gal}(\ov F/F)^{I}$. En fait Cong Xue a montré dans \cite{these-cong} que 
 $H_{I,W}$ est de dimension finie. La démonstration est difficile et écrite seulement quand $G$ est déployé. 
 \end{rem}
 
   Pour toute application  $\zeta:I\to J$, l'isomorphisme de coalescence 
   \eqref{intro-isom-coalescence} respecte trivialement les  parties Hecke-finies et induit donc  un isomorphisme 
    \begin{gather} \label{isom-chi-zeta-intro}
    H_{I,W}= \Big( \varinjlim _{\mu}\restr{\mc H _{N, I, W}^{\leq\mu}}{\Delta(\ov{\eta})}\Big)^{\mr{Hf}}
    \isom 
    \Big( \varinjlim _{\mu}\restr{\mc H _{N, J,W^{\zeta}}^{\leq\mu}}{\Delta(\ov{\eta})}\Big)^{\mr{Hf}}=H_{J,W^{\zeta}}
    \end{gather}
    où $\Delta$ désigne le morphisme diagonal $X\to X^{I}$ ou $X\to X^{J}$. 
    
    \begin{defi} On définit l'isomorphisme
    $$ \chi_{\zeta}:    H_{I,W} \isom H_{J,W^{\zeta}}$$ apparaissant dans le b) de la \propref{prop-a-b-c} par \eqref{isom-chi-zeta-intro}. 
    \end{defi}
    
    L'isomorphisme $ \chi_{\zeta}$  est $\on{Gal}(\ov F/F)^{J}$-équivariant,  où $\on{Gal}(\ov F/F)^{J}$ agit sur le membre de gauche par le   morphisme diagonal  
\begin{gather}\label{morp-diag-Gal-intro}\on{Gal}(\ov F/F)^{J}\to \on{Gal}(\ov F/F)^{I},  \ (\gamma_{j})_{j\in J}\mapsto (\gamma_{\zeta(i)})_{i\in I}. 
\end{gather}
En effet, si $\Delta_{\zeta}:X^{J}\to X^{I}$ est le morphisme diagonal 
    \eqref{morph-giad-X-intro} et si la suite  $(\mf F_{\lambda})_{\lambda\in \N}$ est comme ci-dessus relativement à $I$ et $W$, alors la suite  $(\Delta_{\zeta}^{*}(\mf F_{\lambda}))_{\lambda\in \N}$ vérifie les mêmes propriétés relativement à $J$ et $W^{\zeta}$, donc 
       $$\chi_{\zeta}: H_{I,W}=\bigcup _{\lambda\in \N}  \restr{\mf F_{\lambda}}{\Delta(\ov{\eta})} =\bigcup _{\lambda\in \N}  \restr{\Delta_{\zeta}^{*}(\mf F_{\lambda})}{\Delta(\ov{\eta})}=H_{J,W^{\zeta}}
       $$ est $\on{Gal}(\ov F/F)^{J}$-équivariant. 
    
\noindent {\bf Démonstration de la \propref{prop-a-b-c}.  }    Les propriétés a) et b) ont déjà été expliquées.  
La propriété c) résulte du fait que la partie Hecke-finie de 
\eqref{I-vide-intro} est formée exactement des formes automorphes cuspidales, c'est-à-dire que 
$$\Big(C_{c} (G(F)\backslash G(\mb A)/K_N \Xi,E)\Big)^{\mr{Hf}}=C_{c}^{\rm{cusp}}(G(F)\backslash G(\mb A)/K_N \Xi,E).$$

\noindent {\bf Preuve de $\supset$. } Toute fonction cuspidale est Hecke-finie  
 car  le $\mc O_{E}$-module $ C_{c}^{\rm{cusp}}(G(F)\backslash G(\mb A)/K_N \Xi,\mc O_{E}) $
  est   de type fini et  stable par tous les opérateurs $T(f)$ pour $f\in C_{c}(K_{N}\backslash G(\mb A)/K_{N},\mc O_{E})$. 
 
 \noindent {\bf Preuve de $\subset$. } On raisonne par l'absurde et on suppose qu'une  fonction Hecke-finie $f$  n'est pas  cuspidale. Il existe alors un sous-groupe parabolique $P\subsetneq G$, de quotient de Levi $M$ et de radical unipotent $U$,  tel  que le terme constant $f_{P}: g\mapsto \int_{U(F)\backslash U(\mb A)}f(ug)$ soit non nul.   Soit $v$ une  place de $X\sm N$. 
       Comme l'anneau des représentations (de dimension finie) de   $\wh M$ est un module de type fini sur  l'anneau des représentations de $\wh G$, $f_{P}$ est également Hecke-finie relativement aux opérateurs de  Hecke pour $M$ en $v$. 
       Ceux-ci comprennent comme cas particuliers les translations par les éléments de $Z_{M}(F_{v})$. On possède une application degré (relativement à 
       $M/Z$), de $U(\mb A)M(F)\backslash G(\mb A)/K_{N}\Xi$ à valeurs dans un $\Z$-module libre de type fini, sur lequel $Z_{M}(F_{v})$ agit par des translations non triviales. Or le support de $f_{P}$ est inclus dans le translaté d'un cône dans 
       ce $\Z$-module libre de type fini, et cela contredit le fait que $f_{P}$ appartient à 
       un espace vectoriel de dimension finie stable par  $Z_{M}(F_{v})$.
   On renvoie à la proposition 8.23 de \cite{coh} pour les détails de la preuve.   \cqfd

    \section{Preuve  du  \thmref{intro-thm-ppal}  à l'aide de la  \propref{prop-a-b-c}} 
    \label{intro-idee-heurist}
    L'idée  se résume ainsi. 
     D'après \eqref{egalite-vide-0-cusp} (que l'on avait déduit de    la  \propref{prop-a-b-c}),     on a  
    $$H_{\{0\},\mbf 1}=C_{c}^{\mr{cusp}}(G(F)\backslash G(\mb A)/K_{N}\Xi,E).$$
   Pour obtenir la décomposition \eqref{intro1-dec-canonique}  il est donc  équivalent de construire (quitte à augmenter $E$) 
   une  décomposition canonique 
     \begin{gather}\label{intro3-dec-canonique}
  H_{\{0\},\mbf 1}=\bigoplus_{\sigma}
 \mf H_{\sigma}.\end{gather}
Cette dernière sera obtenue par  décomposition spectrale  d'une famille commutative d'endomorphismes de $H_{\{0\},\mbf 1}$, appelés opérateurs d'excursion, que nous allons  étudier à l'aide 
  des  propriétés a) et  b) de la \propref{prop-a-b-c}. On commence par rappeler leur construction, déjà donnée dans la \defiref{defi-constr-excursion-intro}.

     Soit  $I$   un ensemble fini et  $W$ une représentation  $E$-linéaire de  
    $(\wh G)^{I}$. 
    On note $\zeta_{I}:I\to \{0\}$ l'application évidente, si bien que 
    $ W^{\zeta_{I}}$ est simplement  $W$ muni de l'action diagonale  de $\wh G$. Soit 
$x: \mbf 1\to W^{\zeta_{I}}$ et  $\xi :  W^{\zeta_{I}}\to  \mbf 1$
des morphismes de  représentations de  $\wh G$ (autrement dit  $x\in W$ et  $\xi\in W^{*}$  sont invariants sous  l'action diagonale de  $\wh G$). Soit   $(\gamma_{i})_{i\in I}\in \on{Gal}(\ov F/F)^{I}$. 
 
 \begin{defi}On définit l'opérateur  
  \begin{gather*}S_{I,W,x,\xi,(\gamma_{i})_{i\in I}}\in 
  \on{End}(H_{\{0\},\mbf 1})  \end{gather*}
    comme la composée 
  \begin{gather}\label{excursion-def-intro}
  H_{\{0\},\mbf  1}\xrightarrow{\mc H(x)}
 H_{\{0\},W^{\zeta_{I}}}\isor{\chi_{\zeta_{I}}^{-1}} 
  H_{I,W}
  \xrightarrow{(\gamma_{i})_{i\in I}}
  H_{I,W} \isor{\chi_{\zeta_{I}}} H_{\{0\},W^{\zeta_{I}}}  
  \xrightarrow{\mc H(\xi)} 
  H_{\{0\},\mbf  1}.
    \end{gather}\end{defi}
    
    Cet opérateur sera appelé ``opérateur d'excursion''. En paraphrasant  \eqref{excursion-def-intro} il est la composée  
    \begin{itemize}
    \item d'un  opérateur de création associé à $x$, dont l'effet est de créer des pattes en le même point (générique) de la courbe, 
    \item d'une action de Galois, 
    qui promène les pattes sur la courbe indépendamment les unes des autres, puis les ramène au même point (générique) de la courbe, 
       \item d'un  opérateur d'annihilation associé à $\xi$, qui annihile les pattes. 
         \end{itemize}

   Le lemme suivant va résulter des propriétés a) et  b) de la \propref{prop-a-b-c}. 
   \begin{lem} Les  opérateurs d'excursion $S_{I,W,x,\xi,(\gamma_{i})_{i\in I}}$ vérifient 
    les propriétés suivantes : 
     \begin{gather}
     \label{SIW-p0}      S_{I,W,x,{}^{t} u(\xi'),(\gamma_i)_{i\in I}}=S_{I,W',u(x),\xi',(\gamma_i)_{i\in I}} 
 \end{gather}
où $u:W\to W'$ est un   morphisme $(\wh G)^{I}$-équivariant et  $x\in W$ et $\xi'\in (W')^{*}$ sont $\wh G$-invariants,  
  \begin{gather}   \label{SIW-p1}
     S_{J,W^{\zeta},x,\xi,(\gamma_j)_{j\in J}}=S_{I,W,x,\xi,(\gamma_{\zeta(i)})_{i\in I}},
     \\
      \label{SIW-p2}
 S_{I_{1}\cup I_{2},W_{1}\boxtimes W_{2},x_{1}\boxtimes x_{2},\xi_{1}\boxtimes \xi_{2},(\gamma^{1}_i)_{i\in I_{1}}\times (\gamma^{2}_i)_{i\in I_{2}}}= S_{I_{1},W_{1},x_{1},\xi_{1},(\gamma^{1}_i)_{i\in I_{1}}}\circ 
S_{I_{2},W_{2},x_{2},\xi_{2},(\gamma^{2}_i)_{i\in I_{2}}}, 
\\ 
\label{SIW-p3} 
S_{I,W,x,\xi,(\gamma_i(\gamma'_i )^{-1}\gamma''_i)_{i\in I}}=
    S_{I\cup I \cup I,W\boxtimes W^{*}\boxtimes W,\delta_{W} \boxtimes x,
    \xi \boxtimes \on{ev}_{W},
    (\gamma_i)_{i\in I} \times (\gamma'_i)_{i\in I} \times (\gamma''_i)_{i\in I}
    }
             \end{gather}
  où la plupart des  notations sont évidentes, 
    $I_{1}\cup I_{2}$ et  $I\cup I\cup I$ désignent des réunions disjointes, et  $\delta_{W}: \mbf 1\to W\otimes W^{*}$ et $\on{ev}_{W}:W^{*} \otimes W \to \mbf 1$ sont les  morphismes naturels.      
  \end{lem}
   
     \noindent {\bf Démonstration de  \eqref{SIW-p0}.}  
    On pose  $x'=u(x)\in W'$ et $\xi={}^{t}u(\xi')\in W^{*}$. 
Le diagramme 
$$   \xymatrix{
     &  H_{\{0\},(W')^{\zeta_{I}}}  
    \ar[r]^{\chi_{\zeta_{I}}^{-1}}&  H_{I,W'} \ar[r]^{(\gamma_{i})_{i\in I}} 
     &
    H_{I,W'} \ar[r]^{\chi_{\zeta_{I}}} &  
    H_{\{0\},(W')^{\zeta_{I}}}  \ar[dr]^{\mc H(\xi')} 
     \\
  H_{\{0\},\mbf  1}  \ar[r]_{\mc H(x)}    \ar[ur]^{\mc H(x')}  &    
    H_{\{0\},W^{\zeta_{I}}}   \ar[u]_{\mc H(u)}
    \ar[r]_{\chi_{\zeta_{I}}^{-1}} & 
    H_{I,W}  \ar[u]_{\mc H(u)} \ar[r]_{(\gamma_{i})_{i\in I}} & 
   H_{I,W}  \ar[u]^{\mc H(u)} \ar[r]_{\chi_{\zeta_{I}}} &
     H_{\{0\},W^{\zeta_{I}}} \ar[u]^{\mc H(u)} \ar[r]_{\mc H(\xi)} & 
     H_{\{0\},\mbf  1}
    }$$           
est commutatif. Or la ligne du bas est égale à  
   $S_{I,W,x,\xi,(\gamma_{i})_{i\in I}}
$ et celle du haut est égale à 
   $S_{I,W',x',\xi',(\gamma_{i})_{i\in I}}
$. \cqfd
       
   \noindent {\bf Démonstration  de  \eqref{SIW-p1}.}  
Le diagramme 
$$   \xymatrix{
     & &  H_{J,W^{\zeta}} \ar[rr]^{(\gamma_{j})_{j\in J}}  &&
    H_{J,W^{\zeta}} \ar[dr]^{\chi_{\zeta_{J}}}
     \\
  H_{\{0\},\mbf  1}  \ar[r]_{\mc H(x)}    &    
    H_{\{0\},W^{\zeta_{I}}}   \ar[ur]^{\chi_{\zeta_{J}}^{-1}} 
    \ar[r]_{\chi_{\zeta_{I}}^{-1}} & 
    H_{I,W}  \ar[u]_{\chi_{\zeta}} \ar[rr]_{(\gamma_{\zeta(i)})_{i\in I}} && 
   H_{I,W}  \ar[u]^{\chi_{\zeta}} \ar[r]_{\chi_{\zeta_{I}}} &
     H_{\{0\},W^{\zeta_{I}}}  \ar[r]_{\mc H(\xi)} & 
     H_{\{0\},\mbf  1}
    }$$           
est commutatif. Or la ligne du bas est égale à 
   $S_{I,W,x,\xi,(\gamma_{\zeta(i)})_{i\in I}}
$ et celle du haut est égale à 
   $    S_{J,W^{\zeta},x,\xi,(\gamma_j)_{j\in J}}$.  \cqfd
       
   \noindent {\bf Démonstration  de  \eqref{SIW-p2}.}  
   L'application évidente 
   $\{0\}\cup \{0\} \to \{0\}$ 
  donne un  isomorphisme $ H_{\{0\}\cup \{0\},\mbf  1}\simeq  H_{\{0\},\mbf  1}$. En notant   $\zeta_{1}:I_{1}\to \{0\}$ et  $\zeta_{2}:I_{2}\to \{0\}$ les applications évidentes,  le membre de gauche de  \eqref{SIW-p2} est égal à la composée 
       \begin{gather*}
 H_{\{0\}\cup \{0\},\mbf  1}
 \xrightarrow{\mc H(x_{1}\boxtimes x_{2})}
 H_{\{0\}\cup \{0\},W_{1}^{\zeta_{1}}\boxtimes W_{2}^{\zeta_{2}}}
 \isor{\chi_{\zeta_{1}\times \zeta_{2}}^{-1}} 
  H_{I_{1}\cup I_{2},W_{1}\boxtimes W_{2}}\\
  \xrightarrow{(\gamma^{1}_{i})_{i\in I_{1}}\times (\gamma^{2}_{i})_{i\in I_{2}}}
  H_{I_{1}\cup I_{2},W_{1}\boxtimes W_{2}} \isor{\chi_{\zeta_{1}\times \zeta_{2}}} 
  H_{\{0\}\cup \{0\},W_{1}^{\zeta_{1}}\boxtimes W_{2}^{\zeta_{2}}}  
  \xrightarrow{\mc H(\xi_{1}\boxtimes \xi_{2})} 
  H_{\{0\}\cup \{0\},\mbf  1}. 
   \end{gather*}
En regroupant  $x_{1},\chi_{\zeta_{1}}^{-1},(\gamma^{1}_{i})_{i\in I_{1}} ,\chi_{\zeta_{1}},\xi_{1}$ d'un côté  et 
$x_{2},\chi_{\zeta_{2}}^{-1},(\gamma^{2}_{i})_{i\in I_{2}} ,\chi_{\zeta_{2}},\xi_{2}$ de l'autre on trouve  le membre de droite. 
On peut le faire car  dans le diagramme  suivant  (où l'on note 
$\gamma^{1}=(\gamma^{1}_{i})_{i\in I_{1}}$ et $\gamma^{2}=(\gamma^{2}_{i})_{i\in I_{2}}$) 
tous les carrés et les triangles commutent. 
\cqfd

   {  \resizebox{14cm}{!}{ $$\!\!\!\!\!\!\!\!\!\!\!\!\!\!\!    \xymatrix{
     {  H_{\{0\}\cup \{0\},\mbf 1} }
           \ar[d]_{{\mc H(x_{1}\boxtimes 1)} }
           \ar[rd]^-{{\mc H(x_1\boxtimes x_2)}}
             \\
      {    H_{\{0\}\cup \{0\},W_{1}^{\zeta_{1}}\boxtimes \mbf 1}}
           \ar[d]_{{\chi_{\zeta_1\times \Id}^{-1}}}
           \ar[r]^<<<{{\mc H(\Id\boxtimes x_{2} )} }
       &  { H_{\{0\}\cup \{0\},W_{1}^{\zeta_{1}}\boxtimes W_{2}^{\zeta_{2}}}}
        \ar[d]_{{\chi_{\zeta_1\times \Id}^{-1}} }
        \ar[rd]^-{{\chi_{\zeta_1\times \zeta_{2}}^{-1}}}
            \\
    {     H_{I_{1}\cup \{0\},W_{1}\boxtimes \mbf 1} }
         \ar[d]_{{\gamma^{1}\times 1}}
         \ar[r]^-{{ \mc H(\Id\boxtimes x_{2} )} }
       &   { H_{I_{1}\cup \{0\},W_{1} \boxtimes W_{2}^{\zeta_{2}}}}
        \ar[d]_{{\gamma^{1}\times 1}}
        \ar[r]^-{{\chi_{\Id\times \zeta_2}^{-1}} }
       &{  H_{I_{1}\cup I_{2},W_{1} \boxtimes W_{2} }  }
         \ar[d]_{{\gamma^{1}\times 1} }
         \ar[rd]^-{ {\gamma^{1}\times \gamma^{2}}}
   \\
   { H_{I_{1}\cup \{0\},W_{1}\boxtimes \mbf 1} }
    \ar[d]_{{\chi_{\zeta_1\times \Id}}}
  \ar[r]^-{{\mc H(\Id\boxtimes x_{2} )} }&  
  { H_{I_{1}\cup \{0\},W_{1} \boxtimes W_{2}^{\zeta_{2}}} }
    \ar[d]_{{\chi_{\zeta_1\times \Id}}}
  \ar[r]^-{{\chi_{\Id\times \zeta_2}^{-1}} }
  & 
 {   H_{I_{1}\cup I_{2},W_{1} \boxtimes W_{2}} }
    \ar[d]_{{\chi_{\zeta_1\times \Id}}}
   \ar[r]^-{{ 1\times  \gamma^{2}}  }
    &{  H_{I_{1}\cup I_{2},W_{1} \boxtimes W_{2}}   }
     \ar[d]_{{\chi_{\zeta_1\times \Id}}}
   \ar[rd]^-{{\chi_{\zeta_1\times \zeta_{2}}}       }
 \\
  {  H_{\{0\}\cup \{0\},W_{1}^{\zeta_{1}}\boxtimes \mbf 1} }
  \ar[d]_{{\mc H(\xi_{1}\boxtimes 1)}  }
  \ar[r]^<<<{{\mc H(\Id\boxtimes x_{2} )} }
  &  { H_{\{0\}\cup \{0\},W_{1}^{\zeta_{1}} \boxtimes W_{2}^{\zeta_{2}}}}
     \ar[d]_{{\mc H(\xi_{1}\boxtimes \Id)} } 
  \ar[r]^-{{\chi_{\Id\times \zeta_2} ^{-1}}}
  & { H_{\{0\}\cup I_{2},W_{1}^{\zeta_{1}} \boxtimes W_{2}}}
    \ar[r]^-{{ 1\times  \gamma^{2}}  }
     \ar[d]_{{\mc H(\xi_{1}\boxtimes \Id)} } &
 {  H_{\{0\}\cup I_{2},W_{1}^{\zeta_{1}} \boxtimes W_{2}}}
      \ar[d]_{{\mc H(\xi_{1}\boxtimes \Id)}  }
       \ar[r]^-{{\chi_{\Id\times \zeta_{2}}}   }
  & {  H_{\{0\}\cup \{0\},W_{1}^{\zeta_{1}} \boxtimes W_{2}^{\zeta_{2}} }}
  \ar[d]_{{\mc H(\xi_{1}\boxtimes \Id)}  }
 \ar[rd]^-{{\mc H(\xi_{1} \boxtimes \xi_2)} }
 \\           
 { H_{\{0\}\cup \{0\},\mbf 1} }
  \ar[r]^-{{ \mc H(1\boxtimes x_{2} )}} &   
 { H_{\{0\}\cup \{0\},\mbf 1  \boxtimes W_{2}^{\zeta_{2}}} }
  \ar[r]^-{{ \chi_{\Id\times \zeta_2} ^{-1}}}&  
{  H_{\{0\}\cup I_{2},\mbf 1  \boxtimes W_{2}} }
   \ar[r]^-{{ 1\times  \gamma^{2}}   } &  
 { H_{\{0\}\cup I_{2},\mbf 1  \boxtimes W_{2}} }
 \ar[r]^-{{\chi_{\Id\times \zeta_{2}}}   }
&{  H_{\{0\}\cup \{0\},\mbf 1  \boxtimes W_{2}^{\zeta_{2}}} }
\ar[r]^-{{\mc H(1\boxtimes \xi_{2} )}} & 
{ H_{\{0\}\cup \{0\},\mbf 1}} }$$}
  }
          
  \vskip1mm        
  \noindent{\bf Démonstration  de   \eqref{SIW-p3}.}  
   Pour tout  $(g_i)_{i\in I}\in (\wh G)^{I}$, 
\begin{itemize}
\item   $ \xi \boxtimes \on{ev}_{W}$     est  invariant par $(1)_{i\in I}\times (g_i)_{i\in I} \times  (g_i)_{i\in I} $
\item $\delta_{W} \boxtimes x$   est  invariant par      $(g_i)_{i\in I} \times (g_i)_{i\in I} \times  (1)_{i\in I} $.
    \end{itemize}    
Donc pour tous   $(\alpha_{i})_{i\in I}$ et $(\beta_{i})_{i\in I}$ dans $\on{Gal}(\ov F/F)^{I}$, le membre de droite de \eqref{SIW-p3} est égal à \begin{gather}\label{SIW-p3-alpha-beta}
 S_{I\cup I \cup I,W\boxtimes W^{*}\boxtimes W,\delta_{W} \boxtimes x,
    \xi \boxtimes \on{ev}_{W},
    (\gamma_i\beta_{i})_{i\in I} \times (\alpha_{i}\gamma'_i\beta_{i})_{i\in I} \times (\alpha_{i}\gamma''_i)_{i\in I}
   . } \end{gather}
  Pour le montrer de fa\c con formelle on factorise  le membre de droite de \eqref{SIW-p3} à travers 
  $$H_{I, \mbf 1} \xrightarrow{ \mc H(\delta_{W})}
     H_{ I,(W\boxtimes W^{*})^{\zeta}} \text{ \  et  \   } H_{ I,(W^{*}\boxtimes W)^{\zeta}} \xrightarrow{ \mc H(\on{ev}_{W})}
    H_{I, \mbf 1},$$ où  $\zeta:I\cup I\to I$ est l'application évidente,   
et on utilise le fait que   $\on{Gal}(\ov F/F)^{I}$ agit trivialement sur  $ H_{I,\mbf 1}\simeq  H_{\emptyset,\mbf 1}$. On prend $\alpha_{i}=\gamma_i (\gamma'_i)^{-1}$ et 
$\beta_{i}=(\gamma'_i)^{-1}\gamma''_i $. Alors   \eqref{SIW-p3-alpha-beta} est égal à  
           \begin{gather}\label{eq-lem-b-2-ggg-intro}
                 S_{I\cup I \cup I,W\boxtimes W^{*}\boxtimes W,\delta_{W} \boxtimes x,\xi \boxtimes \on{ev}_{W},
                 (\gamma_i(\gamma'_i )^{-1}\gamma''_i)_{i\in I} \times 
                 (\gamma_i(\gamma'_i )^{-1}\gamma''_i)_{i\in I} \times
                  (\gamma_i(\gamma'_i )^{-1}\gamma''_i)_{i\in I}
                  }.                   \end{gather}
En appliquant   \eqref{SIW-p1} à  l'application évidente  $\zeta:I\cup I \cup I\to I$, on voit que  \eqref{eq-lem-b-2-ggg-intro} est égal à 
  \begin{gather}\label{eq-lem-b-2-ggg2-intro}S_{I,W\otimes W^{*}\otimes W,\delta_{W} \otimes x,\xi \otimes \on{ev}_{W},(\gamma_i(\gamma'_i )^{-1}\gamma''_i)_{i\in I}}.\end{gather}
Finalement on montre que 
 \eqref{eq-lem-b-2-ggg2-intro} est égal au  membre de gauche de \eqref{SIW-p3}
 en appliquant  \eqref{SIW-p0} 
 à   l'injection  $(\wh G)^{I}$-linéaire $$u:W=\mbf  1\otimes W
    \xrightarrow{  \delta_{W}    \otimes\Id_{W}} W\otimes W^{*}\otimes W, $$
        qui vérifie   $\delta_{W} \otimes x =u( x)$ et  
                  ${}^{t}u(\xi\otimes\on{ev}_{W}) =\xi$, puisque 
      \begin{gather}\label{zorro} \text{la composée \ } W
    \xrightarrow{  \delta_{W}    \otimes\Id_{W}} W\otimes W^{*}\otimes W
     \xrightarrow{ \Id_{W}\otimes\on{ev}_{W}} W 
                \text{ \   est égale à   \ } \Id_{W} \end{gather}
                (c'est un lemme facile  dans la catégorie des espaces vectoriels, appelé parfois lemme de Zorro, et  c'est aussi un des axiomes des catégories tannakiennes).             \cqfd

      On note 
$\mc B \subset \on{End}_{C_{c}(K_{N}\backslash G(\mb A)/K_{N},E)}( H_{\{0\},\mbf 1}) $ 
la sous-algèbre engendrée par tous les  opérateurs d'excursion $S_{I,W,x,\xi,(\gamma_{i})_{i\in I}}$.
En vertu de   \eqref{SIW-p2},  $\mc B$ est commutative. 
   
   Dans la suite de ce paragraphe on considère $\wh G$ comme un schéma en groupes défini sur $E$. 
     
 \begin{observation} \label{observation-f-W-x-xi} Les fonctions \begin{gather}\label{intro-def-f}
   f: (g_{i})_{i\in I}\mapsto \s{\xi, (g_{i})_{i\in I}\cdot x}\end{gather} que l'on obtient en faisant varier  $W$, 
 $x$, et  $\xi $ sont exactement les fonctions régulières sur  le quotient grossier     de  $(\wh G)^{I}$ par  translation à gauche et à droite par 
 $\wh G$ diagonal, que l'on notera $\wh G\backslash (\wh G)^{I}/\wh G$. 
 \end{observation}
 
 \begin{lem}\label{f-W,x,xi}
L'opérateur $S_{I,W,x,\xi,(\gamma_{i})_{i\in I}}$ dépend seulement de   $I$, $f$, et $(\gamma_{i})_{i\in I}$, où $f$ est donnée par \eqref{intro-def-f}. 
\end{lem}
\dem Soit $W,x,\xi$ comme précédemment et soit $f\in \mc O(\wh G\backslash (\wh G)^{I}/\wh G)$ donnée  par \eqref{intro-def-f}. 
 On note $W_{f}$ le sous-$E$-espace vectoriel de dimension finie   de 
 $\mc O((\wh G)^{I}/\wh G)$ engendré les  translatées à gauche de $f$ par $(\wh G)^{I}$. On pose  $x_{f}=f\in W_{f}$ et on note  $\xi_{f}$ la forme linéaire  sur $W_{f}$ donnée par l'évaluation en $1\in (\wh G)^{I}/\wh G$. 
Alors  $W_{f}$ est un  sous-quotient de $W$: si $W_{x}$ est la 
sous-$(\wh G)^{I}$-représentation $E$-linéaire  de $W$ engendrée par $x$, 
 $W_{f}$ est le  quotient de $W_{x}$ par la plus grande sous-$(\wh G)^{I}$-représentation $E$-linéaire sur laquelle $\xi$ s'annule. 
 On a alors les   diagrammes $$W \overset{\alpha}{\hookleftarrow} W_{x} \overset{\beta}{\twoheadrightarrow} W_{f}, 
 \ \ x \overset{\alpha}{\longleftarrow\!\shortmid} x \overset{\beta}{\shortmid\!\longrightarrow} x_{f}, 
 \ \ \xi \overset{{}^{t}\alpha}{\shortmid\!\longrightarrow}  \restr{\xi}{W_{x}} \overset{{}^{t}\beta}{\longleftarrow\!\shortmid} \xi_{f}$$
 de  $(\wh G)^{I}$-représentations, de 
  vecteurs $\wh G$-invariants  et de  formes linéaires  $\wh G$-invariantes. 
  En appliquant  \eqref{SIW-p0} 
  à $u=\alpha$ et $u=\beta$, on obtient 
      \begin{gather}\nonumber 
  S_{I,W,x,\xi,(\gamma_{i})_{i\in I}}=
  S_{I,W_{x},x,\xi |_{W_{x}},(\gamma_{i})_{i\in I}}
  =S_{I,W_{f},x_{f},\xi_{f},(\gamma_{i})_{i\in I}}.\end{gather}
Cela montre que  $S_{I,W,x,\xi,(\gamma_{i})_{i\in I}}$ dépend seulement de  $I$, $f$, et $(\gamma_{i})_{i\in I}$. 
\cqfd

Le lemme précédent permet de définir  la notation suivante (déjà introduite dans la \defiref{defi-constr-excursion-intro}).  
 
 \begin{defi}
 Pour toute fonction  $f\in  \mc O(\wh G\backslash (\wh G)^{I}/\wh G)$ on pose  
 \begin{gather} S_{I,f,(\gamma_{i})_{i\in I}}=S_{I,W,x,\xi,(\gamma_{i})_{i\in I}}\in \mc B \end{gather}
  où $W,x,\xi$ sont tels que  $f$ satisfasse  
 \eqref{intro-def-f}. 
 \end{defi}
 
 Cette nouvelle  notation permet  une formulation   plus synthétique des propriétés     \eqref{SIW-p1}, \eqref{SIW-p2} et \eqref{SIW-p3}, 
 sous la forme des  propriétés (i), (ii), (iii) et (iv) de la \propref{prop-SIf-i-ii-iii}  que l'on peut maintenant démontrer.

      \noindent {\bf Démonstration des propriétés (i), (ii), (iii) et (iv) de la \propref{prop-SIf-i-ii-iii}.}    Le fait que    $S_{I,f,(\gamma_{i})_{i\in I}}$ dépend seulement de  l'image de 
 $(\gamma_{i})_{i\in I}$ dans   $\pi_{1}(X\sm N, \ov\eta)^{I}$ est montré dans le lemme suivant. 
On déduit   (ii) de \eqref{SIW-p1}. Pour montrer   (i) on invoque 
\eqref{SIW-p2} avec $I_{1}=I_{2}=I$,  et on applique   \eqref{SIW-p1} à l'application évidente $\zeta :I\cup I\to I$. La propriété  (iii) résulte de   \eqref{SIW-p3}, après avoir remarqué que le lemme de Zorro \eqref{zorro} implique que  
    $$\s{ \xi \boxtimes \on{ev}_{W}, \big((g_{i})_{i\in I} 
  \boxtimes  (g'_{i})_{i\in I} \boxtimes  (g''_{i})_{i\in I} \big)\cdot ( \delta_{W} \boxtimes x)}
  =\s{\xi, (g_{i}(g'_{i})^{-1}g''_{i})_{i\in I} \cdot x}.$$
  Enfin (iv) résulte du fait que $H_{I,W}$ est une limite inductive de représentations continues de dimension finie de $\on{Gal}(\ov F/F)^{I}$.  \cqfd

Le lemme suivant a été trouvé par    B\"ockle,    Harris,  Khare et  Thorne  et sert   dans  leur article \cite{boeckle-harris...}. 

\begin{lem} \label{prop-harris} 
   Pour tout  $I$ et $f\in \mc O(\wh G\backslash (\wh G)^{I}/\wh G)$,    $S_{I,f,(\gamma_{i})_{i\in I}}$ dépend seulement de  l'image de 
 $(\gamma_{i})_{i\in I}$ dans   $\pi_{1}(X\sm N, \ov\eta)^{I}$, et 
 $(\gamma_{i})_{i\in I}\mapsto S_{I,f,(\gamma_{i})_{i\in I}}$ est continue du  groupe profini $\pi_{1}(X\sm N, \ov\eta)^{I}$ vers la $E$-algèbre de  dimension finie   
   $\mc B$ munie de la topologie   $E$-adique. 
\end{lem}
 \dem    
         Soit  $v$ une  place de $X\sm N$. On fixe  un plongement $\ov F\subset \ov {F_{v}}$, d'où une  inclusion 
            $ \on{Gal}(\ov {F_{v}}/F_{v})\subset  \on{Gal}(\ov F/F)$. Soit  $I_{v}=\on{Ker}( \on{Gal}(\ov {F_{v}}/F_{v})\to \wh \Z)$ le groupe d'inertie en   $v$. Alors 
pour $I, W,x,\xi$ comme dans \eqref{excursion-def-intro}, l'image de  
la composée $  H_{\{0\},\mbf  1}\xrightarrow{\mc H(x)}
 H_{\{0\},W^{\zeta_{I}}}\isor{\chi_{\zeta_{I}}^{-1}} 
  H_{I,W}$ (qui est le début de \eqref{excursion-def-intro}) est formée d'éléments invariants par 
    $(I_{v})^{I}$, car les opérateurs 
    de création sont des morphismes de faisceaux  sur $\Delta(X\sm N)$ tout entier (et en particulier en $\Delta(v)$).    Donc pour $(\gamma_{i})_{i\in I}\in \on{Gal}(\ov F/F)^{I}$ et 
    $(\delta_{i})_{i\in I}\in (I_{v})^{I}$ on a 
    \begin{gather}\label{rel-S-gammai-deltai}S_{I,W,x,\xi,(\gamma_{i})_{i\in I}}=S_{I,W,x,\xi,(\gamma_{i}\delta_{i})_{i\in I}}.\end{gather} 
    Cela est vrai pour tout plongement $\ov F\subset \ov {F_{v}}$ 
    (en fait,    grâce à la \remref{rem-gamma-en-plus},  \eqref{rel-S-gammai-deltai}  pour 
    un plongement implique \eqref{rel-S-gammai-deltai} pour tous les plongements). Or   $\pi_{1}(X\sm N, \ov\eta)$ est le quotient topologique    de 
    $\on{Gal}(\ov F/F)$ par le sous-groupe fermé engendré  
        par les $I_{v}$ pour $v\in (X\sm N)$ et leurs conjugués. \cqfd

 On ne sait pas  si $\mc B$ est réduite. On possède néanmoins  une décomposition spectrale 
  (c'est-à-dire une décomposition en espaces propres généralisés, ou ``espaces caractéristiques'') 
  \begin{gather}\label{dec-nu-intro-intro}
  H_{\{0\},\mbf 1} =\bigoplus_{\nu} \mf H_{\nu}\end{gather}
 où dans le membre de droite la somme directe est indexée  par les  caractères  $\nu$ de  $\mc B$. Quitte à augmenter $E$ on suppose que tous les caractères de $\mc B$ sont définis sur $E$. 
 
  La proposition suivante permet d'obtenir la décomposition \eqref{intro3-dec-canonique}  à partir de \eqref{dec-nu-intro-intro} en associant à chaque caractère $\nu$ un paramètre de Langlands $\sigma$.

 \begin{prop} \label{intro-Xi-n}  Pour tout caractère  $\nu$ de  $\mc B$ il existe un    morphisme 
       $\sigma:\pi_{1}(X\sm N, \ov\eta)\to \wh G(\Qlbar)$
tel que  
\begin{itemize}
\item [] (C1) $\sigma$ prend ses valeurs  dans  $\wh G(E')$, où $E'$ est une extension finie  de $E$, et il est continu, 
\item [] (C2) $\sigma$ est semi-simple, 
 c'est-à-dire que si 
son image est incluse dans un parabolique elle est incluse dans un Levi associé
(comme $\Qlbar$ est de caractéristique $0$ cela équivaut à dire que 
 l'adhérence  de  Zariski de son  image est réductive \cite{bki-serre}),  
\item[] (C3) 
pour tout  $I$ et    $f\in \mc O(\wh G\backslash (\wh G)^{I}/\wh G)$,  on a  
$$\nu(S_{I,f,(\gamma_{i})_{i\in I}})= f\big((\sigma(\gamma_{i}))_{i\in I} \big).
$$
\end{itemize}
De plus   $\sigma$ est unique à conjugaison près par   $\wh G(\Qlbar)$. 
\end{prop}
\noindent

   \noindent {\bf Démonstration.}   On renvoie à la preuve de la proposition 10.7 de \cite{coh} pour quelques détails suplémentaires. 
   La preuve utilise uniquement  les propriétés (i), (ii), (iii) et (iv) de la  \propref{prop-SIf-i-ii-iii}. 
    Soit $\nu$ un caractère de   $\mc B$.

Pour tout $n\in \N$ 
  on note  
   $(\wh G)^{n}\modmod \wh G$ le  quotient grossier de  $(\wh G)^{n}$ par l'action  de  $\wh G$  par  conjugaison diagonale, c'est-à-dire   
   $$h.(g_{1},...,g_{n})=(hg_{1}h^{-1},...,hg_{n}h^{-1}).$$ 
   Alors le morphisme 
   \begin{gather*} (\wh G)^{n}\to (\wh G)^{\{0,...,n\}}, (g_{1},...,g_{n})\mapsto 
   (1,g_{1},...,g_{n})\end{gather*}
   induit  un  isomorphisme  $$\beta: (\wh G)^{n}\modmod \wh G\isom \wh G\backslash (\wh G)^{\{0,...,n\}}/\wh G ,$$ d'où un   isomorphisme  d'algèbres   $$\mc O( (\wh G)^{n}\modmod \wh G)\isom \mc O(\wh G\backslash (\wh G)^{\{0,...,n\}}/\wh G ), \ \ f\mapsto f \circ \beta^{-1}. $$ 

On introduit  \begin{align} \nonumber \Theta_{n}^{\nu}: \mc O((\wh G)^{n}\modmod \wh G)&\to 
C(\pi_{1}(X\sm N, \ov\eta)^{n}, E)\\ \nonumber 
f& \mapsto [(\gamma_{1},...,\gamma_{n})\mapsto 
\nu(S_{I,f \circ \beta^{-1},(1,\gamma_{1},...,\gamma_{n})})]\end{align}

La condition (C3) que doit vérifier $\sigma$ se reformule de la fa\c con suivante : pour tout $n$ et pour tout $f\in \mc O((\wh G)^{n}\modmod \wh G)$, 
\begin{gather}\label{rel-sigma-nu}\Theta_{n}^{\nu}(f)=[(\gamma_{1},...,\gamma_{n})\mapsto f((\sigma(\gamma_{1}),...,\sigma(\gamma_{n})))   ]. \end{gather}

On déduit immédiatement des propriétés (i), (ii), (iii) et (iv) de la  \propref{prop-SIf-i-ii-iii} que la suite 
$(\Theta_{n}^{\nu})_{n\in \N^{*}}$ 
vérifie les propriétés suivantes
\begin{itemize} \item 
pour tout  $n$, $\Theta_{n}^{\nu}$ est un  morphisme d'algèbres, 
                                                 \item  la suite  $(\Theta_{n}^{\nu})_{n\in \N^{*}}$ est fonctorielle par rapport à  toutes les applications entre les ensembles  $\{1,...,n\}$, c'est-à-dire  que  pour $m,n\in \N^{*}$, 
               $$\zeta: \{1,...,m\}\to \{1,...,n\}$$ arbitraire,   
               $f\in \mc O((\wh G)^{m}\modmod \wh G)$ et 
               $(\gamma_{1},...,\gamma_{n })\in \pi_{1}(X\sm N, \ov\eta)^{n }$, 
               on a  
             $$\Theta_{n}^{\nu}( f^{\zeta}) ((\gamma_{j})_{j\in \{1,...,n\}})=
             \Theta_{m}^{\nu}(f)((\gamma_{\zeta(i)})_{i\in \{1,...,m\}})$$ 
            où  $f^{\zeta}\in \mc O((\wh G)^{n}\modmod \wh G)$  est définie par 
               $$f^{\zeta}((g_{j})_{j\in \{1,...,n\}})=f((g_{\zeta(i)})_{i\in \{1,...,m\}}), $$
 \item     pour $n\geq 1$, 
  $f\in \mc O((\wh G)^{n}\modmod \wh G)$  
    et   $(\gamma_{1},...,\gamma_{n+1})\in \pi_{1}(X\sm N, \ov\eta)^{n+1}$ on a  
   $$\Theta_{n+1}^{\nu}( \wh f)(\gamma_{1},...,\gamma_{n+1})=
   \Theta_{n}^{\nu}( f)(\gamma_{1},...,\gamma_{n}\gamma_{n+1}) $$
 où $\wh f\in  \mc O((\wh G)^{n+1}\modmod \wh G)$ est définie par 
   $$\wh f(g_{1},...,g_{n+1})=f(g_{1},...,g_{n}g_{n+1}). $$
\end{itemize}
 Pour justifier la dernière propriété, on applique la  propriété 
   (iii) de la  \propref{prop-SIf-i-ii-iii} à $$I=\{0,...,n \}, 
   (\gamma_{i})_{i\in I}=(1,\gamma_{1},...,\gamma_{n}), (\gamma'_{i})_{i\in I}=(1)_{i\in I}, (\gamma''_{i})_{i\in I}=(1, ...,1,\gamma_{n+1})$$
     et on utilise   (ii) pour supprimer tous les   $1$ sauf le premier dans  
   $(\gamma_{i})_{i\in I}\times (\gamma'_{i})_{i\in I}\times (\gamma''_{i})_{i\in I}$.

On va montrer que ces propriétés de 
 la suite 
$(\Theta_{n}^{\nu})_{n\in \N^{*}}$ 
entraînent l'existence et l'unicité de   $\sigma$ vérifiant (C1), (C2)  et (C3) (c'est-à-dire
\eqref{rel-sigma-nu}). 

Pour $G=GL_{r}$ le résultat est déjà connu:  la suite $(\Theta_{n}^{\nu})_{n\in \N^{*}}$  
est déterminée par $\Theta_{1}^{\nu}(\on{Tr})$ 
(qui doit être le caractère de $\sigma$) 
et 
$\Lambda^{r+1}\mr{St}=0$ implique la relation de pseudo-caractère 
d'où l'existence de $\sigma$ par     \cite{taylor}.  On renvoie à la remarque 11.8 de \cite{coh} pour plus de détails. 

En général on utilise des résultats de   \cite{richardson}. 
On dit qu'un  $n$-uplet   $(g_{1},...,g_{n})\in
 \wh G(\Qlbar)^{n}$ est semi-simple si 
 tout parabolique le contenant 
 possède un sous-groupe de Levi associé le contenant. Comme $  \Qlbar$ est de caractéristique $0$ cela équivaut (d'après \cite{bki-serre}) à la condition que 
 l'adhérence de   Zariski  
  $\ov{<g_{1},...,g_{n}>}$ du sous-groupe  $<g_{1},...,g_{n}>$ engendré par   $g_{1},...,g_{n}$ est réductive. 
 D'après le théorème 3.6 de  \cite{richardson} la $\wh G$-orbite (par conjugaison) de $(g_{1},...,g_{n})$ est fermée dans 
 $(\wh G)^{n}$ si et seulement si $(g_{1},...,g_{n})$ est semi-simple. 
 Donc les points sur $\Qlbar$ du quotient  grossier $(\wh G)^{n}\modmod \wh G$ (qui correspondent aux 
    $\wh G$-orbites fermées définies sur    $\Qlbar$ dans   $(\wh G)^{n}$)  sont   en bijection avec  les classes de conjugaison par  $\wh G(\Qlbar)$  de   $n$-uplets  semi-simples $(g_{1},...,g_{n})\in
 \wh G(\Qlbar)^{n}$.

 Pour tout $n$-uplet 
   $(\gamma_{1},...,\gamma_{n})\in \pi_{1}(X\sm N, \ov\eta)^{n}$ on note  
  $\xi_{n} (\gamma_{1},...,\gamma_{n})$ le point défini sur $\Qlbar$ du 
    quotient  grossier $(\wh G)^{n}\modmod \wh G$ associé au caractère 
   $$\mc O((\wh G)^{n}\modmod \wh G)\to \Qlbar,  \ f\mapsto \Theta_{n}^{\nu}( f)(\gamma_{1},...,\gamma_{n}).$$ 
   On note 
 $\xi_{n}^{\mr{ss}}(\gamma_{1},...,\gamma_{n})$ la  classe de conjugaison de   
  $n$-uplets  semi-simples correspondant à  $\xi_{n}(\gamma_{1},...,\gamma_{n})$ par le résultat de  \cite{richardson} rappelé ci-dessus.  

La relation  \eqref{rel-sigma-nu} équivaut à la condition que pour tout $n$ et pour tout 
$(\gamma_{1},...,\gamma_{n})$, 
$ (\sigma(\gamma_{1}),...,\sigma(\gamma_{n}))\in (\wh G(\Qlbar))^{n}$ 
(qui n'est pas en général semi-simple)
est au-dessus de 
$\xi_{n}(\gamma_{1},...,\gamma_{n})$. 

\noindent {\bf Unicité de $\sigma$ (à conjugaison près).}  On choisit $n$ et $(\gamma_{1},...,\gamma_{n})$ tels que $\sigma(\gamma_{1}),...,\sigma(\gamma_{n})$ engendrent un sous-groupe Zariski dense dans  $\ov{ \mr{Im}(\sigma)}$. 
Comme $\sigma$ est supposé semi-simple, $(\sigma(\gamma_{1}),...,\sigma(\gamma_{n}))$ est semi-simple. 
On fixe 
$(g_{1},...,g_{n})$ dans $\xi_{n}^{\mr{ss}}(\gamma_{1},...,\gamma_{n})$. 
Donc $(\sigma(\gamma_{1}),...,\sigma(\gamma_{n}))$ est conjugué à 
$(g_{1},...,g_{n})$ et quitte à conjuguer $\sigma$ on peut supposer qu'il lui est égal. Alors $\sigma$ est déterminé de fa\c con unique car pour tout $\gamma$, 
$\sigma(\gamma)$  appartient à l'adhérence de Zariski du sous-groupe engendré par $(g_{1},...,g_{n})$ et 
$(g_{1},...,g_{n},\sigma(\gamma)) \in 
\xi_{n+1}^{\mr{ss}}(\gamma_{1},...,\gamma_{n},\gamma)$, donc 
la connaissance de $\xi_{n+1} (\gamma_{1},...,\gamma_{n},\gamma)$ détermine uniquement $\sigma(\gamma)$. 

\noindent {\bf Existence de $\sigma$.}  
  Pour tout $n$ et tout $(\gamma_{1},...,\gamma_{n})\in \pi_{1}(X\sm N, \ov\eta)^{n}$ on choisit $(g_{1},...,g_{n})\in \xi_{n}^{\mr{ss}}(\gamma_{1},...,\gamma_{n})$ (bien défini à conjugaison près).  On choisit alors $n$ et 
   $(\gamma_{1},...,\gamma_{n})\in \pi_{1}(X\sm N, \ov\eta)^{n}$ tels que
       \begin{itemize}
   \item    (H1)
  la  dimension de 
  $\ov{<g_{1},...,g_{n}>}$ est la plus grande  possible 
  \item     (H2) le centralisateur   $C(g_{1},...,g_{n})$ de   $<g_{1},...,g_{n}>$ est le plus petit  possible (dimension minimale puis  nombre de composantes connexes minimal).  
  \end{itemize}
  On fixe  $(g_{1},...,g_{n})\in \xi_{n}^{\mr{ss}}(\gamma_{1},...,\gamma_{n})$ pour le reste de la démonstration et on construit une application   $$\sigma: \pi_{1}(X\sm N, \ov\eta)\to \wh G(\Qlbar)$$  en demandant que pour tout  $\gamma\in \pi_{1}(X\sm N, \ov\eta)$, $\sigma(\gamma)$ est l'unique élément $g$ de  $\wh G(\Qlbar)$ tel que  $(g_{1},...,g_{n},g)\in \xi_{n+1}^{\mr{ss}}(\gamma_{1},...,\gamma_{n},\gamma)$.  
  L'existence et l'unicité de $g$ sont justifiées de la fa\c con suivante. 
  \begin{itemize}
  \item  {\bf A) Existence de $g$} : pour $(h_{1},...,h_{n},h)\in \xi_{n+1}^{\mr{ss}}(\gamma_{1},...,\gamma_{n},\gamma)$, $(h_{1},...,h_{n})$ est forcément semi-simple    (car $(h_{1},...,h_{n})$ est au-dessus de $\xi_{n}(\gamma_{1},...,\gamma_{n})$ et $(g_{1},...,g_{n})\in \xi_{n}^{\mr{ss}}(\gamma_{1},...,\gamma_{n})$
  donc d'après le théorème 5.2 de \cite{richardson}, $\ov{<h_{1},...,h_{n}>}$ admet un sous-groupe de Levi isomorphe à $\ov{<g_{1},...,g_{n}>}$, or 
  $\dim(\ov{<h_{1},...,h_{n}>})\leq \dim(\ov{<g_{1},...,g_{n}>})$
  par (H1)), 
  donc quitte à conjuguer $(h_{1},...,h_{n},h)$ on peut supposer que $(h_{1},...,h_{n})=(g_{1},...,g_{n})$ et on prend alors $g=h$.   
  \item {\bf B) Unicité de $g$} : on a $C(g_{1},...,g_{n},g)\subset C(g_{1},...,g_{n})$ et l'égalité a  lieu par (H2), donc $g$ commute avec $C(g_{1},...,g_{n})$ et comme il était bien déterminé modulo conjugaison par $C(g_{1},...,g_{n})$ il est unique. 
  \end{itemize}
  
Puis on montre que l'application  $\sigma$ que l'on vient de construire est un morphisme de groupes. 
En effet soient $\gamma, \gamma'\in \pi_{1}(X\sm N, \ov\eta)$. Le même argument que dans A) ci-dessus montre qu'il existe $g,g'$ tels que 
\begin{gather}\label{semi-simple-g-g'}(g_{1},...,g_{n},g,g') \in \xi_{n+2}^{\mr{ss}}(\gamma_{1},...,\gamma_{n},\gamma,\gamma').\end{gather}
   Grâce aux  propriétés vérifiées par la suite 
      $(\Theta_{n}^{\nu})_{n\in \N^{*}}$ on voit que 
      $\xi_{n+1}(\gamma_{1},...,\gamma_{n},\gamma\gamma')$ est l'image de  $\xi_{n+2}(\gamma_{1},...,\gamma_{n},\gamma,\gamma')$ par le morphisme 
   $$(\wh G)^{n+2}\modmod \wh G \to (\wh G)^{n+1}\modmod \wh G, 
   (h_{1},...,h_{n},h,h')\mapsto (h_{1},...,h_{n},hh').$$
       On en déduit que 
        $   (g_{1},...,g_{n},gg')$ est au-dessus de 
        $  \xi_{n+1} (\gamma_{1},...,\gamma_{n},\gamma\gamma')$. De plus 
        $ (g_{1},...,g_{n},gg')$ est semi-simple par le même argument que dans A), car 
       \begin{gather*}  \dim(\ov{< g_{1},...,g_{n},gg' >}) 
        \leq  \dim(\ov{< g_{1},...,g_{n},g,g' >}) 
         \leq  \dim(\ov{< g_{1},...,g_{n}  >}) \end{gather*}
          (où la dernière inégalité vient de (H1)). Donc 
          $ (g_{1},...,g_{n},gg')$ appartient à  
        $  \xi_{n+1}^{\mr{ss}} (\gamma_{1},...,\gamma_{n},\gamma\gamma')$
          et 
           $gg'=\sigma(\gamma \gamma')$. 
        Les mêmes arguments montrent que 
     \begin{gather}\label{gg'-gamma-gamma'}g=\sigma(\gamma)\text{ \ \   et \ \ }  g'=\sigma(\gamma').\end{gather} 
On a finalement montré que $\sigma(\gamma \gamma') =\sigma(\gamma)\sigma(\gamma')$. 
  
  Donc $\sigma$ est un morphisme de groupes à valeurs dans $\wh G(E')$
  (où $E'$ est une extension finie de $E$ telle que $g_{1},...,g_{n}$ appartiennent à $\wh G(E')$).    
L'argument pour montrer que $\sigma$ est continu est le suivant. 
On sait que  pour toute fonction $f$ sur $(\wh G_{E'})^{n+1}\modmod \wh G_{E'}$, l'application 
$$
\pi_{1}(X\sm N, \ov\eta)\to E',  \ \ 
\gamma\mapsto f(g_{1},...,g_{n}, \sigma(\gamma))= 
\Theta_{n+1}^{\nu}( f)(\gamma_{1},...,\gamma_{n},\gamma)$$
est continue. 
Or le morphisme 
 \begin{align*}\mc O((\wh G_{E'})^{n+1}\modmod \wh G_{E'}) &\to  \mc O(\wh G_{E'}\modmod C(g_{1},...,g_{n})) \\ f & \mapsto  [g\mapsto f(g_{1},...,g_{n},g)]
\end{align*} est surjectif
(parce que $(g_{1},...,g_{n})$ est semi-simple, donc son orbite par conjugaison est une sous-variété affine fermée de 
$(\wh G_{E'})^{n}$, 
isomorphe à $\wh G_{E'}/C(g_{1},...,g_{n})$). 
De plus, en notant 
$D(g_{1},...,g_{n})$ le centralisateur de $C(g_{1},...,g_{n})$ (qui contient l'image de $\sigma$), le morphisme de restriction 
$$\mc O(\wh G_{E'}\modmod C(g_{1},...,g_{n})) =
\mc O(\wh G_{E'} )^{C(g_{1},...,g_{n})} \to \mc O(D(g_{1},...,g_{n}))$$ est surjectif  parce que la   restriction $\mc O(\wh G_{E'}) \to \mc O(D(g_{1},...,g_{n}))$  est  évidemment   surjective, 
qu'elle le reste   lorsqu'on prend  les invariants par le groupe réductif   $C(g_{1},...,g_{n})$ (agissant par conjugaison), et que 
$C(g_{1},...,g_{n})$ agit trivialement sur  $\mc O(D(g_{1},...,g_{n}))$.  
Donc pour toute fonction $h\in \mc O(D(g_{1},...,g_{n}))$, l'application 
$$
\pi_{1}(X\sm N, \ov\eta)\to E',  \ \ \gamma\mapsto h(\sigma(\gamma))$$ est continue, et on a montré que $\sigma$ est continu. 
 
 Il reste à montrer \eqref{rel-sigma-nu}, c'est-à-dire que pour 
   $m\in \N^{*}$, $f\in \mc O((\wh G)^{m}\modmod \wh G)$ et  $(\delta_{1},...,\delta_{m})\in \pi_{1}(X\sm N, \ov\eta)^{m}$, on a  
\begin{gather}\nonumber f(\sigma(\delta_{1}),...,\sigma(\delta_{m}))=
 \big(\Theta^{\nu}_{m}(f)\big)(\delta_{1},...,\delta_{m}).\end{gather}
Par les mêmes  arguments que  pour \eqref{semi-simple-g-g'} et \eqref{gg'-gamma-gamma'} on montre que  
 $$(g_{1},...,g_{n},\sigma(\delta_{1}),...,\sigma(\delta_{m}))\in \xi_{n+m}^{ss}(\gamma_{1},...,\gamma_{n},\delta_{1},...,\delta_{m} ). $$
 Donc 
 $(\sigma(\delta_{1}),...,\sigma(\delta_{m}))$ est au-dessus de  $\xi_{m}(\delta_{1},...,\delta_{m} )$. \cqfd

Donc on a obtenu la décomposition \eqref{intro1-dec-canonique}. Ceci achève la démonstration du \thmref{intro-thm-ppal}, à condition d'admettre le résultat suivant, qui sera justifié  dans le paragraphe \ref {subsection-intro-decomp}. 

 \noindent{\bf Résultat admis provisoirement} (\propref{S-non-ram-concl-intro}). La décomposition \eqref{intro1-dec-canonique}  est compatible avec l'isomorphisme de Satake en toutes les places de $X\sm N$.

            \section{Morphismes de création et d'annihilation}
           \label{subsection-crea-annihil-intro}
             A partir de maintenant et jusqu'au paragraphe \ref{subsection-intro-decomp}  on va justifier 
             \begin{itemize}
             \item   les deux résultats admis provisoirement  dans  le paragraphe \ref{section-esquisse-abc}
 et qui avaient servi à     la preuve  de la \propref{prop-a-b-c}
       \item et   le résultat admis provisoirement à la fin du paragraphe précédent et qui avait permis de terminer la preuve du  \thmref{intro-thm-ppal}. 
                     \end{itemize}

    Dans ce paragraphe notre but est d'utiliser les   isomorphismes de coalescence  \eqref{intro-isom-coalescence} pour construire des morphismes de création et d'annihilation, puis d'exprimer les   opérateurs de Hecke en les  places de  $X\sm N$ comme la composée   
 \begin{itemize}
 \item d'un morphisme de création,
  \item de  l'action d'un  morphisme de Frobenius partiel,
  \item d'un  morphisme d'annihilation,
  \end{itemize}
et d'utiliser cela pour étendre les  opérateurs de Hecke \eqref{defi-Tf} en des morphismes de faisceaux sur   $(X\sm N)^{I}$ tout entier  et pour obtenir les relations d'Eichler-Shimura.

  Soient   $I$ et $J$ des ensembles finis. 
  On va définir maintenant les morphismes de création et d'annihilation, dont l'idée est la suivante. Les pattes indexées par $I$ restent inchangées et on crée (ou on annihile) les pattes indexées par $J$ en un même point de la courbe (indexé par un ensemble à un élément, que l'on note $\{0\}$).

  On a des applications évidentes 
   $$\zeta_{J} :J \to \{0\},  \ 
  \zeta_{J}^{I}=(\Id_{I},\zeta_{J}): I\cup J\to I\cup \{0\} \text{ \   et   \ } \zeta_{\emptyset}^{I}=(\Id_{I},\zeta_{\emptyset}): I \to I \cup \{0\}.$$ 
  Soient $W$ et $U$ des  représentations  $E$-linéaires de dimension finie de 
$(\wh G)^{I}$ et   $(\wh G)^{J}$ respectivement. 
On rappelle que   $U^{\zeta_{J}}$ est 
la représentation de $\wh G$ obtenue en restreignant 
  $U$ à la diagonale   $\wh G\subset (\wh G)^{J}$. 
Soient  $x\in U$ et $\xi\in U^{*}$ invariants sous  l'action diagonale de $\wh G$.  
Alors $W\boxtimes U$ est une représentation de  $(\wh G)^{I\cup J}$ et 
 $W\boxtimes U^{\zeta_{J}}$ et  $W\boxtimes \mbf 1$ sont des  représentations de  $(\wh G)^{I\cup \{0\}}$ reliées par les   morphismes 
 $$W\boxtimes \mbf 1\xrightarrow{\Id_{W}\boxtimes x}W\boxtimes U^{\zeta_{J}}\text{ \  et \ }  W\boxtimes U^{\zeta_{J}}\xrightarrow{\Id_{W}\boxtimes \xi}W\boxtimes \mbf 1.$$ 
 
 On note  $\Delta:X\to X^{J}$ le   morphisme diagonal et on désigne par $E_{X\sm N}$   le faisceau constant sur  $X\sm N$. 
         \begin{defi}
 On définit  le  morphisme de création
  $
\mc C_{x}^{\sharp}$ comme la  composée   \begin{gather*}
  \mc H _{N,I,W}^{\leq\mu} \boxtimes E_{(X\sm N)} 
  \isor{\chi_{\zeta_{\emptyset}^{I}}}
 \mc H _{N,I\cup\{0\},W\boxtimes \mbf 1}^{\leq\mu}  \\
 \xrightarrow{\mc H(\Id_{W}\boxtimes x) }
 \mc H _{N,I\cup\{0\},W\boxtimes U^{\zeta_{J}}}^{\leq\mu} 
 \isor{ \chi_{\zeta_{J}^{I}}^{-1} }
 \restr{ \mc H _{N,I\cup J,W\boxtimes U}^{\leq\mu}}{(X\sm N)^{I}\times \Delta(X\sm N)}  
   \end{gather*}
   où $\chi_{\zeta_{\emptyset}^{I}}$ et $\chi_{\zeta_{J}^{I}}$ sont  les   isomorphismes  de coalescence de \eqref{intro-isom-coalescence}. 
   De même on définit le  morphisme d'annihilation $
\mc C_{\xi}^{\flat}$  comme la  composée 
  \begin{gather*}
  \restr{ \mc H _{N,I\cup J,W\boxtimes U}^{\leq\mu}}{(X\sm N)^{I}\times \Delta(X\sm N)} 
  \isor{ \chi_{\zeta_{J}^{I}}}
 \mc H _{N,I\cup\{0\},W\boxtimes U^{\zeta_{J}}}^{\leq\mu} \\
 \xrightarrow{ \mc H(\Id_{W}\boxtimes \xi) }
  \mc H _{N,I\cup\{0\},W\boxtimes \mbf 1}^{\leq\mu} 
 \isor{\chi_{\zeta_{\emptyset}^{I}}^{-1} }
   \mc H _{N,I,W}^{\leq\mu} \boxtimes E_{(X\sm N)} . 
    \end{gather*}
    \end{defi}
    
   Tous les  morphismes ci-dessus sont des  morphismes de faisceaux sur  
    $(X\sm N)^{I} \times (X\sm N)$.

 Nous allons maintenant utiliser ces morphismes avec $J=\{1,2\}$. 
     Soit  $v$  une  place dans   $|X|\sm |N|$.    On considère  $v$ également  comme un  sous-schéma de  $X$ et on note   $E_{v}$ le faisceau constant sur   $v$.    Soit $V$ une représentation irréductible de $\wh G$.    Comme précédemment on note   $\mbf 1\xrightarrow{\delta_{V}} V\otimes V^{*}$ et $ V\otimes V^{*}\xrightarrow{\on{ev}_{V}} \mbf 1$ les  morphismes naturels. 
        
      Pour  $\kappa$ assez grand  (en fonction de $\deg(v),V$), on définit   $S_{V,v}$ comme la  composée 
  \begin{gather}\label{def-SVv-intro1}
 \mc H _{N,I,W}^{\leq\mu} \boxtimes E_{v} \\ \label{def-SVv-intro2}
  \xrightarrow{ \restr{\mc C_{
    \delta_{V}}^{\sharp}}{(X\sm N)^{I}\times v}}
 \restr{ \mc H _{N,I\cup\{1,2\},W\boxtimes V\boxtimes V^{*}}^{\leq\mu}}{(X\sm N)^{I}\times \Delta(v)} \\ \label{def-SVv-intro3}
 \xrightarrow{ \restr{(F_{\{1\}})^{\deg(v)} }{(X\sm N)^{I}\times \Delta(v)}}
 \restr{ \mc H _{N,I\cup\{1,2\},W\boxtimes V\boxtimes V^{*}}^{\leq\mu+\kappa}}{(X\sm N)^{I}\times \Delta(v)} \\ \label{def-SVv-intro4}
  \xrightarrow{ \restr{\mc C_{
    \on{ev_{V}}}^{\flat }}{(X\sm N)^{I}\times v}}
     \mc H _{N,I,W}^{\leq\mu+\kappa}  \boxtimes E_{v} . 
 \end{gather}
Autrement dit on crée deux nouvelles pattes en  $v$ à l'aide de  
 $ \delta_{V}:\mbf 1\to V\otimes V^{*}$, on applique le  morphisme de Frobenius partiel (à la puissance  $\deg(v)$) à la première, puis on les annihile à l'aide de  $ \on{ev_{V}}: V\otimes V^{*}\to\mbf 1$. 
  
   En tant que  morphisme de faisceaux constructibles sur 
     $(X\sm N)^{I}\times v$,  $S_{V,v}$ commute avec l'action naturelle du morphisme de   Frobenius partiel sur   $E_{v}$ dans  \eqref{def-SVv-intro1} et  \eqref{def-SVv-intro4}, puisque 
     \begin{itemize}
     \item les  morphismes de création et  d'annihilation  entrelacent  cette action avec l'action  de $F_{\{1,2\}}$ sur  \eqref{def-SVv-intro2} et \eqref{def-SVv-intro3}, par la \propref{rem-coalescence-frob-compat}, 
          \item  $F_{\{1\}}$ et donc   $F_{\{1\}}^{\deg(v)}$ commutent avec  $F_{\{1,2\}}=F_{\{1\}}F_{\{2\}}$. 
     \end{itemize}
  \begin{defi} Par abus  on note encore     
        $$S_{V,v}: \mc H _{N,I,W}^{\leq\mu}\to 
         \mc H _{N,I,W}^{\leq\mu+\kappa}$$  le  morphisme 
    de faisceaux sur   $(X\sm N)^{I}$ obtenu par descente relativement à l'action de  $\Z/\deg(v)\Z$  (en prenant les invariants par  l'action naturelle du morphisme de   Frobenius partiel sur   $E_{v}$ dans  \eqref{def-SVv-intro1} et  \eqref{def-SVv-intro4}). 
     \end{defi}
  
 \begin{prop}\label{prop-coal-frob-cas-part-intro} 
  La  restriction de  $S_{V,v}$  à  $(X\sm (N \cup v))^{I}$
est égale, en tant que   morphisme de faisceaux sur  $(X\sm (N \cup v))^{I}$,  à l'opérateur de Hecke $$T(h_{V,v}):  \restr{\mc H _{N,I,W}^{\leq\mu}}{(X\sm (N \cup v))^{I}}\to 
       \restr{  \mc H _{N,I,W}^{\leq\mu+\kappa}}{(X\sm (N \cup v))^{I}}
       . $$
       \end{prop}
\noindent
  Il suffit de le démontrer lorsque   $V$  et $W$ sont irréductibles. La preuve est de nature  géométrique.  
  Nous l'esquissons ici dans un cas simple, où elle se réduit à l'intersection de deux sous-champs lisses  dans un champ de Deligne-Mumford lisse et où cette intersection s'avère être transverse. La preuve est plus compliquée en général à cause des singularités. On renvoie à la preuve de la proposition 6.2 de \cite{coh} pour le cas général (mais une solution alternative consisterait à se ramener à une situation d'intersection transverse lisse à l'aide des résolutions de Bott-Samelson). 
  
\noindent {\bf Démonstration lorsque  $V$ est minuscule  et  $\deg(v)=1$.} 
On rappelle qu'une représentation irréductible de $\wh G$ est dite minuscule si tous ses poids sont conjugués par le groupe de Weyl. Cela équivaut au fait que l'orbite correspondante dans la grassmannienne affine est fermée  (et implique donc que la strate fermée correspondante  est lisse). On note $d$ la dimension de l'orbite associée à $V$. 

 Grâce à l'hypothèse que  $\deg(v)=1$ on peut supprimer   
$\boxtimes E_{v}$ partout. 
On considère le   champ de Deligne-Mumford  
$$\mc Z^{(\{1\},\{2\}, I)}=\restr{\Cht_{N,I \cup \{1,2\},W \boxtimes V\boxtimes V^{*}}^{(\{1\},\{2\},I)}}
 {(X\sm (N\cup v))^{I}\times \Delta(v)}.$$
On va construire deux sous-champs fermés   $\mc Y_{1}$ et $\mc Y_{2}$ dans $\mc Z^{(\{1\},\{2\}, I)}$,    munis de   morphismes 
 $\alpha_{1}$ et $\alpha_{2}$ vers  $$ \mc Z^{(I)}=\restr{\Cht_{N,I ,W }^{(I)}}
 {(X\sm (N\cup v))^{I}}$$ de sorte que 
 \begin{itemize}
 \item {\bf A)} la  restriction à  $(X\sm (N\cup v))^{I}$ de la  composée 
\eqref{def-SVv-intro1}$\to$\eqref{def-SVv-intro2}$\to$\eqref{def-SVv-intro3}
du  morphisme de création
 et de l'action du morphisme de Frobenius partiel  est réalisée par une correspondance cohomologique   supportée par la correspondance $\mc Y_{2}$ de  $\mc Z^{(I)}$ vers    $\mc Z^{(\{1\},\{2\}, I)}$, et dont la restriction aux ouverts de lissité est déterminée par son support, avec un coefficient correctif de $q^{-d/2}$, 
  \item {\bf B)} la  restriction à $(X\sm (N\cup v))^{I}$ du morphisme d'annihilation    \eqref{def-SVv-intro3}$\to$\eqref{def-SVv-intro4}  est réalisée par une  correspondance cohomologique   supportée par    la   correspondance $\mc Y_{1}$ de  $\mc Z^{(\{1\},\{2\}, I)}$ vers $\mc Z^{(I)}$, et dont la restriction aux ouverts de lissité est déterminée par son support.   
 \end{itemize}
Donc $S_{V,v}$ sera réalisée par une correspondance cohomologique supportée par   le produit 
    $\mc Y_{1}\times_{\mc Z^{(\{1\},\{2\}, I)}}\mc Y_{2}$ de ces correspondances. On verra 
    \begin{itemize}
    \item que le  produit $\mc Y_{1}\times_{\mc Z^{(\{1\},\{2\}, I)}}\mc Y_{2}$ n'est autre que la   correspondance  de Hecke $\Gamma^{(I)}$ de  
$ \mc Z^{(I)}$ dans lui-même  (qui est une correspondance finie, et même étale)
\item que $S_{V,v}$, qui est donc une correspondance cohomologique supportée par $\Gamma^{(I)}$ est en fait égale à la correspondance cohomologique évidente  supportée par $\Gamma^{(I)}$ avec un coefficient correctif de $q^{-d/2}$  
 (qui réalise 
  $T(h_{V,v})$ puisque $V$ est minuscule). \end{itemize}
   Grâce à l'hypothèse que $V$ est minuscule il suffira de faire  le  calcul  sur les ouverts de lissité, et le calcul  sera évident car on verra que sur les ouverts de lissité l'intersection $\mc Y_{1}\times_{\mc Z^{(\{1\},\{2\}, I)}}\mc Y_{2}$ 
   est une intersection transverse de deux sous-champs lisses.

 On construit maintenant tous ces objets. 
 La correspondance de Hecke  $\Gamma^{(I)}$ est le  champ 
classifiant la donnée de   $(x_i)_{i\in I}$ et d'un diagramme 
\begin{gather} \label{intro-diag-Gamma}
 \xymatrix{
 (\mc G', \psi') \ar[r]^-{\phi'} & 
 (\ta{\mc G'}, \ta \psi')   \\
  (\mc G, \psi)  \ar[r]^-{\phi}  \ar[u]_-{\kappa}  &
 (\ta{\mc G}, \ta \psi) \ar[u]_-{\ta \kappa}
 } \end{gather}
 tel que   
 \begin{itemize}
 \item la ligne inférieure $\big( (x_i)_{i\in I}, (\mc G, \psi) \xrightarrow{\phi}   (\ta{\mc G}, \ta \psi)
\big)$
et la ligne supérieure  $\big( (x_i)_{i\in I}, (\mc G', \psi') \xrightarrow{\phi'}   (\ta{\mc G'}, \ta \psi')
\big)
$ appartiennent à  $\mc Z^{(I)}$, 
\item  $\kappa:\restr{\mc G}{(X\sm v)\times S}\isom \restr{\mc G'}{(X\sm v)\times S}$  est un  isomorphisme 
tel que la   position relative de $\mc G$ par rapport à  $\mc G'$ en   $v$ est {\it égale} au poids dominant de $V$ (on rappelle que  $V$ est minuscule),
\item la restriction de $\kappa$ à $N\times S$, qui est un isomorphisme,  entrelace les structures de niveau $\psi$ et $\psi'$.   
\end{itemize}
De plus les deux   projections 
 $\Gamma^{(I)}\to \mc Z^{(I)}$ sont les  morphismes qui conservent les lignes inférieures et supérieures de   \eqref{intro-diag-Gamma}. 
  
 Comme les  pattes indexées par   $I$ varient dans   $X\sm (N\cup v)$ et restent   disjointes  des  pattes $1$ et   $2$ fixées en  $v$, on 
 peut changer la partition $(\{1\},\{2\}, I)$ en $(\{1\},I,\{2\})$ et on a donc
    \begin{gather*}\mc Z^{(\{1\},\{2\}, I)}
=
 \restr{\Cht_{N,I \cup \{1,2\},W \boxtimes V\boxtimes V^{*}}^{(\{1\},I,\{2\})}}
 {(X\sm (N\cup v))^{I}\times \Delta(v)}.\end{gather*}
 Autrement dit le champ  $\mc Z^{(\{1\},\{2\}, I)}$ classifie la donnée de  $(x_i)_{i\in I}$ et d'un diagramme 
 \begin{gather} \label{intro-diag-W}
 \xymatrix{
 & (\mc G_{1}, \psi_{1}) \ar[r]^-{\phi'_{2}} \ar[d]^-{\phi_{2}}& 
 (\mc G'_{2}, \psi'_{2}) \ar[d]^-{\phi'_{3}} & 
 (\ta{\mc G_{1}}, \ta \psi_{1}) \\
 (\mc G_{0}, \psi_{0}) \ar[ru]^-{\phi_{1}} &
 (\mc G_{2}, \psi_{2}) \ar[r]^-{\phi_{3}}    &
 (\ta{\mc G_{0}}, \ta \psi_{0})\ar[ru]^-{\ta \phi_{1}} &
 } \end{gather}
 avec 
\begin{gather*}\big( (x_i)_{i\in I}, (\mc G_{0}, \psi_{0}) \xrightarrow{\phi_{1}}   (\mc G_{1}, \psi_{1}) \xrightarrow{\phi_{2}} (\mc G_{2}, \psi_{2}) \xrightarrow{\phi_{3}}    (\ta{\mc G_{0}}, \ta \psi_{0})
\big)\\
\in \restr{\Cht_{N,I \cup \{1,2\},W \boxtimes V\boxtimes V^{*}}^{(\{1\},\{2\},I)}}
 {(X\sm (N\cup v))^{I}\times \Delta(v)}\end{gather*} et  \begin{gather*}\big( (x_i)_{i\in I}, (\mc G_{0}, \psi_{0}) \xrightarrow{\phi_{1}}   (\mc G_{1}, \psi_{1}) \xrightarrow{\phi'_{2}} (\mc G'_{2}, \psi'_{2}) \xrightarrow{\phi'_{3}}    (\ta{\mc G_{0}}, \ta \psi_{0})
\big)\\
\in \restr{\Cht_{N,I \cup \{1,2\},W \boxtimes V\boxtimes V^{*}}^{(\{1\},I,\{2\})}}
 {(X\sm (N\cup v))^{I}\times \Delta(v)}. \end{gather*}  
 Les flèches  obliques, verticales et horizontales du  diagramme  \eqref{intro-diag-W} sont respectivement les   modifications associées à la  patte  $1$, à  la  patte  $2$ et aux pattes indexées par   $I$. La flèche $\ta \phi_{1}$ à droite du diagramme \eqref{intro-diag-W} est déterminée par $\phi_{1}$, mais on l'a dessinée car elle servira pour définir $\mc Y_{2}$ ci-dessous.

 On note  $\mc Y_{1} \overset{\iota_{1}}{\hookrightarrow} 
 \mc Z^{(\{1\},\{2\}, I)}$ le sous-champ fermé défini par la condition que dans le diagramme  \eqref{intro-diag-W},  $\phi_{2}\phi_{1}: 
 \restr{\mc G_{0}}{(X-v)\times S}\to \restr{\mc G_{2}}{(X-v)\times S}$ 
 s'étend en un  isomorphisme sur  $X\times S$. 
On a un   morphisme $$\alpha_{1}: \mc Y_{1}\to \mc Z^{(I)}
 $$
qui envoie    \begin{gather*}
 \xymatrix{
 & (\mc G_{1}, \psi_{1}) \ar[r]^-{\phi'_{2}} \ar[d]^-{\phi_{2}}& 
 (\mc G'_{2}, \psi'_{2}) \ar[d]^-{\phi'_{3}} & 
 (\ta{\mc G_{1}}, \ta \psi_{1}) \\
 (\mc G_{0}, \psi_{0}) \ar[ru]^-{\phi_{1}} \ar[r]^-{\sim}&
 (\mc G_{2}, \psi_{2}) \ar[r]^-{\phi_{3}}    &
 (\ta{\mc G_{0}}, \ta \psi_{0})\ar[ru]^-{\ta \phi_{1}} &
 } \end{gather*}
sur la ligne du bas, c'est-à-dire   
 \begin{gather}\label{intro-ligne-bas}\big( (x_i)_{i\in I}, (\mc G_{0}, \psi_{0}) \xrightarrow{\phi_{3}(\phi_{2}\phi_{1})}    (\ta{\mc G_{0}}, \ta \psi_{0})
\big). \end{gather}

L'assertion  B)  ci-dessus vient d'un énoncé similaire concernant les faisceaux de Mirkovic-Vilonen. En effet 
\begin{itemize}\item par le a) du \thmref{satake-geom-thm} l'image directe de 
 $\mc S_{\{1,2\},V\boxtimes V^{*} }^{(\{1\}, \{2\})}$ (qui est le faisceau constant $E$ décalé) par le morphisme d'oubli
 (de la modification intermédiaire) 
 $\mr{Gr}_{\{1,2\},V\boxtimes V^{*} }^{(\{1\}, \{2\})}\to 
 \mr{Gr}_{\{1,2\},V\boxtimes V^{*} }^{(\{1,2\})}$ est égale à 
 $\mc S_{\{1,2\},V\boxtimes V^{*} }^{(\{1,2\})}$, 
 \item  par le c) du \thmref{satake-geom-thm} la restriction de  $\mc S_{\{1,2\},V\boxtimes V^{*} }^{(\{1,2\})}$ au-dessus de la diagonale 
 (et donc en particulier au-dessus de $\Delta(v)$) est égale à 
  $\mc S_{\{0\},V\otimes V^{*} }^{(\{0\})}$ que l'on envoie dans 
  le faisceau gratte-ciel $\mc S_{\{0\},\mbf 1}^{(\{0\})}$ par $\on{ev}_{V}:V\otimes V^{*}\to \mbf 1$
\end{itemize}
et par le théorème de changement de base propre cela donne lieu à une correspondance cohomologique entre $\restr{\mr{Gr}_{\{1,2\},V\boxtimes V^{*} }^{(\{1\}, \{2\})}}{\Delta(v)}$ et le point, et on vérifie que celle-ci est la correspondance cohomologique  évidente supportée par le sous-schéma fermé lisse de $\restr{\mr{Gr}_{\{1,2\},V\boxtimes V^{*} }^{(\{1\}, \{2\})}}{\Delta(v)}$ formé des 
$ (\mc G_{0} \xrightarrow{\phi_{1}}  
\mc G_{1}\xrightarrow{\phi_{2}}
  \mc G_{2}  \isom G ) $ tels que $\phi_{2}\phi_{1}$ soit un isomorphisme. 
  
  On note  
$\mc Y_{2} \overset{\iota_{2}}{\hookrightarrow} \mc Z^{(\{1\},\{2\}, I)}$ 
le sous-champ fermé  défini par la condition que  dans  le diagramme  \eqref{intro-diag-W},   $\ta \phi_{1}\phi_{3}': 
 \restr{\mc G'_{2}}{(X-v)\times S}\to \restr{\ta{\mc G_{1}}}{(X-v)\times S}$ s'étend en  un  isomorphisme sur $X\times S$. 
 On a un   morphisme $$\alpha_{2}: \mc Y_{2}\to \mc Z^{(I)}
 $$ 
qui envoie  \begin{gather}\label{intro-diag-élément-Y2}
 \xymatrix{
 & (\mc G_{1}, \psi_{1}) \ar[r]^-{\phi'_{2}} \ar[d]^-{\phi_{2}}& 
 (\mc G'_{2}, \psi'_{2}) \ar[d]^-{\phi'_{3}} \ar[r]^-{\sim}& 
 (\ta{\mc G_{1}}, \ta \psi_{1}) \\
 (\mc G_{0}, \psi_{0}) \ar[ru]^-{\phi_{1}} &
 (\mc G_{2}, \psi_{2}) \ar[r]^-{\phi_{3}}    &
 (\ta{\mc G_{0}}, \ta \psi_{0})\ar[ru]^-{\ta \phi_{1}} &
 } \end{gather}
sur la ligne du haut, c'est-à-dire   
 \begin{gather}\label{intro-ligne-haut}\big( (x_i)_{i\in I}, (\mc G_{1}, \psi_{1}) \xrightarrow{(\ta \phi_{1}\phi'_{3})\phi'_{2}}    (\ta{\mc G_{1}}, \ta \psi_{1})
\big). \end{gather}
La justification de l'assertion A)  ci-dessus se fait
\begin{itemize}
\item par un argument similaire à celui donné pour justifier B)
mais concernant cette fois-ci  $\delta_{V}:\mbf 1\to V\otimes V^{*}$ et le champ
$\restr{\Cht_{N,I \cup \{1,2\},W \boxtimes V\boxtimes V^{*}}^{(I,\{2\},\{1\})}}
 {(X\sm (N\cup v))^{I}\times \Delta(v)}$
\item par le fait que la restriction à  $ (X\sm (N\cup v))^{I}\times \Delta(v)$ du morphisme de Frobenius partiel
$$\on {Fr}_{\{1\}} ^{(\{1\},I,\{2\})}:  \Cht_{N,I \cup \{1,2\},W \boxtimes V\boxtimes V^{*}}^{(\{1\},I,\{2\})} \to 
 \Cht_{N,I \cup \{1,2\},W \boxtimes V\boxtimes V^{*}}^{(I,\{2\},\{1\})}$$
envoie \eqref{intro-diag-W} sur 
\begin{gather}  \nonumber 
 \xymatrix{
   (\mc G_{1}, \psi_{1}) \ar[r]^-{\phi'_{2}}  & 
 (\mc G'_{2}, \psi'_{2}) \ar[d]^-{\phi'_{3}} & 
 (\ta{\mc G_{1}}, \ta \psi_{1}) \\
      &
 (\ta{\mc G_{0}}, \ta \psi_{0})\ar[ru]^-{\ta \phi_{1}} &
 } \end{gather}
\end{itemize}

Le fait qu'on n'ait pas besoin d'introduire de signe  
dans A) et B) ci-dessus est justifié dans la remarque 6.9 de \cite{coh}, dont l'idée est la suivante. Il suffit de trouver une situation où apparaissent un opérateur de création et un opérateur de création du même type (avec des pattes crées ou annihilées apparaissant dans le même ordre) et où on sait calculer leur composée directement. 
 Or la composée $V\xrightarrow{\Id_{V}\otimes \on{\delta}_{V}} 
    V\otimes V^{*} \otimes V\xrightarrow{ \on{ev}_{V} \otimes \Id_{V}} V$ est l'identité par le lemme de Zorro \eqref{zorro}. La composée des correspondances cohomologiques  entre champs de Hecke 
    données par les 
    opérateurs de création et d'annihilation associés à $\Id_{V}\otimes \on{\delta}_{V}$ et $ \on{ev}_{V} \otimes \Id_{V}$ 
    se calcule aisément et coïncide avec la correspondance identité avec les choix de signes comme dans A) et B) ci-dessus, ce qui prouve qu'ils étaient bons 
    (ou au moins que le produit des deux normalisations était correct, la normalisation de chacun étant plus subtile et inutile pour nous).

D'autre part on a un  isomorphisme canonique 
\begin{gather}\label{intro-Gamma-produit-fibre}
\Gamma^{(I)}\simeq \mc Y_{1} \times_{\mc Z^{(\{1\},\{2\}, I)}}\mc Y_{2}. 
\end{gather}
En effet un  point  de  $\mc Z^{(\{1\},\{2\}, I)}$ appartenant à   $\mc Y_{1} $ et  à $\mc Y_{2}$ est donné par un  diagramme 
  $$ \xymatrix{
 & (\mc G_{1}, \psi_{1}) \ar[r]^-{\phi'_{2}} \ar[d]^-{\phi_{2}}& 
 (\mc G'_{2}, \psi'_{2}) \ar[d]^-{\phi'_{3}} \ar[r]^-{\sim}& 
 (\ta{\mc G_{1}}, \ta \psi_{1}) \\
 (\mc G_{0}, \psi_{0}) \ar[ru]^-{\phi_{1}} \ar[r]^-{\sim}&
 (\mc G_{2}, \psi_{2}) \ar[r]^-{\phi_{3}}    &
 (\ta{\mc G_{0}}, \ta \psi_{0})\ar[ru]^-{\ta \phi_{1}} &
 } $$
Il équivaut donc  à la donnée d'un point de  $ \Gamma^{(I)}$, car  
 en contractant les deux   isomorphismes du   diagramme précédent on obtient  le diagramme 
  \begin{gather*}
 \begin{CD} 
  (\mc G_{1}, \psi_{1}) 
  @>(\ta \phi_{1}\phi'_{3}) \phi'_{2}>> 
  (\ta{\mc G_{1}}, \ta \psi_{1})  \\
 @AA\phi_{1}A 
 @AA\ta \phi_{1}A \\
(\mc G_{0}, \psi_{0})
 @>\phi_{3}(\phi_{2}\phi_{1})>>  
 (\ta{\mc G_{0}}, \ta \psi_{0}) 
 \end{CD}
 \end{gather*}
que l'on identifie au diagramme \eqref{intro-diag-Gamma}. 

On a des morphismes naturels des champs $\mc Z^{(\{1\},\{2\}, I)}, \mc Z^{(I)} , \mc Y_{1}, \mc Y_{2}$  et $\Gamma^{(I)}$ vers $ \mr{Gr}_{I,W}^{(I)}/
  G_{\sum n_{i} x_i}$. Comme $V$ est minuscule 
  ces morphismes sont lisses. Donc les 
   ouverts de lissité 
  ${}^{\circ}{\mc Z}^{(\{1\},\{2\}, I)}, {}^{\circ}{\mc Z}^{(I)},  {}^{\circ}{\mc Y}_{1}, {}^{\circ}{\mc Y}_{2}, {}^{\circ}{\Gamma}^{(I)}$
sont les images inverses de $ {}^{\circ}{\mr{Gr}}_{I,W}^{(I)}/
  G_{\sum n_{i} x_i}$ où $ {}^{\circ}{\mr{Gr}}_{I,W}^{(I)}$ désigne l'ouvert de lissité de 
$ \mr{Gr}_{I,W}^{(I)}$. 

Un calcul d'espaces tangents  montre que 
 ${}^{\circ}\mc Y_{1}$ et ${}^{\circ}\mc Y_{2}$ sont des sous-champs lisses dans le 
champ de Deligne-Mumford  lisse ${}^{\circ}\mc Z^{(\{1\},\{2\}, I)}$ et s'y s'intersectent transversalement, et de plus il résulte de    \eqref{intro-Gamma-produit-fibre} que leur  intersection est 
${}^{\circ}\Gamma^{(I)}$.   
On a donc l'égalité de  correspondances cohomologiques
entre $S_{V,v}$ et $T(h_{V,v})$ sur ${}^{\circ}\Gamma^{(I)}$ mais comme 
$\Gamma^{(I)}$ est une correspondance   étale  entre 
$\mc Z^{( I)}$ et lui-même l'égalité a lieu partout
(en effet un morphisme du faisceau pervers $\on{IC}_{\Gamma^{(I)}}$ dans lui-même est déterminé par sa restriction à ${}^{\circ}{\Gamma}^{(I)}$). 
\cqfd

    Une  conséquence de la   \propref{prop-coal-frob-cas-part-intro} est que l'on possède,  pour  tout    $f\in C_{c}(K_{N}\backslash G(\mb A)/K_{N},E)$ et  $\kappa$ assez grand,  
    une extension naturelle 
  du  morphisme $T(f)$ (introduit dans  \eqref{defi-Tf})  en un  morphisme  
 $T(f):\mc H_{N,I,W}^{\leq\mu}\to 
 \mc H_{N,I,W}^{\leq\mu+\kappa}$ de faisceaux constructibles sur  $  (X\sm N)^{I}$ tout entier, de fa\c con compatible avec la  composition des   opérateurs de Hecke. En effet, en notant $K_{N}=\prod K_{N,v}$,  il suffit de le montrer pour toute place $v$ et pour   $f\in C_{c}(K_{N,v}\backslash G(F_{v})/K_{N,v},E)$.
  Il n'y a rien à faire  si $v\in N$. Si $v\not\in N$ 
  il suffit de traiter le cas où $f=h_{V,v}$ et alors 
  l'extension est donnée par $S_{V,v}$ grâce à la  \propref{prop-coal-frob-cas-part-intro}. Pour plus de détails, on renvoie au corollaire 6.5 
  de \cite{coh}.

 Pour les variétés de Shimura sur les corps de nombres de tels prolongements ont été définis dans de nombreux cas, de fa\c con modulaire, par adhérence de  Zariski ou à l'aide de  cycles proches (voir \cite{deligne-bki-mod,faltings-chai, genestier-tilouine}). 

     Comme $S_{V,v}$ est le prolongement de 
   $T(h_{V,v})$, la proposition suivante exprime exactement la 
  relation d'Eichler-Shimura. On utilise encore $\{0\}$ pour noter un ensemble à un élement (indexant la patte à laquelle s'applique la relation  d'Eichler-Shimura). 
 
 \begin{prop}\label{Eichler-Shimura-intro}
 Soient $I,W$ comme ci-dessus et $V$  une   représentation irréductible de  $\wh G$. Alors  $$F_{\{0\}}^{\deg(v)}: \varinjlim _{\mu}\restr{\mc H _{N, I\cup\{0\}, W\boxtimes V}^{\leq\mu}}{(X\sm N)^{I}\times v}\to \varinjlim _{\mu}\restr{\mc H _{N, I\cup\{0\}, W\boxtimes V}^{\leq\mu}}{(X\sm N)^{I}\times v}$$ est annulé par un  polynôme de degré 
 $\dim(V)$ dont les  coefficients sont les 
  restrictions à $(X\sm N)^{I}\times v$ des morphismes $S_{\Lambda^{i}V,v}$. 
Plus précisément on a  
  \begin{gather*}  \sum_{i=0}^{\dim V} (-1)^{i} (F_{\{0\}}^{\deg(v)})^{i}\circ \restr{S_{\Lambda^{\dim V-i}V,v}}{(X\sm N)^{I}\times v}=0  . \end{gather*} 
    \end{prop}
    On rappelle que  $S_{\Lambda^{i}V,v}$ étend  l'opérateur de Hecke 
 $T(h_{\Lambda^{i}V,,v})$ $$\text{de \ \ }  (X\sm (N\cup v))^{I\cup \{0\}}\text{ \ \  à \ \ }(X\sm N)^{I\cup \{0\}}$$ et on remarque que cette  extension est absolument nécessaire pour prendre la restriction à  $(X\sm N)^{I}\times v$. Grâce à la définition des morphismes   $S_{\Lambda^{i}V,v}$ par  \eqref{def-SVv-intro1}-\eqref{def-SVv-intro4}, 
 la preuve de la \propref{Eichler-Shimura-intro} est un simple calcul d'algèbre tensorielle  (inspiré d'une démonstration  du théorème de  Hamilton-Cayley, et fondée uniquement sur le fait que $\Lambda^{\dim V+1}V=0$). 
 On renvoie au chapitre 7 de \cite{coh} pour cette preuve. 
 
 \begin{rem} Dans \cite{xiao-zhu}, Liang Xiao et  Xinwen Zhu ont  
défini (dans un cadre un peu différent) un anneau de correspondances cohomologiques entre $\restr{ \Cht_{ N,I\cup \{0\},W\boxtimes V}^{(I,\{0\})}}{(X\sm N)^{I}\times v}$ et lui-même, dans lequel la relation d'Eichler-Shimura résulte formellement de Hamilton-Cayley. 
\end{rem}

\section{ Sous-faisceaux constructibles stables sous l'action des  morphismes de  Frobenius partiels}
\label{section-sous-faisceaux-constr}
Le but de ce paragraphe est de montrer le \lemref{lem-Hf-union-stab} et   la \propref{cor-action-Hf-intro}, qui avaient été admis  dans le paragraphe \ref{section-esquisse-abc}. On renvoie au chapitre 8 de \cite{coh} pour plus de détails. 

  On rappelle que l'on a fixé un  point géométrique  $\ov\eta=\on{Spec}(\ov F)$  au-dessus du  point  générique  $\eta$ de $X$.   
   Soit $I$ un ensemble fini et  $W=\boxtimes_{i\in I}W_{i}$ une   représentation irréductible de $(\wh G)^{I}$. 
 On note  $\Delta:X\to X^{I}$ le  morphisme diagonal. 
 On rappelle qu'on a fixé un  point géométrique $\ov{\eta^{I}}$ au-dessus du point  générique $\eta^{I}$ de $X^{I}$ et une flèche de spécialisation   
 $\on{\mf{sp}}: \ov{\eta^{I}}\to \Delta(\ov \eta)$.

 \begin{lem}\label{lem-Hf-union-stab}
 L'espace   $\Big( \varinjlim _{\mu}\restr{\mc H _{N, I, W}^{\leq\mu}}{\ov{\eta^{I}}}\Big)^{\mr{Hf}}$  est la réunion  de sous-$\mc O_{E}$-modules $\mf M=\restr{\mf G}{\ov{\eta^{I}}}$ où 
   $\mf G$ est un sous-$\mc O_{E}$-faisceau constructible   de 
  $\varinjlim _{\mu}\restr{\mc H _{N, I, W}^{\leq\mu}}{\eta^{I}}$ 
     stable sous l'action des morphismes de  Frobenius partiels. 
 \end{lem}
  \noindent {\bf Démonstration.} On renvoie à la démonstration de   la proposition 8.27 de \cite{coh} pour plus de détails. 
Pour toute famille  $(v_{i})_{i\in I}$ de points fermés de  $X\sm N$, on note  ${\times_{i\in I} v_{i}}$ leur produit, qui est une réunion finie de  points fermés de  $(X\sm N)^{I}$. 
  Soit $\check{\mf M}$ un  sous-$\mc O_{E}$-module de type fini de 
  $\varinjlim _{\mu}\restr{\mc H _{N, I, W}^{\leq\mu}}{\ov{\eta^{I}}}$ stable par  
$\pi_{1}(\eta^{I},\ov{\eta^{I}})$ et   $C_{c}(K_{N}\backslash G(\mb A)/K_{N},\mc O_{E})$. 
On va construire $\mf M\supset \check{\mf M}$ vérifiant les propriétés de l'énoncé. 
Comme  $\check{\mf M}$ est de type fini, il existe $\mu_0$ tel que  $\check{\mf M}$ soit inclus dans  l'image de  $\restr{\mc H _{N, I, W}^{\leq\mu_0}}{\ov{\eta^{I}}}$ dans  $\varinjlim _{\mu}\restr{\mc H _{N, I, W}^{\leq\mu}}{\ov{\eta^{I}}}$. Quitte à augmenter  $\mu_0$, on peut supposer que   $\check{\mf M}$ est un sous-$\mc O_{E}$-module  de  $\restr{\mc H _{N, I, W}^{\leq\mu_0}}{\ov{\eta^{I}}}$. 
Soit $\Omega_0$ un ouvert dense de  $X^{I}$ sur lequel  
$ \mc H _{N, I, W}^{\leq\mu_0}$ est lisse. Il existe un  unique 
sous-$\mc O_{E}$-faisceau lisse $\check{\mf G}\subset  \restr{\mc H _{N, I, W}^{\leq\mu_0}}{\Omega_0}$ sur  $\Omega_0$ tel que  $\restr{\check{\mf G}}{\ov{\eta^{I}}}=\check{\mf M}$. On choisit  
 $(v_{i})_{i\in I}$ tel que  ${\times_{i\in I} v_{i}}$ soit inclus dans  $\Omega_0$. Pour tout  $i$, la  relation d'Eichler-Shimura  
 (\propref{Eichler-Shimura-intro}) 
 implique alors que   \begin{gather}\label{intro-ES-inclusion}(F_{\{i\}}^{\deg(v_{i})})^{\dim W_{i}}(\restr{\check{\mf G}}{\times_{i\in I} v_{i}})
   \subset 
   \sum_{\alpha=0}^{\dim W_{i}-1}  (F_{\{i\}}^{\deg(v_{i})})^{\alpha}(S_{\Lambda^{\dim W_{i}-\alpha}W_{i},v_{i}}(\restr{\check{\mf G}}{\times_{i\in I} v_{i}}))
  \end{gather} dans  $
  \varinjlim_{\mu} \restr{\mc H _{N, I, W}^{\leq\mu}}{\times_{i\in I} v_{i}}$. 
Grâce à la lissité de 
$(\Frob_{\{i\}}^{\deg(v_{i})\dim W_{i}})^{*}(\check{\mf G})$
en ${\times_{i\in I} v_{i}}$, l'inclusion   \eqref{intro-ES-inclusion} se propage  en  $\eta^{I}$, c'est-à-dire  que 
 \begin{gather*}
 F_{\{i\}}^{\deg(v_{i})\dim W_{i}}
 (\restr{(\Frob_{\{i\}}^{\deg(v_{i})\dim W_{i}})^{*}(\check{\mf G}}{\eta^{I}}) )
 \\   \subset 
   \sum_{\alpha=0}^{\dim W_{i}-1}  F_{\{i\}}^{\deg(v_{i})\alpha}(\Frob_{\{i\}}^{\deg(v_{i})\alpha})^{*}(\restr{S_{\Lambda^{\dim W_{i}-\alpha}W_{i},v_i}(\check{\mf G})}{\eta^{I}})
  \end{gather*} dans  $
  \varinjlim_{\mu} \restr{\mc H _{N, I, W}^{\leq\mu}}{\eta^{I}}$. 
 Or   $\restr{\check{\mf G}}{\eta^{I}}$ est stable par  $S_{\Lambda^{\dim W_{i}-\alpha}W_{i},v_i}=T(h_{\Lambda^{\dim W_{i}-\alpha}W_{i},v_i})$ 
 puisque  
 $$h_{\Lambda^{\dim W_{i}-\alpha}W_{i},v_i}\in C_{c}(G(\mc O_{v_{i}})\backslash G(F_{v_{i}})/G(\mc O_{v_{i}}),\mc O_{E})\subset C_{c}(K_{N}\backslash G(\mb A)/K_{N},\mc O_{E}). $$
   Par conséquent 
 \begin{gather*}
 F_{\{i\}}^{\deg(v_{i})\dim W_{i}}
 (\restr{(\Frob_{\{i\}}^{\deg(v_{i})\dim W_{i}})^{*}(\check{\mf G}}{\eta^{I}}) )
   \subset 
   \sum_{\alpha=0}^{\dim W_{i}-1}  F_{\{i\}}^{\deg(v_{i})\alpha}(\Frob_{\{i\}}^{\deg(v_{i})\alpha})^{*}(\restr{\check{\mf G}}{\eta^{I}})
  \end{gather*} dans  $
  \varinjlim_{\mu} \restr{\mc H _{N, I, W}^{\leq\mu}}{\eta^{I}}$. 
  On en déduit que   $$\mf G=\sum_{(n_{i})_{i\in I}\in \prod _{i\in I}\{0,...,\deg(v_{i})\dim(W_{i})-1\}}\restr{\prod _{i\in I}F_{\{i\}}^{n_{i}}\Big(\prod _{i\in I}\Frob_{\{i\}}^{n_{i}}\Big)^{*}(\check{\mf G})}{\eta^{I}}$$
 est  un  sous-$\mc O_{E}$-faisceau constructible  de 
  $\varinjlim _{\mu}\restr{\mc H _{N, I, W}^{\leq\mu}}{\eta^{I}}$ 
qui  est stable sous  l'action des morphismes de  Frobenius partiels. Comme 
$\Big( \varinjlim _{\mu}\restr{\mc H _{N, I, W}^{\leq\mu}}{\ov{\eta^{I}}}\Big)^{\mr{Hf}}$  est la réunion  de sous-$\mc O_{E}$-modules $\check{\mf M}$ comme au début de la démonstration et que $ \mf M =\restr{ \mf G }{\ov{\eta^{I}}}$ contient $\check{\mf M}$ on obtient l'énoncé du lemme. 
 \cqfd

\begin{prop}\label{cor-action-Hf-intro}  
L'espace  $\Big( \varinjlim _{\mu}\restr{\mc H _{N, I, W}^{\leq\mu}}{\ov{\eta^{I}}}\Big)^{\mr{Hf}}$ est muni d'une action naturelle  de  $\pi_{1}(\eta,\ov{\eta})^{I}$. Plus précisément c'est une réunion de sous-$E$-espaces vectoriels   munis d'une action continue de 
$\pi_{1}(\eta,\ov{\eta})^{I}$. 
 \end{prop}
\dem Pour tout 
   sous-$\mc O_{E}$-faisceau constructible $\mf G$   de 
  $\varinjlim _{\mu}\restr{\mc H _{N, I, W}^{\leq\mu}}{\eta^{I}}$ 
     stable sous l'action des morphismes de  Frobenius partiels, 
      le lemme de  Drinfeld \ref{lem-Dr-intro} 
 fournit (grâce à  $\mf{sp}$ et à la \remref{rem-lem-Drinfeld})
  une action continue de   $\pi_{1}(\eta,\ov{\eta})^{I}$ sur 
$\mf M=\restr{\mf G}{\ov{\eta^{I}}}$. D'après le \lemref{lem-Hf-union-stab}, 
  $\Big( \varinjlim _{\mu}\restr{\mc H _{N, I, W}^{\leq\mu}}{\ov{\eta^{I}}}\Big)^{\mr{Hf}}$ est la réunion de  tels  $\mf M$. \cqfd

\begin{rem}\label{rem-WeilF} L'action de $\pi_{1}(\eta,\ov{\eta})^{I}$
sur $\Big( \varinjlim _{\mu}\restr{\mc H _{N, I, W}^{\leq\mu}}{\ov{\eta^{I}}}\Big)^{\mr{Hf}}$  est    déterminée de manière unique par les actions de 
$\pi_{1}(\eta^{I},\ov{\eta^{I}})$ et   des morphismes de Frobenius partiels. Cela résulte du \lemref{lem-Dr-intro}
mais voici une autre fa\c con de le voir  (pour plus de détails on renvoie au chapitre 8 de \cite{coh}). 
 Suivant \cite{drinfeld78}, on va définir 
   un groupe    $\on{FWeil}(\eta^{I},\ov{\eta^{I}})$ 
       \begin{itemize}
     \item qui   est une  extension de $\Z^{I}$ par   $\on{Ker}(\pi_{1}(\eta^{I},\ov{\eta^{I}})\to \wh \Z)$, 
     \item et qui, lorsque  $I$ est un   singleton, 
 s'identifie au groupe de   Weil   usuel 
   $\on{Weil} (\eta,\ov{\eta})=\pi_{1}(\eta,\ov{\eta})\times_{\wh \Z}\Z$.      
       \end{itemize}
       
   On note  $F^{I}$ le corps des fonctions de $X^{I}$, $(F^{I})^{\mr{perf}}$ son 
  perfectisé et $\ov{F^{I}}$ la clôture algébrique de $F^{I}$ telle que 
  $\ov{\eta^{I}}=\on{Spec}(\ov{F^{I}})$. 
   On définit alors 
          \begin{gather*}\on{FWeil}(\eta^{I},\ov{\eta^{I}})=
       \big\{\varepsilon \in \on{Aut}_{\ov\Fq}((\ov{F^{I}})), \exists (n_{i})_{i\in I}\in \Z^{I}, \restr{\varepsilon}{(F^{I})^{\mr{perf}}}=\prod_{i\in I}(\Frob_{\{i\}})^{n_{i}}\big\} 
  .\end{gather*}
Le choix de 
  $\mf{sp}$ fournit 
  une inclusion   $\ov{F}\otimes_{\ov\Fq} \cdots 
 \otimes_{\ov\Fq}\ov{F}\subset \ov{F^{I}} 
  $.  Par restriction des automorphismes, on  en déduit   un 
   morphisme  surjectif 
    \begin{gather}\label{morph-Weil-I}\on{FWeil}(\eta^{I},\ov{\eta^{I}})\to  \on{Weil} (\eta,\ov{\eta}) ^{I}.\end{gather} L'énoncé de la \propref{cor-action-Hf-intro}  se reformule  alors en disant que l'action naturelle de $\on{FWeil}(\eta^{I},\ov{\eta^{I}})$ 
    sur $\Big( \varinjlim _{\mu}\restr{\mc H _{N, I, W}^{\leq\mu}}{\ov{\eta^{I}}}\Big)^{\mr{Hf}}$
    se factorise par le morphisme  \eqref{morph-Weil-I}, et même à travers 
    $\pi_{1}(\eta,\ov{\eta})^{I}$.  
      \end{rem}

  \section{Homorphismes de spécialisation et cohomologie Hecke-finie} 
  \label{para-homom-spe-Hf}
   Le but de ce paragraphe est le \corref{bijectivite-Hecke-fini}, qui avait été admis dans le paragraphe \ref{section-esquisse-abc}. Soit $W=\boxtimes_{i\in I} W_{i}$ une  représentation   $E$-linéaire irréductible de $(\wh G)^{I}$.

  \begin{prop} \label{surjectivite-Hecke-fini-intro}
   L'image de l'homomorphisme de spécialisation  
  \begin{gather}\label{intro-sp}\on{\mf{sp}}^{*}: \varinjlim _{\mu}\restr{\mc H _{N, I, W}^{\leq\mu}}{\Delta(\ov{\eta})}\to 
  \varinjlim _{\mu}\restr{\mc H _{N, I, W}^{\leq\mu}}{\ov{\eta^{I}}}\end{gather}
 contient 
  $\Big( \varinjlim _{\mu}\restr{\mc H _{N, I, W}^{\leq\mu}}{\ov{\eta^{I}}}\Big)^{\mr{Hf}}$. 
  \end{prop}
   \noindent {\bf Démonstration.} On renvoie à la preuve de la proposition 8.31   de \cite{coh} pour plus de détails. 
   D'après le \lemref{lem-Hf-union-stab}, 
    $\Big( \varinjlim _{\mu}\restr{\mc H _{N, I, W}^{\leq\mu}}{\ov{\eta^{I}}}\Big)^{\mr{Hf}}$  est la réunion  de sous-$\mc O_{E}$-modules $\mf M=\restr{\mf G}{\ov{\eta^{I}}}$ où 
   $\mf G$ est un sous-$\mc O_{E}$-faisceau constructible   de 
  $\varinjlim _{\mu}\restr{\mc H _{N, I, W}^{\leq\mu}}{\eta^{I}}$ 
     stable sous l'action des morphismes de  Frobenius partiels. Il suffit donc de montrer qu'un tel $\mf M$ est inclus dans l'image de  \eqref{intro-sp}.  
  Soit $\mu_0$ assez grand  pour que   $\mf G$  soit un  sous-$\mc O_{E}$-faisceau  de    
  $\restr{\mc H _{N, I, W}^{\leq\mu_0}}{\eta^{I}}$. Soit  $\Omega_0$ un ouvert dense  de  $X^{I}$ tel que   $\restr{\mc H _{N, I, W}^{\leq\mu_0}}{\Omega_0}$ soit lisse. Alors    $\mf G$ se prolonge  en un  sous-$\mc O_{E}$-faisceau lisse de 
  $\restr{\mc H _{N, I, W}^{\leq\mu_0}}{\Omega_0}$. 
    Par la preuve du lemma 9.2.1 de \cite{eike-lau}, 
     l'ensemble des   $\prod_{i\in I}\Frob_{\{i\}}^{n_{i}}(\Delta(\eta))$ pour $(n_{i})_{i\in I}\in \N^{I}$ est Zariski dense dans $X^{I}$. Donc  
 on peut trouver   $(n_{i})_{i\in I}\in \N^{I}$ tel que 
  $\prod_{i\in I}\Frob_{\{i\}}^{n_{i}}(\Delta(\eta))$ appartienne à  $\Omega_0$. 
 Alors  
  $\restr{\mf G}{\prod_{i\in I}\Frob_{\{i\}}^{n_{i}}(\ov{\eta^{I}})}$ est inclus dans l'image de    \begin{gather}\label{intro-sp-Frob}\on{\wt{\mf{sp}}}^{*}: \varinjlim _{\mu}\restr{\mc H _{N, I, W}^{\leq\mu}}{\prod_{i\in I}\Frob_{\{i\}}^{n_{i}}(\Delta(\ov{\eta}))}\to 
  \varinjlim _{\mu}\restr{\mc H _{N, I, W}^{\leq\mu}}{\prod_{i\in I}\Frob_{\{i\}}^{n_{i}}(\ov{\eta^{I}})}\end{gather}
 pour toute  flèche de spécialisation $\wt{\mf{sp}}: \prod_{i\in I}\Frob_{\{i\}}^{n_{i}}(\ov{\eta^{I}})\to \prod_{i\in I}\Frob_{\{i\}}^{n_{i}}(\Delta(\ov{\eta}))$, et en particulier pour l'image de  $\mf{sp}$ par $\prod_{i\in I}\Frob_{\{i\}}^{n_{i}}$. Comme   $\mf G$ est stable sous l'action  des  morphismes de  Frobenius partiels, on en déduit que  
  $\mf M=\restr{\mf G}{\ov{\eta^{I}}}$ est inclus dans l'image de  \eqref{intro-sp}. \cqfd

\begin{prop}\label{injectivite-sp}
L'homomorphisme de spécialisation 
 \begin{gather}\label{sp*-sans-Hf-intro2}\on{\mf{sp}}^{*}: 
 \varinjlim _{\mu}\restr{\mc H _{N, I, W}^{\leq\mu}}{\Delta(\ov{\eta})} \to
  \varinjlim _{\mu}\restr{\mc H _{N, I, W}^{\leq\mu}}{\ov{\eta^{I}}}\end{gather} est injectif.   \end{prop}
\noindent{\bf Démonstration. }
Soit $a$ dans le noyau de \eqref{sp*-sans-Hf-intro2}. 
On choisit $\mu_{0}$ et $\wt a\in \restr{\mc H _{N, I, W}^{\leq\mu_{0}}}{\Delta(\ov{\eta})}$ tels que 
 $a$ soit  l'image de $\wt a$ dans 
 $\varinjlim _{\mu}\restr{\mc H _{N, I, W}^{\leq\mu}}{\Delta(\ov{\eta})} $. 
Soit $\Omega_{0}$ un ouvert dense de $X\sm N$ sur lequel 
$\Delta^{*}\big( \mc H _{N, I, W}^{\leq\mu_{0}}\big)$ est lisse. Soit $v\in |\Omega_{0}|$. On pose $d=\deg(v)$ pour raccourcir les formules. Soit $\ov v$ un point géométrique au-dessus de $v$. 
Soit $\on{\mf{sp}}_{v}:\ov \eta\to \ov v$ une flèche de spécialisation. On note encore  $\on{\mf{sp}}_{v}:\Delta(\ov \eta)\to \Delta(\ov v)$ la flèche de spécialisation qui s'en déduit. 
Grâce à la lissité de $\Delta^{*}\big( \mc H _{N, I, W}^{\leq\mu_{0}}\big)$  sur 
 $\Omega_{0}$ on possède un unique élément 
 $\wt b\in \restr{\mc H _{N, I, W}^{\leq\mu_{0}}}{\Delta(\ov{v})}$ tel que $\wt a =\on{\mf{sp}}_{v}^{*}(\wt b)$. 
 On note $b$ l'image de  $\wt b$ dans $\varinjlim _{\mu}\restr{\mc H _{N, I, W}^{\leq\mu}}{\Delta(\ov{v})}$, de sorte que $a$ est l'image de $b$ par 
 $$\on{\mf{sp}}_{v}^{*}: \varinjlim _{\mu}\restr{\mc H _{N, I, W}^{\leq\mu}}{\Delta(\ov{v})}\to 
  \varinjlim _{\mu}\restr{\mc H _{N, I, W}^{\leq\mu}}{\Delta(\ov{\eta})} .$$
  
L'action des morphismes de Frobenius partiels fournit  pour tout $\mu$  et pour tout 
$(n_{i})_{i\in I}\in \N^{I}$  
un morphisme de faisceaux sur $(X\sm N)^{I}$ 
\begin{gather}\label{action-Frob-partiels-preuve}\prod_{i\in I}F_{\{i\}}^{dn_{i}}:(\prod_{i\in I}\Frob_{\{i\}}^{dn_{i}})^{*} ( \mc H _{N, I, W}^{\leq\mu})\to 
 \mc H _{N, I, W}^{\leq\mu+\kappa(\sum n_{i})}\end{gather} 
avec  $\kappa $ assez grand en fonction de $W$ et de $d$.   
Comme $ \prod_{i\in I}\Frob_{\{i\}}^{dn_{i}} $ agit trivialement sur $\Delta(v)$, 
 $\prod_{i\in I}F_{\{i\}}^{dn_{i}}$ agit sur 
$\varinjlim _{\mu}\restr{\mc H _{N, I, W}^{\leq\mu}}{\Delta(\ov{v})}$. 

Pour tout $(n_{i})_{i\in I}\in \N^{I}$ on note 
$$b_{(n_{i})_{i\in I}} =\prod_{i\in I}F_{\{i\}}^{dn_{i}}(b)\in 
\varinjlim _{\mu}\restr{\mc H _{N, I, W}^{\leq\mu}}{\Delta(\ov{v})}. $$
En particulier 
$b_{(0)_{i\in I}}=b$. 

On pose 
\begin{gather}\label{def-a-n-i}a_{(n_{i})_{i\in I}} =\on{\mf{sp}}_{v}^{*}(b_{(n_{i})_{i\in I}} )\in 
\varinjlim _{\mu}\restr{\mc H _{N, I, W}^{\leq\mu}}{\Delta(\ov{\eta})}, \end{gather}
de sorte que $a_{(0)_{i\in I}}=a$. 

La suite $a_{(n_{i})_{i\in I}}$ vérifie  les deux propriétés énoncées dans le lemme suivant. La première affirme que cette suite est ``multirécurrente'', c'est-à-dire récurrente en chaque variable $n_{i}$, et la seconde implique qu'elle est ``presque partout'' nulle. On déduira aisément de la conjonction des deux propriétés que cette suite est partout nulle, et donc en particulier que 
$a=a_{(0)_{i\in I}}$ est nul. 

Pour énoncer la seconde propriété on remarque que   $\on{\mf{sp}}^{*}(\on{\mf{sp}}_{v}^{*}(b))=\on{\mf{sp}}^{*}(a)=0$  
 dans $\varinjlim _{\mu}\restr{\mc H _{N, I, W}^{\leq\mu}}{\ov{\eta^{I}}}
$. 
Donc il existe $\mu_{1}\geq \mu_{0}$ tel que 
$\on{\mf{sp}}^{*}(\on{\mf{sp}}_{v}^{*}(\wt b))\in \restr{\mc H _{N, I, W}^{\leq\mu_{0}}}{\ov{\eta^{I}}}$
ait une image   nulle dans $\restr{\mc H _{N, I, W}^{\leq\mu_{1}}}{\ov{\eta^{I}}}
$. Autrement dit en notant 
$\wh b$ l'image de $\wt b$ dans 
$\restr{\mc H _{N, I, W}^{\leq\mu_{1}}}{\Delta(\ov{v})}$, 
on a 
$\on{\mf{sp}}^{*}(\on{\mf{sp}}_{v}^{*}(\wh b))=0$ dans 
$\restr{\mc H _{N, I, W}^{\leq\mu_{1}}}{\ov{\eta^{I}}}$. 
Soit $\Omega_{1}\subset (X\sm N)^{I}$ un ouvert dense sur lequel 
$\mc H _{N, I, W}^{\leq\mu_{1}}$ est lisse. 

\begin{lem}\label{lem-a-b}
a) Pour tout $j\in I$ et pour tout $(n_{i})_{i\in I}\in \N^{I}$, 
\begin{gather} \label{relation-ani}
\sum_{\alpha=0}^{\dim W_{j}}(-1)^{\alpha} S_{\Lambda^{\dim W_{j}-\alpha}W_{j},v} (a_{(n_{i}+\alpha \delta_{i,j})_{i\in I}})=0
\end{gather}
dans $\varinjlim _{\mu}\restr{\mc H _{N, I, W}^{\leq\mu}}{\Delta(\ov{\eta})}$. 

b)  
Pour tout $(n_{i})_{i\in I}\in \N^{I}$ tel que $\prod_{i\in I}\Frob_{\{i\}}^{dn_{i}}
(\Delta(\ov\eta))\in \Omega_{1}$, on a $a_{(n_{i})_{i\in I}}=0$ dans 
$\varinjlim _{\mu}\restr{\mc H _{N, I, W}^{\leq\mu}}{\Delta(\ov{\eta})} $. 
\end{lem}

\noindent{\bf Démonstration de a). }
Les $b_{(n_{i})_{i\in I}} $ satisfont une  relation identique à 
\eqref{relation-ani} 
(dans $\varinjlim _{\mu}\restr{\mc H _{N, I, W}^{\leq\mu}}{\Delta(\ov{v})}$), à savoir la  relation d'Eichler-Shimura en la patte $j$ (\propref{Eichler-Shimura-intro}). 
Alors \eqref{relation-ani} s'obtient en appliquant 
$\on{\mf{sp}}_{v}^{*}$ à 
 cette relation (ce qui est légitime puisque les $S_{\Lambda^{\dim W_{j}-\alpha}W_{j},v}$ sont des morphismes de faisceaux).  \cqfd

\noindent{\bf Démonstration de b). } Soit $(n_{i})_{i\in I}$ satisfaisant   l'hypothèse  de b). 
Comme 
\eqref{action-Frob-partiels-preuve} 
est un morphisme de faisceaux sur $(X\sm N)^{I}$, 
on peut intervertir les homomorphismes de spécialisation et les morphismes de Frobenius partiels. 
Autrement dit 
on a un diagramme commutatif 
  \begin{gather*} \xymatrixcolsep{5pc} \xymatrix{
     \restr{ \mc H _{N, I, W}^{\leq\mu_{1}}}{\Delta(\ov v)}
=
\restr{(\prod_{i\in I}\Frob_{\{i\}}^{dn_{i}})^{*} ( \mc H _{N, I, W}^{\leq\mu_{1}})}{\Delta(\ov v)}\ar[d]^-{ \on{\mf{sp}}_{v,(n_{i})_{i\in I}}^{*} } \ar[r]^-{\prod_{i\in I}F_{\{i\}}^{dn_{i}}}
 & \varinjlim _{\mu}\restr{\mc H _{N, I, W}^{\leq\mu}}{\Delta(\ov{v})} \ar[d]^-{\on{\mf{sp}}_{v}^{*}} 
 \\
\restr{(\prod_{i\in I}\Frob_{\{i\}}^{dn_{i}})^{*} ( \mc H _{N, I, W}^{\leq\mu_{1}})}{\Delta(\ov \eta)}
 \ar[r]^-{\prod_{i\in I}F_{\{i\}}^{dn_{i}}}
& \varinjlim _{\mu}\restr{\mc H _{N, I, W}^{\leq\mu}}{\Delta(\ov{\eta})}  }\end{gather*}
 où la notation $ \on{\mf{sp}}_{v,(n_{i})_{i\in I}}^{*}$ indique que l'homomorphisme de spécialisation associé à la flèche $\on{\mf{sp}}_{v}:\Delta(\ov \eta)\to \Delta(\ov v)$ {\it est appliqué au  faisceau} $(\prod_{i\in I}\Frob_{\{i\}}^{dn_{i}})^{*} ( \mc H _{N, I, W}^{\leq\mu_{1}})$ (et non pas à 
 $\mc H _{N, I, W}^{\leq\mu_{1}}$). 
 Le diagramme précédent donne lieu à 
   \begin{gather*} \xymatrixcolsep{5pc} \xymatrix{
\wh b \ar@{|->}[d]^-{\on{\mf{sp}}_{v,(n_{i})_{i\in I}}^{*}} \ar@{|->}[r]^-{\prod_{i\in I}F_{\{i\}}^{dn_{i}}}
 &b_{(n_{i})_{i\in I}} \ar@{|->}[d]^-{\on{\mf{sp}}_{v}^{*}} 
 \\
\on{\mf{sp}}_{v,(n_{i})_{i\in I}}^{*}(\wh b) \ar@{|->}[r]^-{\prod_{i\in I}F_{\{i\}}^{dn_{i}}}
&a_{(n_{i})_{i\in I}}  }\end{gather*}

Donc  
pour montrer $a_{(n_{i})_{i\in I}} =0$ (et terminer ainsi la preuve de b)) il suffit de montrer que     
\begin{gather}\label{fibre-Frob-Delta}
  \on{\mf{sp}}_{v,(n_{i})_{i\in I}}^{*}(\wh b) \in 
\restr{ \mc H _{N, I, W}^{\leq\mu_{1}}}{(\prod_{i\in I}\Frob_{\{i\}}^{dn_{i}})(\Delta(\ov \eta))}=
\restr{(\prod_{i\in I}\Frob_{\{i\}}^{dn_{i}})^{*} ( \mc H _{N, I, W}^{\leq\mu_{1}})}{\Delta(\ov \eta)}\end{gather}
est nul. 
 Or \eqref{fibre-Frob-Delta} peut aussi être considéré comme l'image de $\wh b$ par un   homomorphisme de spécialisation pour le faisceau
 $\mc H _{N, I, W}^{\leq\mu_{1 }}$ mais associé  à une flèche de spécialisation $(\prod_{i\in I}\Frob_{\{i\}}^{dn_{i}})(\Delta(\ov \eta))\to \Delta(\ov v)$. On en déduit que  \eqref{fibre-Frob-Delta} est nul car 
 \begin{itemize}
 \item 
  $\prod_{i\in I}\Frob_{\{i\}}^{dn_{i}}
(\Delta(\ov\eta))$ appartient à $ \Omega_{1}$ par hypothèse 
\item pour tout point géométrique $\ov x$ de 
$\Omega_{1}$
et toute flèche de spécialisation $\on{\mf{sp}}_{\ov x}: \ov x\to \Delta(\ov v)$,     $\on{\mf{sp}}_{\ov x}^{*}(\wh  b )$ s'annule dans 
$\restr{\mc H _{N, I, W}^{\leq\mu_{1 }}}{\ov x}$. 
\end{itemize}
 Cette dernière assertion résulte du fait que   $\mc H _{N, I, W}^{\leq\mu_{1}}$ est lisse sur $\Omega _{1}$ 
et que  l'image de $\wh  b$ par tout homomorphisme de spécialisation vers 
$\restr{\mc H _{N, I, W}^{\leq\mu_{1}}}{\ov{\eta^{I}}}
$ est nulle (puisque c'est le cas de $\on{\mf{sp}}^{*}(\on{\mf{sp}}_{v}^{*}(\wh  b))$ et que 
$\pi_{1}(\eta^{I},\ov{\eta^{I}})$ agit transitivement sur les flèches de spécialisation de $\ov{\eta^{I}}$ vers $\Delta(\ov v)$).   \cqfd

\noindent{\bf Fin de la démonstration de la \propref{injectivite-sp}. }
 Comme $\prod _{i\in I}\Frob_{\{i\}}$ est le Frobenius total, 
 $\prod_{i\in I}F_{\{i\}}^{dn }$ agit de fa\c con bijective sur 
 $\varinjlim _{\mu}\restr{\mc H _{N, I, W}^{\leq\mu}}{\Delta(\ov{\eta})} $ et envoie $a_{(n_{i} )_{i\in I}}$ sur $a_{(n_{i}+n)_{i\in I}}$. De ceci et du a) du \lemref{lem-a-b} on  
   déduit facilement que pour montrer que $a=a_{(0)_{i\in I}}$ est nul 
   (et même que toute la suite $a_{(n_{i})_{i\in I}}$ est nulle) il suffit de  trouver 
 $(n_{i})_{i\in I}\in \N^{I}$ tel que 
 \begin{gather*} 
 a_{(n_{i}+\alpha_{i})_{i\in I}}=0 \text{ pour tout } (\alpha_{i})_{i\in I}\in \prod _{i\in I} \{0,..., \dim W_{i}-1\}.
 \end{gather*}
Or cela est possible d'après le b) du \lemref{lem-a-b}, car la densité de l'ouvert $\Omega_{1}$ implique  que l'on peut trouver  $(n_{i})_{i\in I}\in \N^{I}$ tel que 
\begin{gather*} 
 \prod_{i\in I}\Frob_{\{i\}}^{d(n_{i}+\alpha_{i})}
(\Delta(\ov\eta))\in \Omega_{1} \text{ pour tout } (\alpha_{i})_{i\in I}\in \prod _{i\in I} \{0,..., \dim W_{i}-1\}.
 \end{gather*} 
 Ceci termine la preuve de la \propref{injectivite-sp}. \cqfd
 
 Les propositions  \ref{surjectivite-Hecke-fini-intro} et \ref{injectivite-sp}   entraînent le corollaire suivant. 
 
 \begin{cor}\label{bijectivite-Hecke-fini}
L'homomorphisme de spécialisation 
 \begin{gather}\label{sp*-sans-Hf2}\on{\mf{sp}}^{*}: \Big( \varinjlim _{\mu}\restr{\mc H _{N, I, W}^{\leq\mu}}{\Delta(\ov{\eta})} \Big)^{\mr{Hf}}\to 
 \Big( \varinjlim _{\mu}\restr{\mc H _{N, I, W}^{\leq\mu}}{\ov{\eta^{I}}}\Big)^{\mr{Hf}}\end{gather}
 est une bijection.     \end{cor}
\dem
L'injectivité résulte de la \propref{injectivite-sp}. Voici la preuve de la surjectivité. Soit $c\in  \Big( \varinjlim _{\mu}\restr{\mc H _{N, I, W}^{\leq\mu}}{\ov{\eta^{I}}}\Big)^{\mr{Hf}}$. D'après la proposition \ref{surjectivite-Hecke-fini-intro} il existe $a\in  \varinjlim _{\mu}\restr{\mc H _{N, I, W}^{\leq\mu}}{\Delta(\ov{\eta})}$ tel que $\on{\mf{sp}}^{*}(a)=c$. L'injectivité de $\on{\mf{sp}}^{*}$ montrée dans 
la \propref{injectivite-sp} implique que $a$ est Hecke-fini. \cqfd

    \section{Compatibilité avec l'isomorphisme de Satake aux places non ramifiées} 
    \label{subsection-intro-decomp}
    Le but de ce paragraphe est de montrer le   lemme suivant, ainsi que   la \propref{S-non-ram-concl-intro} qui avait été admise à la fin du paragraphe \ref{intro-idee-heurist}.

  Le lemme suivant  montre le (v)  de la \propref{prop-SIf-i-ii-iii}, 
  à savoir que les 
    opérateurs de Hecke en les places non ramifiées sont des cas particuliers d'opérateurs d'excursion. 
    
    Soit 
      $v$ une place dans  $X\sm N$.   On fixe un plongement 
       $\ov F\subset \ov F_{v}$. 
    Comme précédemment  $\mbf 1\xrightarrow{\delta_{V}} V\otimes V^{*}$  et   $ V\otimes V^{*}\xrightarrow{\on{ev}_{V}} \mbf 1$ sont les morphismes naturels.

     \begin{lem} \label{S-non-ram-intro}     Pour tout $d\in \N$ et tout  $\gamma\in \on{Gal}(\ov {F_{v}}/F_{v})\subset  \on{Gal}(\ov F/F)$ tel que $\deg(\gamma)=d$, $
   S_{\{1,2\},V \boxtimes V^{*},\delta_{V},\on{ev}_{V},(\gamma,1)}$  dépend seulement de   $d$, 
  et si    $d=1$ il est égal à  $T(h_{V,v})$. 
     \end{lem}
       \noindent{\bf Démonstration.}
 On fixe un  point géométrique  $\ov v$ au-dessus de  $v$ et une  flèche de spécialisation 
    $\on{\mf{sp}}_{v}:\ov \eta\to \ov v$, associés au plongement $\ov F\subset \ov {F_{v}}$ choisi ci-dessus. On note encore   $\on{\mf{sp}}_{v}$ la  flèche de spécialisation    
    $\Delta(\ov \eta)\to \Delta(\ov v)$ égale à son image   par $\Delta$. Pour que le diagramme suivant tienne dans la page on pose   $I=\{1,2\}$ et $W=V\boxtimes V^{*}$. 
    Le diagramme  
       $$  \xymatrix{  C_{c}^{\mr{cusp}}(G(F)\backslash G(\mb A)/K_N \Xi,E)\ar[d]_-{\restr{\mc C_{\delta_{V}}^{\sharp}}{\ov v}} \ar[dr]^{ \mc C_{\delta_{V}}^{\sharp}} & &
       \\ \Big( \varinjlim _{\mu} \restr{\mc H _{N, I, W}^{\leq\mu}}{\Delta(\ov v)}\Big)^{\mr{Hf}} \ar[r]^-{\mf{sp}_{v}^{*}} \ar[d]^{F_{\{1\}}^{\deg(v)d}} 
       & 
       \Big( \varinjlim _{\mu} \restr{\mc H _{N, I, W}^{\leq\mu}}{\Delta(\ov\eta)}\Big)^{\mr{Hf}} \ar@{=}[r]  
       & H_{I,W} \ar[d]^{(\gamma,1)}
       \\  \Big( \varinjlim _{\mu} \restr{\mc H _{N, I, W}^{\leq\mu}}{\Delta(\ov v)}\Big)^{\mr{Hf}} \ar[r]^-{\mf{sp}_{v}^{*}}\ar[d]_-{\restr{\mc C_{\on{ev}_{V}}^{\flat}}{\ov v}}   & 
        \Big( \varinjlim _{\mu} \restr{\mc H _{N, I, W}^{\leq\mu}}{\Delta(\ov\eta)}\Big)^{\mr{Hf}} \ar@{=}[r]  \ar[dl]^{ \mc C_{\on{ev}_{V}}^{\flat}} & 
      H_{I,W}
       \\ C_{c}^{\mr{cusp}}(G(F)\backslash G(\mb A)/K_N \Xi,E)  & &
       } $$
   est commutatif (la commutativité du grand rectangle est montrée dans le lemme 10.4 de \cite{coh}). 
   Or    $S_{\{1,2\},V \boxtimes V^{*},\delta_{V},\on{ev}_{V},(\gamma,1)}$  est égal par définition à la composée par le chemin le plus à droite.    Donc      il  est égal à la  composée donnée par la colonne de gauche. Par conséquent  il dépend seulement de   $d$. 
   Lorsque $d=1$ la composée donnée par la colonne de gauche est égale par définition à $S_{V,v}$, et donc à $T(h_{V,v})$  par  la \propref{prop-coal-frob-cas-part-intro}. \cqfd

       \begin{rem} On n'a calculé la composée de la colonne de gauche que pour $d=1$ mais pour d'autres valeurs de $d$ elle n'apporte rien de nouveau 
       car on pourrait montrer qu'elle est égale à une combinaison de $S_{W,v}$ avec $W$ représentation irréductible de $\wh G$. 
       \end{rem}
       
   La proposition suivante affirme  la compatibilité de la décomposition \eqref{intro1-dec-canonique} avec 
l'isomorphisme de Satake  en les places de $X\sm N$.           

\begin{prop}\label{S-non-ram-concl-intro} 
Soit  $\sigma$ apparaissant 
 dans \eqref{intro1-dec-canonique}  et  
$v$ une place    de $X\sm N$. Alors   
 pour toute représentation irréductible $V$ de $\wh G$, 
       $T(h_{V,v})$ agit sur   $\mf H_{\sigma}$
       par multiplication par le scalaire  $\chi_{V}(\sigma(\Frob_{v}))$, où $\chi_{V}$ est le caractère de $V$ et $\Frob_{v}\in \pi_{1}(X\sm N, \ov\eta)$ est un   élément de  Frobenius en $v$. 
 \end{prop}
 \dem 
 On reprend les notations du \lemref{S-non-ram-intro}. 
 Comme 
    $  \s{\on{ev}_{V}, (\sigma(\gamma), 1) . \delta_{V} }=\chi_{V}(\sigma(\gamma))$ ce lemme implique que  
        pour  toute   représentation  irréductible $V$ de $\wh G$, 
   $\mf H_{\sigma}$ est inclus dans l'espace propre généralisé (ou espace caractéristique)   de $T(h_{V,v})$ pour la valeur propre 
             $\chi_{V}(\sigma(\Frob_{v}))$. Or on sait  que les opérateurs de Hecke aux places non ramifiées sont diagonalisables
      (car ce sont des opérateurs normaux sur l'espace hermitien des formes automorphes cuspidales à coefficients dans $\C$). Donc $T(h_{V,v})$ agit sur $\mf H_{\sigma}$  par homothétie de rapport 
      $\chi_{V}(\sigma(\Frob_{v}))$. 
\cqfd

 Cela termine  la preuve  du \thmref{intro-thm-ppal}.

      \section{Lien avec le programme de Langlands géométrique}
     \label{subsection-link-langl-geom}
  Il est évident que la coalescence   et la  permutation des pattes 
sont reliées aux  structures de  factorisation introduites par  
  Beilinson et Drinfeld
 \cite{chiral} et en effet notre article utilise de fa\c con essentielle   le produit de fusion  sur la grassmannienne affine de Beilinson-Drinfeld dans l'équivalence de  Satake  géométrique     \cite{mv,ga-de-jong}.  
 Par ailleurs l'idée de décomposition spectrale est familière dans le programme de Langlands géométrique, cf \cite{beilinson-heisenberg} et surtout 
 le corollaire  4.5.5 de 
 \cite{dennis-laumon} qui affirme (dans le cadre du programme de Langlands géométrique pour les $D$-modules où la courbe $X$ est définie sur un corps algébriquement clos de caractéristique $0$) que la DG-catégorie des $D$-modules sur 
 $\Bun_{G}$ est ``au-dessus'' du champ des systèmes locaux pour 
 $\wh G$. On notera curieusement que l'on ne sait pas formuler d'énoncé analogue avec les faisceaux $\ell$-adiques lorsque $X$ est sur $\Fq$
 (la conjecture d'annulation montrée par Gaitsgory \cite{ga-vanishing} apparaissant comme une partie émergée de l'iceberg), et cependant notre article peut être considéré comme une version ``classique'' ou plutôt ``arithmétique'' d'un tel  énoncé.

 En fait le lien est bien plus direct qu'une simple analogie: 
 nous allons voir   que les conjectures du programme de Langlands géométrique $\ell$-adique  permettent  de comprendre les opérateurs d'excursion et fournissent une explication très éclairante de notre approche grâce à une construction de Braverman et Varshavsky \cite{brav-var} qui généralise le fait qu'un faisceau sur $\Bun_{G}$ donne par les  traces de Frobenius une fonction sur $\Bun_{G}(\Fq)$. 
On prend ici  $N$ vide, c'est-à-dire  $K_{N}=G(\mathbb O)$ mais les considérations qui suivent resteraient valables pour tout niveau $N$ (et même pour les groupes réductifs non déployés).

   Les conjectures du programme de Langlands   géométrique sont formulées à l'aide des foncteurs de  Hecke:  pour toute représentation $W$ de $(\wh G)^{I}$  le foncteur de  Hecke
     $$\phi_{I,W}:D^{b}_{c}(\Bun_{G},\Qlbar)\to  D^{b}_{c}(\Bun_{G}\times X^{I}, \Qlbar)$$      
    est donné  par 
     $$\phi_{I,W}(\mc F)=q_{1,!}\big(q_{0}^{*}(\mc F)\otimes \mc F_{I,W}\big)$$
 où $\Bun_{G}\xleftarrow{q_{0}}\Hecke_{I,W}^{(I)}\xrightarrow{q_{1}}\Bun_{G}\times X^{I}$ est la   correspondance de Hecke    et  
 \begin{itemize}
 \item quand  $W$ est irréductible, $\mc F_{I,W}$ est égal,  à un décalage près,  au faisceau d'intersection  de $\Hecke_{I,W}^{(I)}$   
 \item en général il est défini, de fa\c con fonctorielle en $W$,  comme  l'image inverse de 
 $\mc S_{I,W}^{(I)}$ par le morphisme lisse naturel 
 $\Hecke_{I,W}^{(I)}\to \mr{Gr}_{I,W }^{(I)}/G_{\sum n_{i}x_{i}}$ (où les $n_{i}$ sont assez grands).  
 \end{itemize}
Soit $\mc E$  un $\wh G$-système local sur $X$. Alors  $\mc F\in D^{b}_{c}(\Bun_{G},\Qlbar)$ est dit propre pour  $\mc E$  si l'on possède,  pour tout ensemble fini $I$ et toute représentation $W$ de 
  $(\wh G)^{I}$,  un isomorphisme 
  $\phi_{I,W}(\mc F)\isom \mc F\boxtimes W_{\mc E}$,   fonctoriel en  $W$, et compatible aux produits extérieurs et à la  fusion (c'est-à-dire  à l'image inverse par le   morphisme diagonal $X^{J}\to X^{I}$ associé à n'importe quelle application $I\to J$). Les  conjectures du programme de Langlands   géométrique impliquent l'existence d'un objet $\mc F$ propre pour  $\mc E$   (et vérifiant une condition de normalisation de Whittaker qui l'empêche en particulier d'être nul). Dans le programme de Langlands   géométrique  $X$ et $\Bun_{G}$ sont habituellement définis sur un corps algébriquement clos, mais ici nous les prenons définis sur $\Fq$. 
 
Soit $\mc F$ propre pour  $\mc E$. On note  $f\in C(\Bun_{G}(\Fq),\Qlbar)$ la fonction associée à $\mc F$ par le dictionnaire  faisceaux-fonctions, c'est-à-dire que  pour $x\in \Bun_{G}(\Fq)$, $f(x)=\on{Tr}(\Frob_{x}, \restr{\mc F}{x})$. 
Soit $\Xi\subset Z(F)\backslash Z(\mathbb A)$ un réseau. On suppose que  $\mc F$ est $\Xi$-équivariant, si bien que  $f\in  C(\Bun_{G}(\Fq)/\Xi,\Qlbar)$ (quitte à diminuer $\Xi$ cela est impliqué par une condition   sur $\mc E$, en fait sur son image par le morphisme de $\wh G$ vers son abélianisé). 
Il est bien connu que  $f$ est un vecteur propre pour tous les  opérateurs de Hecke: pour toute place $v$ et toute représentation  irréductible $V$ de $\wh G$, 
$T(h_{V,v})(f)=\on{Tr}(\Frob_{v},\restr{V_{\mc E}}{v})f$, où 
$\Frob_{v}$ est un élément de  Frobenius en $v$. 

La proposition suivante (reposant sur un résultat non encore rédigé) exprime la compatibilité entre  le programme de Langlands géométrique et la décomposition \eqref{intro1-dec-canonique}. 

\begin{prop}\label{prop-compat-langl-geom}
Etant donné $\mc F$ propre pour $\mc E$ et $\Xi$-équivariant tel que   la fonction $f$ associée à $\mc F$ soit cuspidale, alors $f$ appartient à $\mc H_{\sigma}$  où 
 $\sigma:\pi_{1}(X,\ov\eta)\to 
\wh G(\Qlbar)$ est la   représentation galoisienne correspondant au système local $\mc E$. 
\end{prop}
\dem 
Dans  \cite{brav-var}, Braverman et Varshavsky utilisent un morphisme de trace très général, et le fait que  
$\Cht_{I,W}^{(I)}$ est l'intersection de la  correspondance de Hecke avec le  graphe de l'endomorphisme de  Frobenius  de $\Bun_{G}$, pour construire  un  morphisme de faisceaux sur $X^{I}$
\begin{gather}\label{trace-brav-var-piIW}\pi^{\mc F,\mc E}_{I,W}: \varinjlim_{\mu}\mc H_{N,I,W}^{\leq\mu}\to W_{\mc E}.\end{gather}
Ces morphismes sont fonctoriels et $\Qlbar$-linéaires en  $W$, et compatibles avec la  coalescence des pattes et avec l'action des morphismes de  Frobenius partiels (ce dernier point n'a pas encore été rédigé). 
De plus   $\pi^{\mc F,\mc E}_{\emptyset,\mbf 1}:C_{c}(\Bun_{G}(\Fq)/\Xi,\Qlbar)\to \Qlbar$ n'est autre que  $h\mapsto \int_{\Bun_{G}(\Fq)/\Xi}     fh$. 
Alors on déduit des propriétés de ces  morphismes $\pi^{\mc F,\mc E}_{I,W}$ que pour tout  $I,W,x,\xi,(\gamma_{i})_{i\in I}$, 
on a 
$$S_{I,W,x,\xi,(\gamma_{i})_{i\in I}}(f)=
\s{\xi, (\sigma(\gamma_{i}))_{i\in I}.x}f.$$ 
Ceci termine la preuve de la \propref{prop-compat-langl-geom}. \cqfd

\section{Rapports avec les travaux antérieurs}\label{intro-previous-works}

      Les méthodes utilisées dans ce travail sont complètement différentes de celles fondées sur la formule des traces qui ont été développées notamment 
  par Drinfeld 
   \cite{drinfeld78,Dr1,drinfeld-proof-peterson,drinfeld-compact}, Laumon, Rapoport et Stuhler \cite{laumon-rapoport-stuhler}, 
  Laumon \cite{laumon-drinfeld-modular,laumon-cetraro},    Laurent  Lafforgue \cite{laurent-asterisque,laurent-jams,laurent-inventiones,laurent-tata}, 
 Ngô Bao Châu \cite{ngo-jacquet-ye-ulm,ngo-modif-sym},  Eike Lau \cite{eike-lau,eike-lau-duke}, Ngo Dac Tuan \cite{ngo-dac-ast,ngo-dac-09,ngo-dac-11},  Ngô Bao Châu et  Ngo Dac Tuan  \cite{ngo-ngo-elliptique}, Kazhdan et  Varshavsky~\cite{kvar,var-SANT} et 
 Badulescu et Roche \cite{badulescu}. 
  
      Cependant  l'action sur la cohomologie des groupes de permutations des pattes des chtoucas apparaît déjà dans 
les travaux de  Ngô Bao Châu,  Ngo Dac Tuan et Eike Lau que nous venons de citer. Ces actions des groupes de permutations jouent par ailleurs un rôle essentiel  dans le programme de Langlands géométrique, et notamment dans la preuve par Gaitsgory de la conjecture d'annulation \cite{ga-vanishing}. 
D'autre part
 la coalescence des pattes  apparaît dans la thèse de Eike Lau \cite{eike-lau} et elle est aussi  utilisée dans le preprint   \cite{brav-var} de Braverman et Varshavsky  
 (afin de montrer la non nullité de certains des morphismes \eqref{trace-brav-var-piIW}).  
 L'article \cite{var} de Varshavsky sur  les champs  
 de $G$-chtoucas et le preprint  très éclairant  \cite{brav-var} de Braverman et  Varshavsky,   ont été fondamentaux pour nous.  Enfin on a mentionné dans le paragraphe précédent le lien avec le   corollaire  4.5.5 de 
 \cite{dennis-laumon} (qui généralise d'ailleurs la conjecture d'annulation).

  Dans la suite de ce paragraphe (qui ne  fournit aucun résultat nouveau et peut 
   être sautée par le lecteur)
on explique 
les quelques arguments nécessaires
  pour donner une nouvelle preuve de l'ingrédient de récurrence de 
  \cite{laurent-inventiones} à l'aide du \thmref{intro-thm-ppal}, en faisant attention à ne pas utiliser les résultats connus comme conséquences de  \cite{laurent-inventiones} pour éviter toute circularité. Bien que nous ayons en vue le cas de $GL_{r}$ il est plus naturel d'énoncer le lemme suivant pour $G$ arbitraire. 
      
                \begin{lem}\label{lem-sigma-dans-coho} 
 Soit $\sigma$ 
 apparaissant dans la décomposition 
      \eqref{intro1-dec-canonique} (c'est-à-dire tel que $\mf H_{\sigma}\neq 0$). 
    Soit  $V$ une  représentation irréductible de $\wh G$ et 
     $V_{\sigma}=\oplus_{\tau} \tau \otimes \mf V_{\tau}$ la  décomposition de la   représentation semi-simple $V_{\sigma}$ indexée par les 
     classes d'isomorphisme  de  représentations  irréductibles $\tau $ de $\pi_{1}(\eta, \ov\eta)$. Alors si  $\mf V_{\tau}\neq 0$, $\tau\boxtimes \tau ^{*}$ apparaît comme un sous-quotient de la   représentation 
     \begin{gather}\label{rep-12-VV*}  H_{\{1,2\}, V\boxtimes V^{*}}=    \Big(\varinjlim _{\mu}\restr{\mc H _{N, \{1,2\}, V\boxtimes V^{*}}^{\leq\mu}}{\ov{\eta^{\{1,2\}}}}\Big)^{\mr{Hf}}\end{gather}  de 
  $(\pi_{1}(\eta,\ov\eta))^{2}$. 
        De plus $\tau$ est $\iota$-pure pour tout isomorphisme $\iota:\Qlbar \isom \C$.      \end{lem}

\begin{rem} Bien sûr la dernière assertion n'est pas un résultat nouveau  car d'après le théorème VII.6 de   \cite{laurent-inventiones}   toute représentation  irréductible (définie sur une extension finie de $\Ql$ et continue) de $\pi_{1}(X\sm N, \ov\eta)$ est $\iota$-pure pour tout $\iota$. 
\end{rem}

   \noindent{\bf Démonstration. } Quitte à augmenter $E$ on suppose $\sigma$ et $\mf H_{\sigma}$ définis sur $E$. Soit $h\neq 0$ dans le sous-espace de $\mf H_{\sigma}$ sur lequel $\mc B$ agit par le caractère $\nu$ associé à $\sigma$ par \eqref{relation-fonda} (on sait que ce sous-espace est non nul bien que l'on ne sache pas si $\mc B$ est réduite). 
  Soit  $\check h\in   C_{c}^{\rm{cusp}}(G(F)\backslash G(\mb A)/K_{N}\Xi,E)$ tel que \begin{gather}\label{int-check-h-h}\int_{G(F)\backslash G(\mb A)/K_{N}\Xi} \check h \,  h=1.\end{gather} 
  On note 
  \begin{itemize}
  \item $f$ l'élément de \eqref{rep-12-VV*} image de $h$ par la composée 
 $$C_{c}^{\rm{cusp}}(G(F)\backslash G(\mb A)/K_{N}\Xi,E)=
  H_{\{0\},\mbf  1}\xrightarrow{\mc H(\delta_{V})}
 H_{\{0\},V\otimes V^{*}}\isor{\chi_{\zeta_{\{1,2\}}}^{-1}} 
 H_{\{1,2\}, V\boxtimes V^{*}},$$
  \item $\check f$ la forme linéaire   sur  \eqref{rep-12-VV*}
  égale à la composée de  
 $$H_{\{1,2\}, V\boxtimes V^{*}} \isor{\chi_{\zeta_{\{1,2\}}}} H_{\{0\},V\otimes V^{*}}  \xrightarrow{\mc H(\on{ev}_{V})} 
  H_{\{0\},\mbf  1} =C_{c}^{\rm{cusp}}(G(F)\backslash G(\mb A)/K_{N}\Xi,E)$$
     et de la forme linéaire
     $$C_{c}^{\rm{cusp}}(G(F)\backslash G(\mb A)/K_{N}\Xi,E)\to E, \ \ g\mapsto \int_{G(F)\backslash G(\mb A)/K_{N}\Xi} \check h g.$$
  \end{itemize}
        Alors $f$  et $\check f$ 
      sont invariants par l'action diagonale de $\pi_{1}(\eta,\ov\eta)$. 
    Pour tout  $(\gamma, \gamma')\in (\pi_{1}(\eta,\ov\eta))^{2}$ on a  
        \begin{gather*}\s{\check f , (\gamma, \gamma') \cdot f}=
       \int_{G(F)\backslash G(\mb A)/K_{N}\Xi} \check h S_{\{1,2\}, V\boxtimes V^{*},\delta_{V},\on{ev}_{V},(\gamma,\gamma')}(h)\\
    =\nu(S_{\{1,2\}, V\boxtimes V^{*},\delta_{V},\on{ev}_{V},(\gamma,\gamma')})       =
        \chi_{V}(\sigma(\gamma\gamma'^{-1}))=  \chi_{V_{\sigma}}(\gamma\gamma'^{-1}),\end{gather*}  où 
   \begin{itemize}
   \item  la première égalité  vient de la définition des opérateurs d'excursion donnée dans \eqref{excursion-def-intro}, 
      \item la deuxième égalité résulte de   l'hypothèse que $h$ est vecteur propre pour $\mc B$ relativement au caractère $\nu$, et de \eqref{int-check-h-h}, 
      \item la troisième égalité vient du fait que $\nu$ est associé à $\sigma$ par \eqref{relation-fonda}. 
        \end{itemize}
       Le   quotient de la   représentation de $(\pi_{1}(\eta,\ov\eta))^{2}$ engendrée par  $f$ par la plus grande sous-représentation sur laquelle  $\check f$ s'annule  est alors isomorphe à la sous-représentation  engendrée par 
       $\chi_{V_{\sigma}}$  dans  
     $C( \pi_{1}(\eta,\ov\eta), E)$ muni  de l'action  par translations à gauche et à droite de $(\pi_{1}(\eta,\ov\eta))^{2}$. D'après \cite{Bki-A8} chapitre 20.5 théorème 1, cette dernière représentation est isomorphe à 
   $\oplus_{\tau, \mf V_{\tau}\neq 0} \tau\boxtimes \tau ^{*}$.   
    On a donc montré que  si   $\mf V_{\tau}\neq 0$, $\tau\boxtimes \tau ^{*}$ est un  quotient d'une  sous-représentation de  \eqref{rep-12-VV*}.      
     On en déduit maintenant que $\tau$ est $\iota$-pure.  Comme $\tau\boxtimes \tau^{*}$ est un sous-quotient de \eqref{rep-12-VV*}, 
          il résulte de Weil II \cite{weil2} que $\tau\boxtimes \tau ^{*}$ est 
  $\iota$-pure de poids $\leq 0$ comme représentation de $\pi_{1}((X\sm N)^{2}, \Delta(\ov\eta))$. 
 Donc pour presque toute place $v$ les valeurs propres de $\tau(\Frob_{v})$ ont  des $\iota$-poids égaux, et déterminés par le $\iota$-poids de $\det(\tau)$. 
     \cqfd
       
    A partir de maintenant on prend  $G=GL_{r}$.  
                                         On rappelle que dans \cite{laurent-inventiones} la correspondance de Langlands est obtenue par récurrence sur $r$,  à l'aide du   ``principe de récurrence'' de Deligne, qui combine 
                                    \begin{itemize}
                                  \item
    les équations fonctionnelles
de fonctions $L$ de Grothendieck \cite{sga5}, 
\item la formule du produit de Laumon \cite{laumon-produit}, 
\item les théorèmes de multiplicité $1$ \cite{piat71,shalika}  
et les  théorèmes réciproques de Hecke, Weil,  Piatetski-Shapiro et Cogdell \cite{inverse-thm}.  
\end{itemize}
    La récurrence est expliquée dans le paragraphe 6.1 et l'appendice B de 
    \cite{laurent-inventiones} et admet comme ingrédient
    \begin{gather}\label{hyp-laurent}
  \text{  l'hypothèse de la   proposition VI.11 (ii) de   \cite{laurent-inventiones}}
  \end{gather}
   à savoir qu'à toute représentation  automorphe cuspidale $\pi$ pour $GL_{r}$ de niveau $N$ on peut associer   $\sigma: \pi_{1}(X\sm N, \ov\eta)\to GL_{r}(\Qlbar) $ défini sur une extension finie de  $\Ql$, continu, 
    pur de poids $0$ et correspondant à $\pi$ au sens de Satake en toutes les places de $X\sm N$. Notre théorème
    \ref{intro-thm-ppal} fournit une nouvelle démonstration de  \eqref{hyp-laurent}, grâce au  lemme suivant. 
    
    \begin{lem} On prend $G=GL_{r}$. 
    On suppose la correspondance de Langlands connue pour  $GL_{r'}$ pour tout $r'<r$. 
         Alors tout $\sigma$ 
    apparaissant dans la décomposition \eqref{intro1-dec-canonique} est irréductible et pur  de poids  $0$. 
    \end{lem}
    \begin{rem} Ce lemme n'apporte pas de résultat nouveau car il 
    résulte {\it a posteriori} de  \cite{laurent-inventiones}. \end{rem}
 
 \noindent{\bf Démonstration (extraite de \cite{laurent-inventiones}). }   
 Soit $(H_{\pi},\pi)$ une représentation automorphe cuspidale telle que $(H_{\pi})^{K_{N}}$ soit non nul et apparaisse dans $\mf H_{\sigma}$. 
      On note  
   \begin{gather}\label{dec-sigma}\sigma=\oplus_{\tau} \tau \otimes \mf V_{\tau}\end{gather} 
   la décomposition de la représentation semi-simple $\sigma$ suivant les classes d'équivalence de  représentations irréductibles $\tau $ de $\pi_{1}(X\sm N, \ov\eta)$. 
    On suppose par l'absurde que $\sigma$ n'est pas irréductible. 
   Toute représentation   $\tau$ telle que $\mf V_{\tau}\neq 0$ est alors de rang  $r_{\tau}<r$ et,  par 
   hypothèse, il existe une   représentation automorphe cuspidale $\pi_{\tau}$  pour $GL_{r_{\tau}}$  associée à $\tau$ par la correspondance de Langlands pour $GL_{r_{\tau}}$.  On choisit un ensemble fini de places $S$ en dehors duquel  les représentations $\pi,\sigma$, et $\tau $, $\pi_{\tau}$ telles que $\mf V_{\tau}\neq 0$ (qui sont en nombre fini)  sont non ramifiées et se correspondent par l'isomorphisme de Satake. 
 D'après la méthode de Rankin-Selberg pour $GL_{r}\times GL_{r_{\tau}}$ (due à Jacquet, Piatetski-Shapiro, Shalika \cite{jacquet-shalika-euler,rankin-selberg}), 
$L(\check{\pi} \times \pi_{\tau},Z)$ est un polynôme en $Z$  
donc a fortiori 
 $L_{S}(\check{\pi} \times \pi_{\tau},Z)$ est un polynôme en $Z$ 
 (puisque les facteurs locaux ont des pôles mais jamais de zéros). 
  Donc 
$$L_{S}(\check{\pi}  \times \pi,Z)=L_{S}(\check{\pi} \times \sigma,Z)=\prod_{\tau }L_{S}(\check{\pi} \times \tau,Z) ^{\dim \mf V_{\tau}}=\prod_{\tau }L_{S}(\check{\pi} \times \pi_{\tau},Z)^{\dim \mf V_{\tau}}$$  n'a pas de pôle. 
 Pourtant d'après le théorème B.10 de \cite{laurent-inventiones} (dû à Jacquet, Shahidi, Shalika), 
 $L_{S}(\check{\pi}  \times \pi,Z)$ a un pôle en $Z=q^{-1}$, ce qui amène  une contradiction. 
     Donc on sait maintenant que $\sigma$ est irréductible. Pour tout $\iota$ on sait d'après le \lemref{lem-sigma-dans-coho} que $\sigma$ est $\iota$-pur et  la connaissance de son déterminant (par la théorie du corps de classe) implique alors que le $\iota$-poids est nul. Donc $\sigma$ est pur de poids $0$. \cqfd
  
  \begin{rem} On mentionne pour le lecteur que les conséquences déjà très importantes de la correspondance de Langlands pour $GL_{r}$, expliquées dans le chapitre VII de \cite{laurent-inventiones},  ont été encore  étendues dans des travaux récents de Deligne et Drinfeld 
  \cite{deligne-finitude,esnault-deligne,drinfeld-del-conj}
  sur les questions de rationalité 
  des traces des Frobenius et sur l'indépendance de $\ell$ pour les faisceaux $\ell$-adiques sur les  variétés lisses sur $\Fq$. 
\end{rem}


\begin{thebibliography}{hartl}

\bibitem[SGA4-2-VIII]{grothendieck-sga4-2-VIII}
A. Grothendieck. 
\newblock Foncteurs fibres, supports, étude cohomologique des morphismes finis, exposé VIII pages 366--412 dans: 
{\em Théorie des topos et cohomologie étale des schémas, }  SGA 4 (Tome 2),  Lecture Notes in Mathematics {\bf 270}, Springer, (1972)

\bibitem[SGA5]{sga5}
\newblock {Cohomologie $\ell$-adique et fonctions L}. 
Séminaire de Géométrie
Algébrique du Bois-Marie 1965--1966 (SGA 5). 
 Edité par Luc Illusie. Lecture Notes
in Mathematics {\bf 589} Springer  (1977)
 
\bibitem[BR17]{badulescu}
A. Badulescu and  Ph. Roche. 
\newblock Global Jacquet-Langlands correspondence for division algebras in characteristic $p$. 
\newblock 
{\em Int Math Res Notices}   {\bf 7}, 2172--2206  (2017)  



\bibitem[BMR05]{bmr}
 M. Bate,  B. Martin and   G. Röhrle. \newblock 
  A geometric approach to complete reducibility. \newblock  {\em  Invent. Math.} {\em 161}  (1) 177--218 (2005) 

\bibitem[BL95]{BL}
   A. Beauville and Y. Laszlo. 
   \newblock Un lemme de descente. 
   \newblock  {\em C. R. Acad. Sci. Paris Sér 
I Math. } {\bf 320}, 335--340 (1995)
  



\bibitem[BD99]{hitchin} A. Beilinson and V. Drinfeld. 
 \newblock
 Quantization of Hitchin's integrable system and Hecke eigensheaves (1999),   available at the address  \newline
  \url{http://math.uchicago.edu/~mitya/langlands.html} 

\bibitem[BD04]{chiral} A. Beilinson and V. Drinfeld. 
 \newblock {\em Chiral algebras. } 
 \newblock American Mathematical Society Colloquium Publications {\bf 51} (2004)

\bibitem[Bei06]{beilinson-heisenberg} 
A. Beilinson. 
  \newblock  Langlands parameters for Heisenberg modules. 
   \newblock {\em Studies in Lie theory} 51--60, Progr. Math. {\bf   243}, Birkhäuser  (2006)

  \bibitem[Bla94]{blasius}
D. Blasius. 
\newblock  On multiplicities for $SL(n)$. 
\newblock {\em Israel J. Math.} {\bf  88}, 237--251 (1994)

\bibitem[BR94]{blasius-rogawski-pspm}
D. Blasius  and J. Rogawski. 
\newblock Zeta functions of Shimura varieties. 
\newblock  Proc. Symp.
Pure Math. {\bf  55}, AMS, Part 2, 525--571 (1994)

\bibitem[BHKT16]{boeckle-harris...}
G.  B\"ockle, M. Harris, C. Khare, and J. Thorne. 
\newblock
$\wh G$-local systems on smooth projective curves are potentially automorphic. 
 \newblock 
Preprint, 	arXiv:1609.03491 (2016)

\bibitem[Bor79]{borel-corvallis}
A. Borel. 
\newblock Automorphic L-functions. 
\newblock {\em Automorphic forms, representations and L-functions (Corvallis)}, Part 2, Proc. Sympos. Pure Math., XXXIII, AMS,  27--61 (1979)


\bibitem[Bou12]{Bki-A8}
N. Bourbaki. 
\newblock Algèbre 8 (deuxième édition)
\newblock Springer, 2012

\bibitem[BT84]{bruhat-tits}
F. Bruhat et J. Tits. 
\newblock  Groupes réductifs sur un corps local. II. Schémas en groupes. Existence d'une
donnée radicielle valuée. \newblock
{\em Publ. Math. IHES} {\bf  60}, 197--376 (1984)


\bibitem[BG02]{brav-gaitsgory}
A. Braverman and D. Gaitsgory. 
\newblock 
Geometric Eisenstein series. 
\newblock  {\em Invent. Math. } {\bf 150} (2) 287--384  (2002) 


\bibitem[BV06]{brav-var}
A. Braverman and Y. Varshavsky. \newblock
 From automorphic sheaves to the cohomology of the moduli spaces of $F$-bundles. \newblock Unpublished preprint 
  (2006)
  
   
    
  \bibitem[BG11]{buzzard-gee}
 K. Buzzard and T. Gee. 
  \newblock The conjectural connections between automorphic representations and Galois representations. 
   \newblock    {\em Automorphic Forms and Galois Representations: Volume 1}, 
   London Math. Soc. Lect. Notes Series 414, 135--187 (2014)



 \bibitem[Car79]{cartier-satake}
  P. Cartier. 
   Representations of $p$-adic groups: a survey. 
  \newblock {\em   Automorphic forms, representations and L-functions},    Corvallis, Proc. Sympos. Pure Math. XXXIII, Amer. Math. Soc., 111--155 (1979)
     
      
\bibitem[CPS94]{inverse-thm}
J.W. Cogdell and  I.I. Piatetski-Shapiro. \newblock Converse theorems for $GL_{n}$. 
\newblock  {\em Publ. Math. IHES} {\bf  79}, 157--214 (1994)

\bibitem[Del71]{deligne-bki-mod}
P. Deligne. 
\newblock  Formes modulaires et representations $\ell$-adiques. \newblock 
 Séminaire Bourbaki 68/69 no. 355. {\em Lecture Notes in Mathematics} {\bf 179},  136-172. Springer (1971) 
 
\bibitem[Del80]{weil2}
P. Deligne.
\newblock  La conjecture de Weil. II. 
\newblock {\em Publ. Math. IHES} {\bf 52}, 137--252 (1980)
 
 \bibitem[Del12]{deligne-finitude}
 P. Deligne.   
 \newblock  Finitude de l'extension de $\mathbb Q$ engendrée par des traces de Frobenius, en caractéristique finie.  
 {\em  Mosc. Math. J.} {\bf 12} (3)  497--514 (2012)  


  \bibitem[Dri78]{drinfeld78}
V. G. Drinfeld. \newblock 
 Langlands' conjecture for GL(2) over functional fields.  \newblock  {\em Proceedings of the International Congress of Mathematicians (Helsinki, 1978)},  565--574, Acad. Sci. Fennica, Helsinki, 1980. 
 
 \bibitem[Dri87]{Dr1}
V. G. Drinfeld. \newblock   Moduli varieties of $F$-sheaves.  \newblock {\em Func. Anal. and Appl.}  {\bf 21},
107--122 (1987)

\bibitem[Dri88]{drinfeld-proof-peterson}
V. G. Drinfeld. 
\newblock   Proof of the Petersson conjecture for ${\rm GL}(2)$ over a global field of characteristic $p$. 
\newblock 
 {\em Funct. Anal. Appl.}  {\bf  22} (1), 28--43 (1988) 
 
\bibitem[Dri89]{drinfeld-compact}
V. G. Drinfeld. 
\newblock  Cohomology of compactified manifolds of modules of F-sheaves of rank 2.
\newblock 
 {\em
 J. Soviet Math.}  {\bf 46},  1789--1821 (1989)
 
\bibitem[Dri12]{drinfeld-del-conj}
V. G. Drinfeld.  \newblock 
On a conjecture of Deligne. 
\newblock  {\em Mosc. Math. J.}  {\bf 12} (3) 515--542 (2012). 
\newblock Voir \url{http://www.ihes.fr/jsp/site/Portal.jsp?document_id=3084&portlet_id=1140} pour un exposé à l'IHES.  

\bibitem[Dri15]{drinfeld-pro-completion}
V. G. Drinfeld.  \newblock 
On the pro-semisimple completion of the fundamental group of a smooth variety over a finite field. 
\newblock  Preprint,  	arXiv:1509.06059 (2015)


\bibitem[EK12]{esnault-deligne}
  H. Esnault and M. Kerz. \newblock 
 A finiteness theorem for Galois representations of function fields over finite fields (after Deligne). 
 \newblock  {\em Acta Math. Vietnam.} {\bf  37} (4) 531--562 (2012)   


   
  \bibitem[FC90]{faltings-chai}
G. Faltings and C.-L. Chai. 
 \newblock 
{\em Degeneration of abelian varieties. }
\newblock 
Ergebnisse der Mathematik  {\bf 22}, 
Springer
(1990)
 
  \bibitem[FL10]{finkelberg-lysenko}
 M.  Finkelberg and S. Lysenko.
  \newblock 
  Twisted geometric Satake equivalence. 
   \newblock 
 {\em J. Inst. Math. Jussieu} {\bf  9} (4) 719--739 (2010)
 
\bibitem[Gai01]{ga-iwahori}
D. Gaitsgory. 
\newblock   
Construction of central elements in the affine Hecke algebra via nearby cycles. 
\newblock   
{\em Invent. Math.} {\bf 144}, 253--280 (2001) 


 \bibitem[Gai04]{ga-vanishing}
D. Gaitsgory. 
\newblock    On a vanishing conjecture appearing in the geometric Langlands correspondence. 
\newblock   
{\em Ann. of Math.}  {\bf 160}  (2), 617--682 (2004)

\bibitem[Gai07]{ga-de-jong}
D. Gaitsgory. 
\newblock    On de Jong's conjecture. 
\newblock {\em Israel J. Math. } {\bf 157}, 155--191 (2007)

\bibitem[Gai15]{dennis-laumon}
D. Gaitsgory. 
\newblock
Outline of the proof of the geometric Langlands conjecture for $GL(2)$. 
\newblock {\em De la Géométrie Algébrique aux Formes Automorphes (II): Une collection d'articles en l'honneur du soixantième anniversaire de Gérard Laumon.}
Astérisque 370 (2015)

\bibitem[GL16]{dennis-sergey}
D. Gaitsgory and S. Lysenko. 
\newblock
Parameters and duality for the metaplectic geometric Langlands theory. 
\newblock Preprint,  arXiv:1608.00284 (2016)


\bibitem[GT05]{genestier-tilouine}
A. Genestier et  J. Tilouine. \newblock  
 Systèmes de Taylor-Wiles pour $\rm GSp\sb 4$.
 \newblock   {\em Formes automorphes. II. Le cas du groupe $\rm GSp(4)$. } Astérisque {\bf  302}, 177--290  (2005)

\bibitem[GL17]{genestier-lafforgue}
A. Genestier et  V. Lafforgue. 
 \newblock 
Chtoucas restreints pour les groupes réductifs et paramétrisation  de Langlands  locale. 
  \newblock  {\em Preprint. } (2017)
  
\bibitem[Gin95]{ginzburg}
 V. Ginzburg.
 \newblock   Perverse sheaves on a loop group and Langlands duality. 
 \newblock   Preprint, 
 	arXiv:alg-geom/9511007 (1995)


\bibitem[Gro98]{gross}
B.~Gross. 
\newblock On the Satake isomorphism.
\newblock {\em Galois Representations in Arithmetic Algebraic Geometry}, Cambridge University Press, 223--237, 1998, 
\newblock also available at the address   \newline 
\url{http://www.math.harvard.edu/~gross/preprints/sat.pdf}

\bibitem[Hei10]{heinloth-unif} J. Heinloth. 
\newblock  Uniformization of G-bundles. \newblock 
{\em Math. Ann.} {\bf  347} (3), 499--528 (2010) 
 

 
\bibitem[JS81]{jacquet-shalika-euler}
H. Jacquet and  J.A. Shalika. 
\newblock On Euler products and the classification of automorphic representations
I and II. 
\newblock 
 {\em American Journal of Mathematics} {\bf 103}, 499--558 and 777--815  (1981)

\bibitem[JPS83]{rankin-selberg}
H. Jacquet, I.I. Piatetski-Shapiro, J.A. Shalika. 
\newblock
Rankin-Selberg convolutions. 
\newblock {\em  American
Journal  of Math.} {\em   105}, 367--464  (1983)


\bibitem[KV13]{kvar}
D. Kazhdan and  Y. Varshavsky.
\newblock On the cohomology of the moduli spaces of $ F$ -bundles: stable cuspidal Deligne-Lusztig part, in preparation.
\newblock Work annonced in \cite{varshavsky-deligne}. 


\bibitem[Kos59]{kostant-betti}
B. Kostant. 
\newblock 
The principal three-dimensional subgroup and the Betti numbers of a complex simple Lie group. 
\newblock  {\em Amer. J. Math.} {\bf 81}, 973--1032 (1959) 


\bibitem[Kot84]{kottwitz1}
R. E. Kottwitz.
\newblock 
 Stable trace formula: cuspidal
tempered terms.
 \newblock {\em  Duke Math. J.} {\bf  51} (3), 611--650 (1984)
 
 \bibitem[Kot86]{kottwitz2}
R. E. Kottwitz.
\newblock   Stable trace formula: elliptic
singular terms. 
\newblock {\em  Math. Ann.} {\bf  275} (3), 365--399 (1986)

 \bibitem[Kot90]{kottwitz3}
R. E. Kottwitz.
\newblock Shimura varieties and $\lambda$-adic representations. 
\newblock {\em Automorphic forms, Shimura varieties, and L-functions (Ann Arbor)}, Vol. I, 161--209, Perspect. Math., 10, Academic Press  (1990) 


\bibitem[Laf97]{laurent-asterisque}
L. Lafforgue. \newblock  Chtoucas de Drinfeld et conjecture de Ramanujan-Petersson. \newblock 
{\em Ast\'erisque} {\bf 243} (1997) 

 
\bibitem[Laf98]{laurent-jams}
L. Lafforgue. \newblock  Une compactification des champs classifiant les chtoucas de Drinfeld, 
\newblock 
 {\em J. Amer. Math. Soc.}  {\bf 11}, 1001--1036 (1998)

\bibitem[Laf02a]{laurent-inventiones}
L. Lafforgue. 
\newblock  Chtoucas de Drinfeld et correspondance de Langlands.
\newblock {\em  Invent. Math.}  {\bf 147} (1), 1--241 (2002) 

\bibitem[Laf02b]{laurent-tata}
L. Lafforgue. \newblock 
Cours à l'Institut Tata sur les chtoucas de Drinfeld et la correspondance de Langlands.  
\newblock 
 {\em
 Prépublication IHÉS}, M/02/45 (2002)

  \bibitem[Laf12]{coh}
V. Lafforgue. 
 \newblock Chtoucas  pour les groupes réductifs et paramétrisation  de Langlands  globale.   \newblock  {\em Preprint  }  \url{http://arxiv.org/abs/1209.5352} (2012)

 \bibitem[Laf18]{texte-ICM}
V. Lafforgue. 
 \newblock Shtukas for reductive groups and Langlands correspondence for functions fields.   \newblock  {\em Preprint  } Text submitted for ICM 2018 (2018)



 \bibitem[Lan70]{langlands67}
 R. P. Langlands. 
 \newblock  Problems in the theory of automorphic forms. 
 \newblock {\em Lectures in modern analysis and applications, III},  18--61. 
 Lecture Notes in Math. {\bf 170} Springer  (1970) 

\bibitem[Lap99]{lapid}
E. Lapid. 
 \newblock Some results on multiplicities for $SL(n)$.  
  \newblock {\em Israel J. Math.} {\bf 112}, 157--186 (1999) 
  
  \bibitem[Lar94]{larsen1}
 M.  Larsen. 
 \newblock  On the conjugacy of element-conjugate homomorphisms.
 \newblock {\em  Israel J. Math.} {\bf 88},   253--277 (1994) 
  
 \bibitem[Lar96]{larsen2}
 M.  Larsen. 
 \newblock  On the conjugacy of element-conjugate homomorphisms. II. 
 \newblock {\em Quart. J. Math. Oxford Ser.} (2) {\bf 47} (185), 73--85   (1996)
  
 \bibitem[LO08]{laszlo-olsson}
  Y.   Laszlo and M. Olsson. 
 \newblock   The six operations for sheaves on Artin stacks. I. Finite coefficients. and II. Adic coefficients. 
  \newblock {\em 
 Publ. Math. Inst. Hautes Études Sci.}  {\bf 107}, 109--168 and 169--210 (2008)
  
  

 \bibitem[Lau87]{laumon-produit}
G. Laumon. 
\newblock  Transformation de Fourier, constantes d'équations fonctionnelles et conjecture de Weil. 
\newblock {\em Publ. Math. IHES} {\bf 65}, 131--210  (1987)

\bibitem[Lau96]{laumon-drinfeld-modular}
G. Laumon. 
\newblock 
{\em  Cohomology of Drinfeld modular varieties. Part I,II.} Cambridge Studies in Advanced Mathematics, {\bf  41} and {\bf  56}. Cambridge University Press (1997)

\bibitem[Lau97]{laumon-cetraro}
G. Laumon. 
\newblock 
  Drinfeld shtukas. Vector bundles on curves -- new directions (Cetraro, 1995). 
\newblock 
 {\em Lecture Notes in Math.} {\bf 1649}, 50--109,  Springer  (1997)


  
 \bibitem[Lau04]{eike-lau}
E.  Lau. 
\newblock  On generalised $\mc D$-shtukas. \newblock Dissertation, Rheinische Friedrich-Wilhelms-Universität Bonn, 2004. {\em Bonner Mathematische Schriften} {\bf 369}. Universität Bonn, Mathematisches Institut, Bonn (2004), available at the address  \newline
\url{http://www.math.uni-bielefeld.de/~lau/publ.html}

\bibitem[Lau07]{eike-lau-duke}
 E. Lau. 
 \newblock  On degenerations of $\mc D$-shtukas. 
 \newblock    {\em Duke Math. J.} {\bf 140} (2), 351--389 (2007)


\bibitem[LMB99]{laumon-moret-bailly} 
G. Laumon et L. Moret-Bailly. 
\newblock 
 {\em Champs algébriques.}  Ergebnisse der Math. {\bf 39}, Springer (1999) 


\bibitem[LRS93]{laumon-rapoport-stuhler} G.  Laumon, M. Rapoport, U. Stuhler. 
\newblock 
D-elliptic sheaves and the Langlands correspondence.
\newblock 
 {\em  Invent. Math.} {\bf  113} (2), 217--338 (1993) 
 
\bibitem[Lus82]{lusztig-satake}
 G. Lusztig. 
 \newblock Singularities, character formulas, and a $q$-analogue for weight multiplicities. 
 \newblock
Analyse et topologie sur les espaces singuliers, {\em  Astérisque} {\bf 102}, 
208--229 (1982) 


\bibitem[Lys06]{sergey-theta}
 S. Lysenko. 
 \newblock   Moduli of metaplectic bundles on curves and theta-sheaves. 
 \newblock  {\em Ann. Sci. École Norm. Sup.} {\bf 39} (3), 415--466  (2006)  

 \bibitem[Lys11]{sergey-theta-SO-Sp}
 S. Lysenko. 
 \newblock   Geometric theta-lifting for the dual pair $SO_{2m}$,$Sp_{2n}$. 
  \newblock  {\em Ann. Sci. École Norm. Sup.} {\bf 44} (3), 427--493 (2011)
  
   \bibitem[Lys14]{lysenko-red}
 S. Lysenko. 
 \newblock   Twisted geometric Satake equivalence: reductive case. 
  \newblock  {\em preprint},   	arXiv:1411.6782 [math.RT] (2014)
     

  
  
\bibitem[MV07]{mv} 
I. Mirkovic and K. Vilonen.  \newblock 
Geometric Langlands duality and representations of algebraic groups over commutative rings. \newblock 
{\em Annals of Math.} {\bf 166}, 95--143 (2007)



\bibitem[NBC99]{ngo-jacquet-ye-ulm}  
 Ngô Bao Châu.  \newblock 
 Faisceaux pervers, homomorphisme de changement de base et lemme fondamental de Jacquet et Ye.
\newblock 
 {\em Ann. Sci. École Norm. Sup. } {\bf 32} (5), 619--679 (1999) 

\bibitem[NBC06a]{ngo-modif-sym} 
 Ngô Bao Châu.  \newblock 
 D-chtoucas de Drinfeld à modifications symétriques et identité de changement de base. 
\newblock 
 {\em
  Ann. Sci. École Norm. Sup.} {\bf 39} (2), 197--243 (2006)
  



\bibitem[NN08]{ngo-ngo-elliptique} 
Ngô Bao Châu  et  Ngo Dac Tuan. 
\newblock 
 Comptage de $G$-chtoucas: la partie régulière elliptique.
\newblock 
 {\em J. Inst. Math. Jussieu}  {\bf 7} (1), 181--203 (2008)



\bibitem[NDT07]{ngo-dac-ast}
Ngo Dac Tuan. 
\newblock 
 {\em
 Compactification des champs of chtoucas et théorie géométrique des invariants. } \newblock 
Astérisque {\bf 313 } (2007)


\bibitem[NDT09]{ngo-dac-09}
Ngo Dac Tuan. 
\newblock  Sur le développement spectral de la formule des traces d'Arthur-Selberg pour les corps of fonctions. 
          \newblock 
   {\em  Bulletin de la Soc. Math. France}  {\bf 137},  545--586 (2009)



\bibitem[NDT11]{ngo-dac-11} 
Ngo Dac Tuan. 
\newblock  Formule des traces d'Arthur-Selberg pour les corps de fonctions II, preprint available at the address   \newline 
\url{http://www.math.univ-paris13.fr/~ngodac/traceformula2.pdf}


\bibitem[NQT11]{thang}
Nguyen Quoc Thang. 
\newblock On Galois cohomology and weak approximation of connected reductive groups over fields of positive characteristic. 
\newblock  {\em Proc. Japan Acad. Ser. A Math. Sci.} {\bf 87}  (10), 203--208 (2011)  

  
\bibitem[Pia75]{piat71}
I.I. Piatetski-Shapiro.
\newblock 
 Euler subgroups. 
 \newblock 
{\em Lie groups and their representations. } 
Proc. Summer School, Bolyai János Math. Soc., Budapest, 1971,   597--620. Halsted, New York, 1975.
 

\bibitem[Ric88]{richardson}
R. W. Richardson. 
\newblock 
 Conjugacy classes of $n$-tuples in Lie algebras and
algebraic groups. 
 \newblock 
 {\em  Duke Math. J.} {\bf 57}, 1--35  (1988)
       
       \bibitem[Ric14]{richarz}
 T.  Richarz. \newblock 
 A new approach to the geometric Satake equivalence
  	\newblock {\em Doc. Math. } {\bf 19}, 209--246  (2014)  
	  
  \bibitem[Sat63]{satake}
   I. Satake. \newblock  Theory of spherical functions on reductive algebraic groups over $p$-adic
fields.  \newblock {\em Publ.  Math. IHES} {\bf 18},  5--69 (1963)
  
     
  


\bibitem[Ser05]{bki-serre} 
J.-P. Serre.  \newblock 
Complète réductibilité.  \newblock
Séminaire Bourbaki. Vol. 2003/2004.
Astérisque No. 299 (2005)

   
  \bibitem[Sha74]{shalika}
 J.A. Shalika. 
\newblock  The multiplicity one theorem for $GL_{n}$. 
\newblock 
{\em Ann. of Math.} {\bf  100} 171--193 (1974)  

 
\bibitem[Tay91]{taylor}
R. Taylor. \newblock  Galois representations associated to Siegel modular forms of low weight. 
\newblock {\em  Duke Math.
J.} {\bf 63},  281--332 (1991)

\bibitem[Var04]{var}
Y. Varshavsky. \newblock
 Moduli spaces of principal $F$-bundles. 
 \newblock {\em  Selecta Math. (N.S.) }  {\bf 10} (1), 131--166 (2004)
  
  
\bibitem[Var05]{varshavsky-deligne}
Y. Varshavsky.  \newblock
A proof of a generalization of Deligne's conjecture. 
\newblock {\em Electron. Res. Announc. Amer. Math. Soc. } {\bf 11},  78--88
 (2005)


\bibitem[Var09]{var-SANT}
Y. Varshavsky. 
\newblock
Towards Langlands correspondence over function fields for split reductive groups. 
\newblock Symmetries in algebra and number theory (SANT), 139--147, Universitätsverlag Göttingen,  available at the address 
\newline \url{http://webdoc.sub.gwdg.de/univerlag/2009/SANT.pdf} 
(2009).
  
   \bibitem[Xue17]{these-cong}
Cong Xue. \newblock  
Cohomologie cuspidale des champs de chtoucas. 
\newblock Thèse de doctorat,  Université Paris-Saclay (2017)
   
 \bibitem[XZ17]{xiao-zhu}
Liang Xiao and  Xinwen Zhu. 
\newblock Cycles on Shimura varieties via geometric Satake. 
\newblock  Preprint.  	arXiv:1707.05700 (2017)
   
\bibitem[Zhu14]{zhu}
Xinwen Zhu.  
\newblock The geometric Satake correspondence for ramified groups (with an appendix by T. Richarz
and X. Zhu).  \newblock {\em Ann. Sci. École Norm. Sup.} {\bf 47} (6),  (2014)

   
\end{thebibliography}
\end{document}